\documentclass{amsart}
\usepackage{amsthm,amsmath,amsfonts}
\usepackage{mathrsfs}
\usepackage{euscript}
\usepackage[colorlinks=true,linktocpage=true]{hyperref}
\usepackage{tocvsec2}
\usepackage{accents}
\usepackage{bbm,dsfont}
\usepackage{hyperref}
\hypersetup{urlcolor=blue, citecolor=blue, linkcolor=red}
\usepackage{framed}
\usepackage{color}
\usepackage{xcolor}
\numberwithin{equation}{section}
\numberwithin{table}{section}
\font\tenscrpt=eusm10
\font\sevenscrpt=eusm10 scaled 700
\font\fivescrpt=eusm10 scaled 500
\newfam\eusmfam
\textfont\eusmfam=\tenscrpt
\scriptfont\eusmfam=\sevenscrpt
\scriptscriptfont\eusmfam=\fivescrpt

\newtheorem{thm}{Theorem}[section]
\newtheorem{cor}{Corollary}[section]
\newtheorem{lem}{Lemma}[section]
\newtheorem{prop}{Proposition}[section]

\theoremstyle{definition}
\newtheorem{defn}{Definition}[section]

\newtheorem{rem}{Remark}[section]
\newtheorem{notn}{Notation}[section]
\newcommand{\thmref}[1]{Theorem~\ref{#1}}
\newcommand{\secref}[1]{Section~\ref{#1}}

\newcommand{\lemref}[1]{Lemma~\ref{#1}}
\newcommand{\coref}[1]{Corollary~\ref{#1}}

\def\qed{\quad\vcenter{\hrule\hbox{\vrule height.6em\kern.6em\vrule}\hrule}}
\newenvironment{pf}{{\bigskip\textit{\newline Proof.}\quad}}{$\qed$\bigskip\newline}
\newenvironment{pf*}[1]{{\bigskip\textit{\newline#1.}\quad}}{$\qed$\bigskip\newline}
\def\ds{\displaystyle}
\def\ts{\textstyle}

\def\psoxz{{K}^{\text{\tiny{\sc{BM}}}^d}_{s_{1};x}}

\def\K{{\mathbb K}}

\def\KBtx{{\K}^{\text{\tiny{\sc{BTBM}}}^d}_{t;x}}
\def\FKBtxi{\hat{\K}^{\text{\tiny{\sc{BTBM}}}^d}_{t;\xi}}

\def\Kbetatx{{\K}^{(\beta)}_{t;x}}

\def\Kbetatx{{\K}^{{(\beta,d)}}_{t;x}}

\def\FKbetatxi{\hat{\K}^{(\beta,d)}_{t;\xi}}
\def\LFKbetathxi{\tilde{\hat{\K}}^{(\beta,d)}_{\theta;\xi}}
\def\Kbetatxy{{\K}^{(\beta,d)}_{t;x,y}}
\def\Kbetatsxy{{\K}^{(\beta,d)}_{t-s;x,y}}
\def\Kbetatrxy{{\K}^{(\beta,d)}_{t-r;x,y}}
\def\Kbetasrxy{{\K}^{(\beta,d)}_{s-r;x,y}}
\def\FKbetasrxxi{\hat{\K}^{(\beta,d)}_{s-r;x,\xi}}
\def\Kbetatrxz{{\K}^{(\beta,d)}_{t-r;x,z}}
\def\FKbetatrxxi{\hat{\K}^{(\beta,d)}_{t-r;x,\xi}}
\def\Kbetatryz{{\K}^{(\beta,d)}_{t-r;y,z}}
\def\FKbetatryxi{\hat{\K}^{(\beta,d)}_{t-r;y,\xi}}

\def\Kbetatrposxy{{\K}^{(\beta,d)}_{(t-r)_{+};x,y}}
\def\Kbetasrposxy{{\K}^{(\beta,d)}_{(s-r)_{+};x,y}}
\def\FKbetatrposxxi{\hat{\K}^{(\beta,d)}_{(t-r)_{+};x,\xi}}
\def\FKbetasrposxxi{\hat{\K}^{(\beta,d)}_{(s-r)_{+};x,\xi}}
\def\Kbetamrposxy{{\K}^{(\beta,d)}_{(-r)_{+};x,y}}

\def\psx{{K}^{\text{\tiny{\sc{BM}}}^d}_{s;x}}

\def\psoxz{{K}^{\text{\tiny{\sc{BM}}}^d}_{s_{1};x}}

\def\psoxy{{K}^{\text{\tiny{\sc{BM}}}^d}_{s_{1};x,y}}
\def\puoxy{{K}^{\text{\tiny{\sc{BM}}}^d}_{u_{1};x,y}}

\def\ptsz{{K}^{\text{\tiny{\sc{BM}}}}_{t;s}}

\def\ptzsk{{K}^{\text{\tiny{\sc{BM}}}}_{t;0,\frac{s_{k}}{\sqrt2}}}
\def\ptrzsk{{K}^{\text{\tiny{\sc{BM}}}}_{t-r;0,\frac{s_{k}}{\sqrt2}}}

\def\psrzuk{{K}^{\text{\tiny{\sc{BM}}}}_{s-r;0,\frac{u_{k}}{\sqrt2}}}

\def\ptzskiimo{{K}^{\text{\tiny{\sc{BM}}}}_{s_{k-i};0,\frac{s_{k-i-1}}{\sqrt2}}}

\def\ptzukiimo{{K}^{\text{\tiny{\sc{BM}}}}_{u_{k-i};0,\frac{u_{k-i-1}}{\sqrt2}}}

\def\peptsz{{K}^{\text{\tiny{\sc{BM}}}}_{\vep t;s}}

\def\KItzs{{K}^{{\Lambda_{\beta}}}_{t;s}}
\def\LKIlzs{\tilde{K}^{{\Lambda_{\beta}}}_{\theta;s}}

\def\KIhalftzxsqrt2{K^{\Lambda_{\frac12}}_{t;0,\sqrt2 x}}
\def\KabBmtzxsqrt2{K^{{|B|}}_{t;0,\frac{x}{\sqrt{2}}}}

\def\K{{\mathbb K}}

\def\KKSepthtx{{\K}^{\text{\tiny{\sc{LKS}}}^d_{\vep,\vth}}_{t;x}}

\def\KKSepthtxy{{\K}^{\text{\tiny{\sc{LKS}}}^d_{\vep,\vth}}_{t;x,y}}

\def\KKSepthtsxy{{\K}^{\text{\tiny{\sc{LKS}}}^d_{\vep,\vth}}_{t-s;x,y}}
\def\KKSepthtrxy{{\K}^{\text{\tiny{\sc{LKS}}}^d_{\vep,\vth}}_{t-r;x,y}}
\def\FKKSepthtrxxi{\hat{\K}^{\text{\tiny{\sc{LKS}}}^d_{\vep,\vth}}_{t-r;x,\xi}}
\def\KKSepthtrxz{{\K}^{\text{\tiny{\sc{LKS}}}^d_{\vep,\vth}}_{t-r;x,z}}
\def\KKSepthtryz{{\K}^{\text{\tiny{\sc{LKS}}}^d_{\vep,\vth}}_{t-r;y,z}}
\def\FKKSepthtryxi{\hat{\K}^{\text{\tiny{\sc{LKS}}}^d_{\vep,\vth}}_{t-r;y,\xi}}
\def\KKSepthtrposxy{{\K}^{\text{\tiny{\sc{LKS}}}^d_{\vep,\vth}}_{(t-r)_{+};x,y}}
\def\FKKSepthtrposxxi{\hat{\K}^{\text{\tiny{\sc{LKS}}}^d_{\vep,\vth}}_{(t-r)_{+};x,\xi}}
\def\KKSepthsrposxy{{\K}^{\text{\tiny{\sc{LKS}}}^d_{\vep,\vth}}_{(s-r)_{+};x,y}}
\def\FKKSepthsrposxxi{\hat{\K}^{\text{\tiny{\sc{LKS}}}^d_{\vep,\vth}}_{(s-r)_{+};x,\xi}}
\def\KKSepthmrposxy{{\K}^{\text{\tiny{\sc{LKS}}}^d_{\vep,\vth}}_{(-r)_{+};x,y}}
\def\KKSepthsrxy{{\K}^{\text{\tiny{\sc{LKS}}}^d_{\vep,\vth}}_{s-r;x,y}}
\def\FKKSepthsrxxi{\hat{\K}^{\text{\tiny{\sc{LKS}}}^d_{\vep,\vth}}_{s-r;x,\xi}}

\def\KKSepthtpry{{\K}^{\text{\tiny{\sc{LKS}}}^d_{\vep,\vth}}_{t+r;y}}
\def\FKKSepthtprxi{\hat{\K}^{\text{\tiny{\sc{LKS}}}^d_{\vep,\vth}}_{t+r;\xi}}

\def\FKKSepthtxi{\hat{\K}^{\text{\tiny{\sc{LKS}}}^d_{\vep,\vth}}_{ t;\xi}}

\def\unx{u_0(x)}

\def\uny{u_0(y)}
\def\un{u_0}

\def\ILb{\Lambda_{\beta}}

\def\sW{\mathscr W}

\def\P{{\mathbb P}}

\def\E{{\mathbb E}}
\def\N{{\mathbb N}}

\def\Rd{{\mathbb R}^{d}}

\def\Rpop{\mathring{\R}_{+}}
\def\R{\mathbb R}

\def\B{{\mathbb B}}
\def\S{\mathbb S}
\def\T{\mathbb T}

\def\Rp{{\R}_+}

\def\sF{{\mathscr F}}

\def\sFt{{\mathscr F}_t}

\def\OFFtP{(\Omega,\sF,\{\sFt\},\P)}

\def\C{\mathrm C}
\def\H{\mathrm H}

\def\sL{\mathscr L}
\def\eqL{\overset{{\sL}}{=}}

\def\D{\Delta}
\def\lap{\Delta}

\def\df#1#2{\ds{\frac{#1}{#2}}}
\def\tf#1#2{\ts{\frac{#1}{#2}}}
\def\lbl#1{\label{#1}}
\def\intrd{\int_{\Rd}}
\def\intrdzt{\int_{\Rd}\int_0^t}

\def\pa{\partial}

\def\lang{\left<}
\def\rang{\right>}
\def\lab{\left|}
\def\rab{\right|}
\def\lpa{\left(}
\def\rpa{\right)}

\def\lbk{\left[}
\def\rbk{\right]}
\def\lbr{\left\{}
\def\rbr{\right\}}
\def\bdf{\begin{defn}}
\def\edf{\end{defn}}
\def\bcr{\begin{cor}}
\def\ecr{\end{cor}}
\def\bnt{\begin{notn}}
\def\ent{\end{notn}}
\def\brm{\begin{rem}}
\def\erm{\end{rem}}
\def\blm{\begin{lem}}
\def\elm{\end{lem}}
\def\bpf{\begin{pf}}
\def\bpfs{\begin{pf*}}
\def\epf{\end{pf}}
\def\epfs{\end{pf*}}

\def\beq{\begin{equation}}
\def\beqs{\begin{equation*}}
\def\eeq{\end{equation}}
\def\eeqs{\end{equation*}}
\def\bsp{\begin{split}}
\def\esp{\end{split}}
\def\bc{\begin{cases}}
\def\ec{\end{cases}}
\def\bt{\begin{tabular}}
\def\et{\end{tabular}}

\def\bthm{\begin{thm}}
\def\ethm{\end{thm}}
\def\bpr{\begin{prop}}
\def\epr{\end{prop}}
\def\bfr{\begin{framed}}
\def\efr{\end{framed}}
\def\bsh{\begin{shaded}}
\def\esh{\end{shaded}}
\def\bcm{\iffalse}
\def\babs{\begin{abstract}}
\def\eabs{\end{abstract}}
\def\bit{\begin{itemize}}
\def\eit{\end{itemize}}
\def\ben{\begin{enumerate}}
\def\rencomalp{\renewcommand{\labelenumi}{(\alph{enumi})}}
\def\rencomrom{\renewcommand{\labelenumi}{(\roman{enumi})}}
\def\een{\end{enumerate}}
\def\babs{\begin{abstract}}
\def\eabs{\end{abstract}}

\def\ig{\iffalse}

\def\e{\mathrm{e}}
\def\p{\partial}
\def\tsig{\sigma^{(d)}_{\mbox{\tiny time}}}

\def\tint{I_{\mbox{\tiny time}}}
\def\ssig{\sigma^{(d)}_{\mbox{\tiny space}}}

\def\sint{I^{(d)}_{\mbox{\tiny space}}}
\def\sintth{I^{(3)}_{\mbox{\tiny space}}}
\def\sintt{I^{(2)}_{\mbox{\tiny space}}}
\def\sinto{I^{(1)}_{\mbox{\tiny space}}}

\def\jsig{\sigma^{(d)}_{\mbox{\tiny joint}}}

\def\var{\mathrm{Var}}
\def\sB{\mathscr{B}}

\def\i{\mathbf i}
\def\vep{\varepsilon}

\def\vth{\vartheta}

\def\pa{\partial}
\def\pat{\pa_{t}}
\def\I{\mathbb{I}}
\def\ind{\mathbbm{1}}

\def\pax{\pa_{x}}
\def\Fr{\mathcal{F}}
\def\La{\mathcal{L}}
\begin{document}
\title[Continuity moduli of time-fractional SPIDEs and L-KS SPDE{\scriptsize s}]
{L-Kuramoto-Sivashinsky SPDE{\scriptsize s} vs.~time-fractional SPIDE{\scriptsize s}: 
exact continuity and gradient moduli, 1/2-derivative criticality, and laws
}
\author{Hassan Allouba}
\address{Department of Mathematical Sciences, Kent State University, Kent, Ohio 44242
\\email: allouba@math.kent.edu}
\author{Yimin Xiao}
\address{Department of Statistics, Michigan State University, East Lansing, Michigan 44242
\\email: xiao@stt.msu.edu}

\thanks{Research of H. Allouba is partially supported by a Kent State Academic Year Research award.}\thanks{Research of Y. Xiao is partially
supported by NSF grants DMS-1307470 and DMS-1309856.}
\begin{abstract}
We establish exact, dimension-dependent, spatio-temporal, uniform and local moduli of continuity for (1) the fourth order L-Kuramoto-Sivashinsky (L-KS) SPDEs and for (2) the time-fractional stochastic partial integro-differential equations (SPIDEs), driven by space-time white noise in one-to-three dimensional space. Both classes were introduced---with Brownian-time-type kernel formulations---by Allouba in a series of articles starting in 2006, where he presented class (2) in its rigorous stochastic integral equations form. He proved existence, uniqueness, and sharp spatio-temporal H\"older regularity for the above two classes of equations in $d=1,2,3$. We show that both classes are $(1/2)^-$ H\"older continuously differentiable in space when $d=1$, and we give the exact uniform and local moduli of continuity for the gradient in both cases. This is unprecedented for SPDEs driven by space-time white noise. Our results on exact moduli show that the half-derivative SPIDE is a critical case. It signals the onset of rougher modulus regularity in space than both time-fractional SPIDEs with time-derivatives of order $<1/2$ and L-KS SPDEs. This is despite the fact that they all have identical spatial H\"older regularity, as shown earlier by Allouba. Moreover, we show that the temporal laws governing (1) and (2) are fundamentally different. We relate L-KS SPDEs to the Houdr\'e-Villa bifractional Brownian motion, yielding a Chung-type law of the iterated logarithm for these SPDEs. We use the underlying explicit kernels and spectral/harmonic analysis to prove our results. On one hand, this work builds on the recent works on delicate sample path properties of Gaussian random fields. On the other hand, it builds on and complements Allouba's earlier works on (1) and (2). Similar regularity results hold for Allen-Cahn nonlinear members of (1) and (2) on compacts via change of measure.
\end{abstract}
\maketitle
\newpage
\tableofcontents
\section{Introduction, statement of results, and preliminaries}\lbl{intro}
\subsection{Two sides of the Brownian-time coin}
We delve into delicate regularity properties of paths of fourth order pattern formation stochastic PDEs (SPDEs) and time-fractional slow diffusion stochastic partial integro-differential equations (SPIDEs).  The fundamental kernels associated with the deterministic versions of these two different classes are both built on the Brownian-time processes in \cite{Abtp1,Abtp2,Aks}\footnote{See also the multi-time-parameter interacting versions of these equations in \cite{Abtbs,AN}} and extensions thereof.  We thus think of these two classes of equations as ``two sides of the Brownian-time coin''.  It is therefore often useful and efficient to study both simultaneously and compare and contrast their various properties.  In this article, we unveil a rather detailed set of results giving the exact dimension-dependent uniform and local modulus of continuity, in time and space, for two important classes of stochastic equations:
\ben
\item the fourth order L-Kuramoto-Sivashinsky (L-KS) SPDEs connected to pattern formation phenomena accompanying the appearance of turbulence  (see \cite{Alksspde,Abtpspde,Aks} for the L-KS class and for its connection to many classical and new examples of pattern formation deterministic and stochastic PDEs, and see \cite{DuanWei14,T} for classical examples of stochastic and deterministic pattern formation PDEs); and
\item time-fractional SPIDEs connected to slow diffusion or diffusion in material with memory (see \cite{Abtp1,Abtp2,DOOToa14,GOP15,LunSin,Main1,MBS,MeerSik} for connected PDEs in the deterministic setting and see \cite{Atfhosie,Abtbmsie,Abtpspde} for the associated stochastic integral equations (SIEs), followed later by the articles \cite{LeC,ZCHKmPKm,MN}, in the stochastic setting).
\een
We also characterize the temporal laws for these two classes of equations.  More specifically, we prove our results on the exact uniform and local moduli of continuity for the canonical equations\footnote{In addition to equations \eqref{nlks} and \eqref{spide}, the general L-KS SPDEs and time-fractional SPIDEs classes include many nonlinear equations (both well known as well as new).  We refer the reader to \thmref{SHcom} and Sections \ref{sec:fvq} and \ref{rigdisc} below for more on that.  The constants in \eqref{nlks} and \eqref{spide} can easily be changed by scaling.  We will alternate freely between the notations $\pa_{u}f(u,v)$ and $\pa f(u,v)/\pa u$.}
\begin{equation} \label{nlks}
 \begin{cases} \displaystyle\frac{\partial U}{\partial t}=
-\tfrac\vep8\lpa\lap+2\vth\rpa^{2}U+\frac{\partial^{d+1} W}{\partial t\partial x}, & (t,x)\in\Rpop\times\Rd;
\cr U(0,x)=\unx, & x\in\Rd,
\end{cases}
\end{equation}
and
\beq\lbl{spide}
\begin{cases} \displaystyle
{}^{\mbox{\tiny C}}\pa_{t}^{\beta}U_{\beta}=\tfrac12\D U_{\beta}+I_{t}^{1-\beta}
\lbk \frac{\partial^{d+1} W}{\partial t\partial x}\rbk,& (t,x)\in\Rpop\times\Rd;
\cr U_{\beta}(0,x)=\unx, & x\in\Rd,
\end{cases}
\eeq	
where $\Rpop=(0,\infty)$; $(\vep,\vth)\in\Rpop\times\R$ is a pair of parameters; 
$\beta\in(0,1/2]$; the noise term ${\partial^{d+1} W}/{\partial t\partial x}$ is 
the space-time white noise corresponding to the real-valued Brownian 
sheet\footnote{As in Walsh \cite{W}, we treat space-time white noise as 
a continuous orthogonal martingale measure, and we denote it by $\sW$.} 
$W$ on $\Rp\times\Rd$, $d=1,2,3$;  the time fractional derivative of order 
$\beta$, ${}^{\mbox{\tiny C}}\pa_{t}^{\beta}$, is the Caputo fractional operator
\beq\lbl{caputo}
{}^{\mbox{\tiny C}}\pa_{t}^{\beta} f(t):=\bc\df{1}{\Gamma(1-\beta)}\int_{0}^{t}\df{f'(\tau)}{(t-\tau)^{\beta}}d\tau,& 
\mbox{ if } 0<\beta<1;\\
\df{d}{dt}f(t),& \mbox{ if } \beta=1,\ec
\eeq
and the time fractional integral of order $\alpha$, $I_{t}^{\alpha}$, is the 
Riemann-Liouville fractional integral of order $\alpha$:
\beq\lbl{RiemLiou}
I_{t}^{\alpha}f:=\df{1}{\Gamma(\alpha)}\int_{0}^{t}\df{f(s)}{(t-s)^{1-\alpha}}ds, 
\mbox{ for }t>0 \mbox{ and }\alpha>0,
\eeq
and $I_{t}^{0} =\mbox{Id}$, the identity operator.  The initial data $\un$ 
here is assumed Borel measurable, deterministic, and suitably regular.  
For convenience and for the sake of comparing \eqref{nlks} to \eqref{spide}, 
when $\beta\in\{1/2^{k};k\in\N\}$, we assume throughout this article 
that there is a $0<\gamma\le1$ such that\footnote{This is for convenience and we may relax these conditions \`a la those in \cite{AN}.}
\beq\lbl{init}
\bc
(a)\ \un\in\mathrm{C}_{b}^{2,\gamma}(\Rd;\R),&\mbox{ for }\eqref{nlks};\\
(b)\ \un\in\C_{b}^{2^{k+1}-2,\gamma}(\Rd;\R),&\mbox{ for }\eqref{spide},
2^{k-1}<\beta^{-1}\le2^{k}, k\in\N
\ec
\eeq
Of course, equations \eqref{nlks} and \eqref{spide} are formal (and nonrigorous) equations. Their rigorous formulations, which we work with in this article, are given in mild form as kernel stochastic integral equations (SIEs).  These SIEs were first introduced and treated by Allouba \cite{Abtpspde,Abtbmsie,Atfhosie,Alksspde}, with their genesis in   \cite{Abtp1,Abtp2,Aks}.   We give them below in \secref{rigdisc}, along with the relevant details.  
  
The results here build on the following works:  (1) Allouba 
\cite{Abtpspde,Abtbmsie,Atfhosie,Alksspde} who established the existence/uniqueness as well as sharp dimension-dependent $L^{p}$ and H\"older regularity of the linear and nonlinear noise versions of \eqref{nlks} and \eqref{spide} (he presented and treated the later in its stochastic integral equation form); and (2) Xiao \cite{X07, X09}; Meerschaert, Wang, and Xiao \cite{MWX}; 
Wu and Xiao \cite{WX06}; Xiao and Xue \cite{XX11}; who established several 
delicate analytic and geometric path properties of Gaussian processes 
and random fields (see also the related works in \cite{TuX07,TuX15,LX12,LX10}).

\subsection{Five questions}\lbl{sec:fvq}
In a series of articles \cite{Abtpspde,Abtbmsie,Atfhosie}, Allouba introduced 
and investigated the regularity of the rigorous kernel stochastic integral 
form\footnote{See \secref{rigdisc} below for the rigorous kernel or mild 
stochastic integral equation formulation for both L-KS SPDEs and time-fractional 
SPIDEs.} of the formal time-fractional SPIDEs in \eqref{spide} with diffusion 
coefficient $a$:
\beq\lbl{spidea}
\begin{cases} \displaystyle
{}^{\mbox{\tiny C}}\pa_{t}^{\beta}U_{\beta}=\tfrac12\D U_{\beta}+I_{t}^{1-\beta}\lbk a(U_{\beta})\frac{\partial^{d+1} W}{\partial t\partial x}\rbk, \quad & (t,x)\in\Rpop\times\Rd;
\cr U_{\beta}(0,x)=\unx, & x\in\Rd.
\end{cases}
\eeq	
He called these stochastic integral equations time-fractional and Brownian-time Brownian motion ($\beta=1/2$) SIEs.  Starting with the 2006 article \cite{Abtpspde}, he proved the existence  of a pathwise unique, continuous, and  $L^{p}$ bounded random field solution on $\{\Rp\times\Rd\}$, $d=1,2,3$, to the stochastic integral equation formulation of \eqref{spide} when $\beta=1/2$ (the Brownian-time process or Brownian-time Brownian motion (BTBM) SIE).  He proved in \cite{Abtpspde} that, in the case $a\equiv1$, the solution $U$ satisfies the $L^{p}$ bound
\beq\lbl{btbmsielpbound}
\sup_{x\in\Rd}\E|U(t,x)|^{2p}\le C\lbk1+t^{\tf{(4-d)p}{4}}\rbk; \quad t>0, d=1,2,3, \mbox{ and }p\ge1.
\eeq
He further proved in \cite{Abtbmsie,Atfhosie} that, under a Lipschitz assumption on the nonlinear $a$, there is a pathwise H\"older continuous solution $U_{\beta}$ to the SIE formulation of \eqref{spidea} such that, for any arbitrary $T>0$ and $\T=[0,T]$,
\beq\lbl{holregspide}
U_{\beta}\in\H^{\lpa{\tf{2\beta^{-1}-d}{4{\beta^{-1}}}}\rpa^{-},{\lpa\tf{4-d}{2}\wedge 1\rpa}^{-}}(\T\times\Rd;\R),
\eeq
for every $d=1,2,3$, $\beta^{-1}\in\{2^{k};k\in\N\}$; where $\H^{\gamma_t^{-},\gamma_s^{-}}(\T\times\Rd;\R)$ is the space of real-valued locally H\"older functions on $\T\times\Rd$ whose time and space H\"older exponents are in $(0,\gamma_t)$ and $(0,\gamma_s)$, respectively.  He also proved in \cite{Abtbmsie,Atfhosie}, under just continuity (no Lipschitz assumption) and linear growth conditions on the nonlinear $a$, the existence of lattice limits solutions to the SIE corresponding to \eqref{spidea}, with the same H\"older regularity as in \eqref{holregspide}.   In \cite{Abtpspde,Alksspde}, motivated by \cite{Aks}, he introduced and gave the explicit kernel stochastic integral equation formulation for a large class of stochastic equations he called L-KS SPDEs.  This class includes stochasric versions of prominent nonlinear equations like the Swift-Hohenberg PDE, variants of the Kuramoto-Sivashinsky PDE, as well as many new ones (see \cite{Alksspde}). He established in \cite{Alksspde}, among other things, the existence of a pathwise unique solution $U$ to the nonlinear L-KS SPDE \eqref{nlks} with Lipschitz diffusion coefficient $a$:
\begin{equation} \label{nlksa}
 \begin{cases} \displaystyle\frac{\partial U}{\partial t}=
-\tfrac\vep8\lpa\lap+2\vth\rpa^{2}U+a(U)\frac{\partial^{d+1} W}{\partial t\partial x}, \quad & (t,x)\in\Rpop\times\Rd;
\cr U(0,x)=\unx, & x\in\Rd,
\end{cases}
\end{equation}
with the same H\"older regularity for $U$ as the $\beta=1/2$ case in \eqref{spidea} (BTBM SPIDE or SIE\footnote{\lbl{FN:btbmpdefo} We remind the reader that when $\beta=1/2$ the scaled BTBM kernel, which is the fundamental solution to the $a\equiv0$ version of \eqref{spidea}, is also the fundamental solution to the fourth order \emph{memoryful} PDE
\beq\lbl{btppdefo}
\bc
\pa_{t} u= \ds{\df{\Delta\un }{\sqrt{4\pi t}}+\df14\Delta^2 u};&(t,x)\in(0,\infty)\times\Rd, \cr
u(0,x) = \unx;&x\in\Rd.
\ec
\eeq
first obtained in \cite{Abtp1,Abtp2}.  Of course, these sharp H\"older exponents for \eqref{nlksa} and \eqref{spidea} play
a crucial role in our exact moduli of continuity for both the L-KS SPDE and the $\beta$ time fractional SPIDEs, as is clear
from \thmref{tempmodlks}--\thmref{spatmodtfhalf} below.}) obtained by plugging $\beta=1/2$ in \eqref{holregspide}.  In addition,
the articles \cite{Abtbmsie,Atfhosie,Alksspde} are the first to obtain solutions to space-time white noise driven equations
that are smoother in time or space---twice as smooth in space in $d=1,2$, as is clear from \eqref{holregspide}---than the
Brownian sheet $W$ corresponding to the driving white noise.  Moreover, the kernels in these time-fractional SIEs, when
$\beta\in\{1/2^{k};k\in\N\}$, are fundamental solutions to higher order PDEs (where $\beta^{-1}$ is the order of the Laplacian as detailed in \cite{Atfhosie}).  Thus, the regularity results in \cite{Atfhosie}, given in \eqref{holregspide}, mean that the maximum
integer number of dimensions for the existence of random field solutions for space-time white noise driven equations is $3$,
no matter how high the Laplacian order is.  They also mean that the solutions for such equations are spatially
$\gamma_{s}$-H\"older for all $\gamma_{s}\in(0,1)$ (nearly locally Lipschitz) in dimensions $d=1,2$ and $\gamma_{s}\in(0,1/2)$ (nearly locally H\"older $1/2$) in $d=3$.
As observed in \cite{Atfhosie}, when $\beta\in\{1/2^{k};k\in\N\}$, letting $\beta\searrow0$ (the order of the Laplacian $\beta^{-1}\nearrow\infty$) \emph{does not} increase the spatial H\"older regularity and the extra H\"older regularizing force is manifested entirely temporally.  

These results in \cite{Alksspde}--\cite{Abtbmsie} naturally lead to the following list of motivating questions:
\bit
\item[(Q1)] Consider the L-KS SPDE and time-fractional SPIDEs in spatial dimension $d=1$.
\ben\rencomalp
\item Are the solutions to \eqref{nlksa} and \eqref{spidea} \emph{actually} spatially locally Lipschitz (not just nearly locally Lipschitz as in $d=2$)?  This would be unprecedented in SPDEs driven by space-time white noise, and is suggested by the sharp $L^{2}$ upper bounds on the kernels spatial differences in Lemma 2.4 in \cite{Abtbmsie,Atfhosie} and Lemma 3.3 in \cite{Alksspde} and is alluded to in Remark 1.2 in \cite{Abtbmsie}.
\item Even more, are the solutions to \eqref{nlksa} and \eqref{spidea} spatially continuously differentiable? and is the one dimensional spatial exponent $3/2$ in the solutions H\"older exponent $(3/2\wedge1)^{-}$ in \eqref{holregspide} indicating that the gradient of these solutions is nearly locally H\"older $1/2$ in space?  Also, what is the temporal H\"older regularity of the gradient?
\item If the answer to the first two parts of (b) is yes, what are the moduli of continuity of the gradient of the solutions to \eqref{nlksa} and \eqref{spidea}, in the space and time variables, respectively?
\een
\item[(Q2)] In \cite{Atfhosie}, it was established that, for each $d\in\{1,2,3\}$,  the time-fractional SIEs (time fractional SPIDEs in \eqref{spidea}) all have the same spatial H\"older regularity---$(\tf{4-d}{2}\wedge 1)^{-}$---for all $\beta\in\{1/2^{k};k\in\N\}$, as is clear from \eqref{holregspide}.  Is the spatial modulus of continuity a more discriminating measure of regularity that depends on $\beta\in\{1/2^{k};k\in\N\}$, and is $\beta=1/2$ critical on $(0,1/2]$?
\item[(Q3)]  It was established in \cite{Abtbmsie,Atfhosie,Alksspde} that the L-KS SPDE \eqref{nlksa} and the $\beta=1/2$ time-fractional SIEs (time fractional SPIDEs in \eqref{spidea}) have identical spatio-temporal H\"older regularity.  Does the continuity modulus capture the rougher regularity for the case $\beta=1/2$ time-fractional SPIDEs \eqref{spidea} (since, by footnote \ref{FN:btbmpdefo}, \eqref{spidea} is also associated with the rougher \emph{positive} bi-Laplacian PDE \eqref{btppdefo})?
\item[(Q4)] What are the exact spatio-temporal moduli of continuity for \eqref{nlks} and \eqref{spide} in $d=1,2,3$.
\item[(Q5)] What are the temporal probability laws associated with L-KS SPDEs and time-fractional SPIDEs?
\eit

\subsection{Main results: answering the questions}
We answer all of the above questions at length in the $a\equiv1$ Gaussian case for our two classes of equations in our main results\footnote{In $d=2$, we obtain a sharp upper bound on the uniform and local spatial moduli of continuity for the two classes of equations.}, which we now present.  First, we deal with the L-KS SPDEs.
\subsubsection{Exact moduli of continuity of L-KS SPDEs and their gradient,
and the bifractional Brownian motion link}

We start with the temporal regularity and probability law for L-KS SPDE \eqref{nlks} in spatial dimensions $d=1,2,3$. Recall that, given constants $H\in (0,1)$ and $K\in (0,1]$,
the bifractional Brownian motion $(B^{H,K}_{t})_{t \in [0,T]}$, introduced by Houdr\'e and Villa in \cite{HoudVilla03}, is a centered Gaussian process with covariance
  \begin{equation}
    \label{cov-bi}
    R^{H,K}(t,s) := R(t,s)= \frac{1}{2^{K}}\left( \left(
        t^{2H}+s^{2H}\right)^{K} -\vert t-s \vert ^{2HK}\right), \hskip0.5cm s,t \in [0,T].
  \end{equation}
We refer to \cite{LN09, RT, TuX07}  for various properties
of this process.
\bfr
\bthm[Temporal moduli of continuity and bi-fBM connection for the L-KS SPDE
in $d=1,2,3$]\lbl{tempmodlks}
Fix $(\vep,\vth)\in\Rp\times\R$ and $x\in\Rd$, and assume $d\in\{1,2,3\}$.
Assume that $(U,\sW)$ is the unique solution to \eqref{nlks} on $\OFFtP$,
with $\un$ satisfying \eqref{init} (a).
\ben\rencomrom
\item There are dimension-dependent constants $k_{1}^{(d)}>0$ and $k_{2}^{(d)}>0$,
independent of $x$, such that
\ben\rencomalp
\item $($Uniform temporal modulus$)$  for any compact interval $\tint\subset\Rp$
\beq
\P\lbk\lim_{\delta\searrow0} \sup_{\substack{|t-s|<\delta\\s,t\in\tint} }\frac{\lab U(t,x)-U(s,x)\rab}{|t-s|^{\tf{4-d}{8}}\sqrt{\log\lbk1/|t-s|\rbk}}=k_{1}^{(d)}\rbk=1,
\eeq
\item $($Local temporal modulus$)$  and for any fixed $t\ge 0$
\beq \label{Eq:time-LIL}
\P\lbk\lim_{\delta\searrow0} \sup_{\substack{|t-s|<\delta} }\frac{\lab U(t,x)-U(s,x)\rab}{\delta^{\tf{4-d}{8}}\sqrt{\log\log\lbk1/\delta\rbk}}=k_{2}^{(d)}\rbk=1.
\eeq
\een
\item $(B^{(H,K)} link)$ Assume $\vth=0$ in \eqref{nlks}, then, $U(\cdot,x)\eqL c_{d}B^{\lpa\frac12,\frac{4-d}{4}\rpa}$, where
\beq\lbl{bifbmlksconst}
c_{d}=(2\pi)^{-d/2} \lpa\tf{8}{\varepsilon}\rpa^{d/8} \frac{2^{\tf{d-4}{8}}}{\sqrt{2-d/2}}\sqrt{\int_{\mathbb{R}^{d}} \e^{-|\xi|^{4}} d\xi};\ d=1,2,3.
\eeq
In particular, we have the following Chung's law of the iterated logarithm for $U(\cdot,x)$:
\beq\lbl{chunglil}
\liminf_{r\searrow0}\frac{\max_{t\in[0,r]}\lab U(t,x)\rab}
{\lbk r^{{(4 -d)}/{8}}\rbk/\lbk\log\log(1/r)\rbk^{(4 -d)/8}}=k^{(d)}_{3}
\eeq
for every $x\in\Rd$ and for some positive finite $d$-dependent constant $k^{(d)}_{3}$.
\een
\ethm
\efr
\brm\lbl{dimensionconst}
\thmref{tempmodlks} establishes the temporal modulus of continuity part of Q4 
and answers Q5 (when $\vth=0$) for  L-KS SPDEs.
We observe that since $\un$ is assumed sufficiently smooth and deterministic, 
the deterministic part of \eqref{nlks} is $\C^{1,4}(\Rp,\Rd)$ smooth (see 
\cite{Aks,Alksspde}) and the modulus is controlled by the random parts of 
the SPDEs \eqref{nlks} (or their associated SIEs \eqref{ibtbapsol} below, 
with $a\equiv1$).  In addition to giving the precise dimension-dependent 
temporal modulus of continuity, \thmref{tempmodlks} says that, up to a constant, 
the simple ($\vth=0$) L-KS SPDE solution process $\lbr U(t,x), t\ge 0 \rbr$ has 
the same law as a bifractional Brownian motion with indices $H=\frac12$ and 
$K= 1 - \frac{d}{4}$. Thus, $U$ shares \emph{all} the temporal sample path 
properties with a $B^{\lpa1/2,(4-d)/{4}\rpa}$, in spatial dimensions $d=1,2,3$, 
which can be found in \cite{LN09, RT, TuX07}.
\erm

We next state our spatial modulus result for the L-KS SPDE \eqref{nlks}. 
 \thmref{spatmodlks}, along with Theorems 
\ref{spatmodtf} and \ref{spatmodtfhalf} below, give the first
instance of space-time white noise driven SPDEs that have a H\"older 
continuous gradient\footnote{See Appendix \ref{glossary} for the definition of $\C^{k,\gamma}(\R,\R)$ and other notations.}.
\bfr
\bthm[Spatial moduli of continuity for the L-KS SPDE in $d=1,2,3$]\lbl{spatmodlks}
Fix $(\vep,\vth)\in\Rp\times\R$ and fix $t\in\Rp$.  Assume that $(U,\sW)$ is the unique
solution to \eqref{nlks} on $\OFFtP$, with $\un$ satisfying \eqref{init} (a).
In the following, $k_{i}^{(d)}>0$ ($i = 4, 5$) are positive and finite constants
depending on $d$, $\vep,\, \vth$ and $t$.
\ben\rencomrom
\item If $d=1$, then $U(t,\cdot)\in\C^{1,\gamma}(\R;\R)$, almost surely, with the
H\"older exponent $\gamma \in(0,1/2)$.  Moreover,
\ben\rencomrom
\item $($Uniform spatial modulus$)$  for any compact rectangle $\sinto\subset\R$
\beq
\P\lbk\lim_{\delta\searrow0} \sup_{\substack{|x-y|<\delta\\x,y\in\sintth} }\frac{\lab\pa_{x} U(t,x)-\pa_{y}U(t,y)\rab}{|x-y|^{{1}/{2}}\sqrt{\log\lbk1/|x-y|\rbk}}=k_{4}^{(1)}\rbk=1,
\eeq
\item $($Local spatial modulus$)$  and for any fixed $x\in\R$
\beq\lbl{lsm1}
\P\lbk\limsup_{\delta\searrow0} \sup_{\substack{|x-y|<\delta} }\frac{\lab\pa_{x} U(t,x)-\pa_{y}U(t,y)\rab}{\delta^{1/{2}}\sqrt{\log\log\lbk1/\delta\rbk}}=k_{5}^{(1)}\rbk=1.
\eeq
\een
\item If $d =3$, then
\ben\rencomrom
\item $($Uniform spatial modulus$)$  for any compact rectangle $\sintth\subset\R^{3}$
\beq
\P\lbk\lim_{\delta\searrow0} \sup_{\substack{|x-y|<\delta\\x,y\in\sintth} }\frac{\lab U(t,x)-U(t,y)\rab}{|x-y|^{{1}/{2}}\sqrt{\log\lbk1/|x-y|\rbk}}=k_{4}^{(3)}\rbk=1,
\eeq
\item $($Local spatial modulus$)$  and for any fixed $x\in\R^{3}$
\beq\lbl{lsm3}
\P\lbk\limsup_{\delta\searrow0} \sup_{\substack{|x-y|<\delta} }
\frac{\lab U(t,x)-U(t,y)\rab}{\delta^{1/{2}} \sqrt{\log\log\lbk1/\delta\rbk}}=k_{5}^{(3)}\rbk=1.
\eeq
\een
\item If $d=2$, then
\ben\rencomrom
\item $($Uniform spatial modulus$)$  for any compact rectangle $\sintt\subset\R^{2}$
\beq \label{Eq:d2univ}
\P\lbk\lim_{\delta\searrow0} \sup_{\substack{|x-y|<\delta\\x,y\in\sintt} }
\frac{\lab U(t,x)-U(t,y)\rab}{|x-y| \log\lbk1/|x-y|\rbk}\le k_{4}^{(2)}\rbk=1,
\eeq
\item $($Local spatial modulus$)$ and for any fixed $x\in\R^{2}$
\beq\lbl{lsm2}
\P\lbk\limsup_{\delta\searrow0} \sup_{\substack{|x-y|<\delta} }
\frac{\lab U(t,x)-U(t,y)\rab}{\delta\sqrt{\log \lbk1/\delta\rbk \log\log\lbk1/\delta\rbk}}
\le k_{5}^{(2)}\rbk=1.
\eeq
\een
\een
\ethm
\efr

\brm\lbl{}
\thmref{spatmodlks} answers the spatial modulus of continuity part of Q4 and
answers Q1 for L-KS SPDEs. When $d=1$, \thmref{spatmodlks} significantly
refines the H\"older conclusion of Theorem 1.1 in \cite{Alksspde} from $U$
being nearly locally Lipschitz in space to continuously differentiable (and
hence locally Lipschitz) in space. Moreover, it also says that, for
$d=1$ and for any fixed time $t$, the spatial derivative of the solution
to the L-KS SPDE \eqref{nlks}, $\pa_{x} U(t,x)$, is nearly locally H\"older
$1/2$ (has H\"older exponent $\gamma\in(0,1/2)$) in space. In addition,
spatially, \thmref{spatmodlks} gives the exact uniform and local moduli of
continuity for the gradient $\pax U$ in $d=1$; the exact uniform
and local moduli of continuity of $U$ in $d=3$; and sharp upper bounds on these
moduli of continuity of $U$ in $d=2$. We note that moduli
of continuity of $U$ in the $d=2$ case are different from those for $d=3$
and the sample functions are nearly locally Lipschitz. However, we believe
that, unlike in the case of $d=1$, the sample function $x \mapsto U(t,x)$ is
nowhere differentiable in the case of $d=2$.

For the case of $d=2$, proving the nondifferentiability and the exact
spatial moduli of continuity will need substantial extra
work because, as a main technical tool for studying these problems, the
property of strong local nondeterminism has only been proved in \cite{X07,X09,XX11}
for Gaussian random fields with (directional) H\"older exponents smaller than 1.
See Remarks below for further information. We will study these and some related
problems in a subsequent paper.

The comparative question Q3 will be answered completely after stating the
corresponding results for time-fractional SPIDEs (Theorems \ref{spatmodtf}
and \ref{spatmodtfhalf} below).
\erm

The last main result for L-KS SPDEs gives the sharp temporal H\"older and the 
exact temporal continuity modulus regularity for the spatial gradient of the 
L-KS SPDE.  Let $\H^{\gamma_{*}^{-}}(\Rp;\R)$ be the space of locally H\"older 
continuous functions $f:\Rp\to\R$ whose H\"older exponent $\gamma\in(0,\gamma_{*})$.
\bfr
\bthm[Sharp temporal H\"older and exact continuity moduli for the L-KS SPDE gradient]
\lbl{thm:gradtemplks}Assume $d=1$ and fix $(\vep,\vth)\in\Rp\times\R$ and $x\in\R$.    
Assume that $(U,\sW)$ is the unique solution to \eqref{nlks} on $\OFFtP$, with $\un$ 
satisfying \eqref{init} (a).  Then, $\pax U(\cdot.x)\in\H^{(1/8)^{-}}(\Rp;\R)$, 
almost surely.  Moreover, there exist constants $k_i \in (0, \infty)$ ($i = 6, 7$)
such that
\ben\rencomrom
\item $($Uniform temporal modulus$)$ for any compact interval $\tint\subset\Rp$
\beq
\P\lbk\lim_{\delta\searrow0} \sup_{\substack{|t-s|<\delta\\t,s\in\tint} }\frac{\lab\pa_{x} U(t,x)-\pa_{x}U(s,x)\rab}{|t-s|^{{1}/{8}}\sqrt{\log\lbk1/|t-s|\rbk}}= k_6\rbk=1,
\eeq
\item $($Local temporal modulus$)$  and for any fixed $t\in\Rp$
\beq\lbl{ltmgradlks}
\P\lbk\lim_{\delta\searrow0} \sup_{\substack{|t-s|<\delta} }\frac{\lab\pa_{x} U(t,x)-\pa_{x}U(s,x)\rab}{\delta^{1/{8}}\sqrt{\log\log\lbk1/\delta\rbk}}=k_7\rbk=1.
\eeq
\een
\ethm
\efr
\brm\lbl{rm:lksgradvsheat}
Theorems \ref{spatmodlks} and \ref{thm:gradtemplks} not only tell us that, when $d=1$, 
the L-KS SPDE \eqref{nlks} gradient $\pax U$ exists and is continuous, but they also 
give us thorough spatio-temporal regularity results for $\pax U$, in both the 
H\"older and modulus senses.  This contrasts starkly with the standard second order 
heat SPDE whose solution is only spatially H\"older continuous with exponent 
$\gamma\in(0,1/2)$.  The spatial gradient spatio-temporal H\"older regularity in 
Theorems \ref{spatmodlks} and \ref{thm:gradtemplks} tell us that the gradient of 
L-KS SPDEs, $\pax U$, is rougher ($\gamma$-H\"older with $\gamma\in(0,1/2)$) in 
space than the continuously differentiable  (and hence Lipschitz) solution $U$.  
More surprisingly, $\pax U$ is also rougher in \emph{time} than $U$ ($\gamma$-H\"older 
with $\gamma\in(0,1/8)$ vs.~$\gamma\in(0,3/8)$ as in Allouba \cite{Alksspde}).  
Compared to the second order heat SPDE, the L-KS gradient $\pax U$ has the same 
spatial H\"older regularity as that of the solution to the heat SPDE; and $\pax U$ 
is twice as rough (half as smooth) as the heat SPDE solution in time, with H\"older 
exponent $\gamma\in(0,1/8)$ vs.~the well-known $\gamma\in(0,1/4)$ for the heat SPDE.   
Similar comments apply with respect to the moduli of continuity of $\pax U$ as compared 
to those of $U$ and to the heat SPDE (see Meerschaert, Wang, and Xiao \cite{MWX} 
for the heat SPDE moduli of continuity).
\erm

\subsubsection{Exact moduli of continuity for the time-fractional SPIDEs
and their gradient, and their temporal fractional laws}
We now turn to the regularity of the $\beta$-time-fractional SPIDEs \eqref{spide}
(and their corresponding time-fractional SIEs) and to their temporal fractional law.
The  uniform and local temporal continuity moduli for $U_{\beta}$, as well as the
law governing the behavior of the solution process $U_{\beta}(\cdot,x)=\lbr 
U_{\beta}(t,x);t\ge0\rbr$ of \eqref{spide}, are given by the next result.

To fully state the next theorem, we need to recall the definition of the
generalized hypergeometric (or simply the hypergeometric) function $_{p}F_{q}(a_{1},\ldots,a_{p};b_{1},\ldots,b_{q};z)$:
\beq\lbl{hypergeo}
_{p}F_{q}(a_{1},\ldots,a_{p};b_{1},\ldots,b_{q};z):=\sum_{n=0}^{\infty}\frac{(a_{1})_{n}
\cdots(a_{p})_{n}}{(b_{1})_{n}
\cdots(b_{q})_{n}}\frac{z^{n}}{n!};\ \ \ a_{i}, b_{i}, \ \mbox{ and }z\in\R,
\eeq
whenever the series in the right hand side of \eqref{hypergeo} converges, where
\beqs
\bsp
(u)_{0}=1\mbox{ and }
(u)_{n}=u(u+1)\cdots(u+n-1); u\in\R,\mbox{ and } n\ge1.
\end{split}
\eeqs
We are now ready for our result.
\bfr
\bthm[Temporal moduli of continuity and laws of the $\beta$-time-fractional
SPIDEs in $d=1,2,3$]\lbl{tempmodtf}
Let $\beta\in(0,1/2]$ and let $(U_{\beta},\sW)$ be the unique solution
to \eqref{spide} on $\OFFtP$, with $\un$ satisfying \eqref{init} (b).
Fix $x\in\Rd$ ($d=1,2,3$) and let $H = \frac{2 - \beta d}{4}$.
\ben\rencomrom
\item There exist constant $k_{i}^{(\beta,d)}>0$ ($i=8, 9$), depending on 
$\beta$ and $d$ but independent of $x$, such that
\ben\rencomalp
\item $($Uniform temporal modulus$)$ for any compact interval $\tint\subset\Rp$
\beq \label{Eq:SPI_umod}
\P\lbk\lim_{\delta\searrow0} \sup_{\substack{|t-s|<\delta\\s,t\in\tint} }\frac{\lab U_{\beta}(t,x)-U_{\beta}(s,x)\rab}
{|t-s|^{H}\sqrt{\log\lbk1/|t-s|\rbk}}
=k_{8}^{(\beta,d)}\rbk=1,
\eeq
\item $($Local temporal modulus$)$  and for any fixed $t\ge 0$
\beq \label{Eq:SPI_LIL}
\P\lbk\lim_{\delta\searrow0} \sup_{\substack{|t-s|<\delta} }
\frac{\lab U_{\beta}(t,x)-U_{\beta}(s,x)\rab}
{\delta^{H}\sqrt{\log\log\lbk1/\delta\rbk}}
=k_{9}^{(\beta,d)}\rbk=1.
\eeq
\een
\item $($Law of the $\beta$-time-fractional SPIDE$)$ The $\beta$-time-fractional
SPIDE solution process $\lbr U_{\beta}(t,x), t\ge 0 \rbr$ is a mean-zero Gaussian
process with covariance $\E\lbk U_{\beta}(t,x)U_{\beta}(s,x)\rbk$ given by
\beq\lbl{tftempcov}
(2\pi)^{-d} \int_{\mathbb{R}^{d}}  \sum_{k=0}^{\infty}\lbk\sum_{j=0}^{k}
\frac{t^{\beta j}s^{\beta(k-j)+1}\;
_{2}F_{1}\lpa1,-\beta j;2+\beta(k-j);\frac{s}t\rpa}{\lbk\beta(k-j)+1\rbk\Gamma(1+\beta j)
\Gamma(1+\beta(k-j))}\rbk\frac{(-1)^{k}\lab\xi\rab^{2k}}{2^{k}}d\xi.
\eeq
 In particular, $\lbr U_{\beta}(t,x), t\ge 0 \rbr$ is self-similar with index $H$,
 but it is not a bifractional Brownian motion. When $\beta=1/2$, the BTBM
 SPIDE has a fundamentally different law from that of the L-KS SPDE.
\een
\ethm
\efr

\brm\lbl{lawcontrast}  \thmref{tempmodtf} answers the temporal modulus of continuity
part of Q4 and answers Q5 for time-fractional SPIDEs. In addition to the precise
temporal continuity moduli of L-KS SPDEs \eqref{nlks}, \thmref{tempmodtf} gives the
first contrasting behaviors of \eqref{nlks} and the half-derivative or Brownian-time
Brownian motion SPIDE (\eqref{spide} with $\beta=1/2$).

The fundamental difference between the Gaussian laws of the time-fractional SPIDEs and
the L-KS SPDE, even at $\beta=1/2$, is most easily seen in the fourth order Brownian-time
Brownian motion ($\beta=1/2$) PDE, obtained first in \cite{Abtp1,Abtp2}
\footnote{We alternate freely between the notations $\pa^{n}_{x_{i}}f(x_{1},\ldots,x_{N})$
and $\pa^{n} f/\pa x^{n}_{i}$, $i=1,\ldots,N$.}:
\beq\lbl{btppdefo1}
\bc
\pa_{t} u= \ds{\df{\Delta\un }{\sqrt{8\pi t}}+\df18\Delta^2 u}, \ \
&(t,x)\in(0,\infty)\times\Rd; \cr
u(0,x) = \unx, &x\in\Rd.
\ec
\eeq
The memory term $\ds\df{\Delta\un }{\sqrt{8\pi t}}$ in the deterministic BTBM PDE
\eqref{btppdefo}, which is not shared with the deterministic version of
\eqref{nlksa} ($a\equiv0$), and the opposite sign of the bi-Laplacian in \eqref{btppdefo1}
vs.~that in the L-KS PDE are manifestations of the fundamental reason why the L-KS
SPDE and BTBM SPIDE have different laws.
\erm

We next state our spatial modulus result for the $\beta$-time-fractional SPIDEs
\eqref{spide}. We will distinguish the cases $0<\beta<1/2$ and $\beta = 1/2$, where
subtle differences arise. It is interesting to notice that, for $0<\beta<1/2$, the
spatial moduli of SPIDEs \eqref{spide} are identical, modulo constants, to those of
the L-KS SPDEs \eqref{nlks}.
\bfr
\bthm[Spatial moduli of continuity for the $\beta$-time-fractional SPIDEs for
$0<\beta<1/2$ and $d=1,2,3$]\lbl{spatmodtf}
Assume that $(U_{\beta},\sW)$ is the unique solution to \eqref{spide} on $\OFFtP$,
with $\un$ satisfying \eqref{init} (b). We assume $t\in\Rp$ is fixed and $0<\beta<1/2$.
In the following, $k_{i}^{(\beta,d)}>0$ ($i = 10, 11$) are constants depending
on $d$, $t$ and $\beta$.
\ben\rencomrom
\item If $d=1$, then $U_{\beta}(t,\cdot)\in\C^{1,\gamma}(\R;\R)$, almost surely,
with the H\"older exponent $\gamma\in(0,1/2)$.
Moreover,
\ben\rencomalp
\item $($Uniform spatial modulus$)$  for any compact rectangle $\sinto\subset\R$
\beq
\P\lbk\lim_{\delta\searrow0} \sup_{\substack{|x-y|<\delta\\x,y\in\sinto} }
\frac{\lab\pa_{x} U_{\beta}(t,x)-\pa_{y}U_{\beta}(t,y)\rab}
{|x-y|^{{1}/{2}}\sqrt{\log\lbk1/|x-y|\rbk}}=k_{10}^{(\beta,1)}\rbk=1,
\eeq
\item $($Local spatial modulus$)$  and for any fixed $x\in\R$
\beq\lbl{lsmtf1}
\P\lbk\limsup_{\delta\searrow0} \sup_{\substack{|x-y|<\delta} }
\frac{\lab\pa_{x} U_{\beta}(t,x)-\pa_{y}
U_{\beta}(t,y)\rab}{\delta^{1/{2}}\sqrt{\log\log\lbk1/\delta\rbk}}
=k_{11}^{(\beta,1)}\rbk=1.
\eeq
\een
\item If $d=3$, then
\ben\rencomrom
\item $($Uniform spatial modulus$)$  for any compact rectangle $\sintth\subset\R^{3}$
\beq
\P\lbk\lim_{\delta\searrow0} \sup_{\substack{|x-y|<\delta\\x,y\in\sintth} }\frac{\lab U_{\beta}(t,x)-U_{\beta}(t,y)\rab}
{|x-y|^{{1}/{2}}\sqrt{\log\lbk1/|x-y|\rbk}}=k_{10}^{({\beta},3)}\rbk=1,
\eeq
\item $($Local spatial modulus$)$  and for any fixed $x\in\R^{3}$
\beq\lbl{lsmtf3}
\P\lbk\limsup_{\delta\searrow0} \sup_{\substack{|x-y|<\delta} }
\frac{\lab U_{\beta}(t,x)-U_{\beta}(t,y)\rab}
{\delta^{1/{2}}\sqrt{\log\log\lbk1/\delta\rbk}}=k_{11}^{(\beta,3)}\rbk=1.
\eeq
\een
\item If $d=2$, then
\ben\rencomrom
\item $($Uniform spatial modulus$)$  for any compact rectangle $\sintt\subset\R^{2}$
\beq
\P\lbk\lim_{\delta\searrow0} \sup_{\substack{|x-y|<\delta\\x,y\in\sintt} }\frac{\lab U_{\beta}(t,x)-U_{\beta}(t,y)\rab}
{|x-y|  \log\lbk1/|x-y|\rbk }\le k_{10}^{(\beta,2)}\rbk=1,
\eeq
\item $($Local spatial modulus$)$  and for any fixed $x\in\R^{2}$
\beq\lbl{lsmtf2}
\P\lbk\limsup_{\delta\searrow0} \sup_{\substack{|x-y|<\delta} }
\frac{\lab U_{\beta}(t,x)-U_{\beta}(t,y)\rab}
{\delta \sqrt{ \log\lbk1/\delta\rbk \log\log\lbk1/\delta\rbk}}\le
k_{11}^{(\beta,2)}\rbk=1.
\eeq
\een
\een
\ethm
\efr

When $\beta=1/2$, the next result shows that the (BTBM) SPIDE  is critical, signaling the onset of rougher spatial sample paths, in $d=1,2,3$, than both the time-fractional SPIDEs with $\beta<1/2$ and the L-KS SPDE\footnote{Carefully examining \thmref{spatmodtfhalf},
we see the extra $1/\sqrt{\log(1/|x-y|)}$ (or $1/\sqrt{\log(1/\delta)}$) term in each
modulus expression as compared to the corresponding expressions in both \thmref{spatmodtf}
and \thmref{spatmodlks} above.}.  This is despite the fact that they all have the same spatial H\"older regularity as established first in Allouba \cite{Abtbmsie,Atfhosie,Alksspde}.
\bfr
\bthm[Spatial continuity modulus for the critical half-derivative BTBM SPIDEs,
$\beta=1/2$ and $d=1,2,3$]\lbl{spatmodtfhalf}
Assume that $(U_{1/2},\sW)$ is the unique solution to \eqref{spide} on $\OFFtP$,
with $\beta=1/2$ and $\un$ satisfying \eqref{init} (b).  In the following,
$k_{i}^{(1/2,d)}>0$ ($i = 10, 11$) are constants depending on $d$ and $t$.
\ben\rencomrom
\item If $d=1$, then $U_{1/2}(t,\cdot)\in\C^{1,\gamma}(\R;\R)$, almost surely,
with the H\"older exponent $\gamma\in(0,1/2)$. Moreover,
\ben\rencomrom
\item $($Uniform spatial modulus$)$  for any compact rectangle $\sinto\subset\R$
\beq
\P\lbk\lim_{\delta\searrow0} \sup_{\substack{|x-y|<\delta\\x,y\in\sinto} }\frac{\lab\pa_{x} U_{1/2}(t,x)-\pa_{y}
U_{1/2}(t,y)\rab}{|x-y|^{{1}/{2}}{\log\lbk1/|x-y|\rbk}}=k_{10}^{(1/2,1)}\rbk=1,
\eeq
\item $($Local spatial modulus$)$  and for any fixed $x\in\R$
\beq\lbl{lsmhalfder}
\P\lbk\limsup_{\delta\searrow0} \sup_{\substack{|x-y|<\delta} }\frac{\lab\pa_{x} U_{1/2}(t,x)-\pa_{y}U_{1/2}(t,y)\rab}
{\delta^{1/{2}}\sqrt{\log[1/\delta]\log\log\lbk1/\delta\rbk}}=k_{11}^{(1/2,1)}\rbk=1.
\eeq
\een
\item If $d=3$, then
\ben\rencomrom
\item $($Uniform spatial modulus$)$  for any compact rectangle $\sintth\subset\R^3$
\beq
\P\lbk\lim_{\delta\searrow0} \sup_{\substack{|x-y|<\delta\\x,y\in\sintth} }\frac{\lab U_{1/2}(t,x)-U_{1/2}(t,y)\rab}
{|x-y|^{{1}/{2}}{\log\lbk1/|x-y|\rbk}}=k_{10}^{(1/2,3)}\rbk=1,
\eeq
\item $($Local spatial modulus$)$  and for any fixed $x\in\R^3$
\beq\lbl{lsmtf3half}
\P\lbk\limsup_{\delta\searrow0} \sup_{\substack{|x-y|<\delta} }\frac{\lab U(t,x)-U(t,y)\rab}
{\delta^{1/{2}}\sqrt{\log[1/\delta]\log\log\lbk1/\delta\rbk}}=k_{11}^{(1/2,3)}\rbk=1.
\eeq
\een
\item If $d=2$, then
\ben\rencomrom
\item $($Uniform spatial modulus$)$  for any compact rectangle $\sintt\subset\R^2$
\beq
\P\lbk\lim_{\delta\searrow0} \sup_{\substack{|x-y|<\delta\\x,y\in\sintt} }\frac{\lab U(t,x)-U(t,y)\rab}{|x-y|{(\log\lbk1/|x-y|\rbk)^{3/2}}}\le k_{10}^{(1/2,2)}\rbk=1,
\eeq
\item $($Local spatial modulus$)$  and for any fixed $x\in\R^2$
\beq\lbl{lsmtf2half}
\P\lbk\limsup_{\delta\searrow0} \sup_{\substack{|x-y|<\delta} }\frac{\lab U(t,x)-U(t,y)\rab}
{\delta\,\log[1/\delta]\sqrt{\log\log\lbk1/\delta\rbk}}\le k_{11}^{(1/2,2)}\rbk =1.
\eeq
\een
\een
\ethm
\efr

The last main result for time-fractional SPIDEs gives the sharp temporal 
H\"older and the exact temporal continuity modulus regularity for the spatial 
gradient of these time-fractional SPIDEs.
\bfr
\bthm[Sharp temporal H\"older and exact continuity moduli for the time-fractional 
SPIDEs gradient]\lbl{thm:gradtempspide}Assume $d=1$, $x\in\R$, and let $\beta\in(0,1/2]$.  
Assume that $(U_{\beta},\sW)$ is the unique solution to \eqref{spide} on $\OFFtP$, 
with $\un$ satisfying \eqref{init} (b).  Then, $\pax U_{\beta}(\cdot.x)\in\H^{\lpa{\lpa2 -3\beta \rpa}/{4}\rpa^{-}}(\Rp;\R)$,
almost surely.  Moreover, there exist constants $ k_i \in (0, \infty)$ ($i = 12, 13$) such that
\ben\rencomrom
\item $($Uniform temporal modulus$)$  for any compact interval $\tint\subset\Rp$
\beq
\P\lbk\lim_{\delta\searrow0} \sup_{\substack{|t-s|<\delta\\t,s\in\tint} }\frac{\lab\pa_{x} U_{\beta}(t,x)-\pa_{x}U_{\beta}(s,x)\rab}{|t-s|^{{\tf{2-3\beta}{4}}}
\sqrt{\log\lbk1/|t-s|\rbk}}=k_{12}\rbk=1,
\eeq
\item $($Local temporal modulus$)$  and for any fixed $t\in\Rp$
\beq\lbl{lsm1}
\P\lbk\lim_{\delta\searrow0} \sup_{\substack{|t-s|<\delta} }\frac{\lab\pa_{x} U_{\beta}(t,x)-\pa_{x}U_{\beta}(s,x)\rab}{\delta^{{\tf{2-3\beta}{4}}}
\sqrt{\log\log\lbk1/\delta\rbk}}= k_{13}\rbk=1.
\eeq
\een
\ethm
\efr
\brm\lbl{rm:spidegrad}
For $d = 1$, Theorems \ref{spatmodtf}, \ref{spatmodtfhalf}, and \ref{thm:gradtempspide} 
give us the existence, as well as thorough results on the spatio-temporal moduli of 
continuity for the gradient $\pax U_{\beta}$ of the $\beta$-time-fractional SPIDEs 
\eqref{spide}. Spatially, Theorems \ref{spatmodtf} and
\ref{spatmodtfhalf} say that, even though the H\"older exponent of the gradient 
$x \mapsto \pax U_{\beta}(t, x)$ is the same for all $\beta\in(0,1/2]$, the exact 
uniform and local moduli of continuity for $\beta \in (0, 1/2)$ and $\beta = 1/2$.
The SPIDEs gradient $\pax U_{\beta}$ is spatially rougher in the modulus sense at 
$\beta=1/2$ than it is for $\beta<1/2$; and $\pax U_{\beta}$ has the same spatial modulus 
of continuity as that of the L-KS SPDE gradient for $\beta<1/2$. Theorem \ref{thm:gradtempspide}
shows that the time-fractional SPIDE gradient $\pax U_{\beta}$,
for $\beta =1/2$, has the same H\"older and modulus regularity in the time variable as
the L-KS SPDE gradient. Moreover, the temporal H\"older exponent $\gamma \nearrow 1/2$ 
[i.e., the temporal H\"older regularity increases] as $\beta\searrow0$. This is 
consistent with the similar phenomenon for the
time-fractional solutions observed by Allouba in \cite{Atfhosie}.  
\erm

Theorems \ref{tempmodlks}--\ref{thm:gradtempspide} together answer all the
questions Q1--Q5 above except that, in the case of $d=2$, extra work
will be needed for completely establishing the exact spatial moduli of
continuity.
\subsection{The strong local nondeterminism property and modulus of continuity}
The proofs of Theorems 1.1 -- 1.5 depend on the results and methods in Meerschaert,
Wang, and Xiao \cite{MWX}, Xue and Xiao \cite{XX11} which, in turn,  are based on
general Gaussian methods (cf. e.g., \cite{MarRos}) and the properties of strong
local nondeterminism in \cite{X07,X09}.  More specifically we obtain an expression
for the spectral measure/density associated with the solution and use it
to prove the exact uniform and local moduli of continuity.

To determine many sample path properties of our Gaussian solution $U$ to our
SIEs separately and jointly in time and space, the following second moments of
spatial and temporal differences are crucial\footnote{In the case of $\beta$
time-fractional SIEs, these quantities depend also on $\beta$.}:
 \beq\lbl{smom}
 \bsp
 \tsig(s,t;x)^{2}&=\E[U(t,x)-U(s,x)]^{2}; \quad s,t\in\Rp, x\in\Rd,
 \\ \ssig(t;x,y)^{2}&=\E[U(t,x)-U(t,y)]^{2}; \quad t\in\Rp, x,y\in\Rd,
 \\ \jsig(s,t;x,y)^{2}&=\E[U(t,x)-U(s,y)]^{2}; \quad s,t\in\Rp, x,y\in\Rd.
 \end{split}
 \eeq
Xiao \cite{X09} gave some general conditions for effectively studying several
analytic and geometric properties of Gaussian random fields.
For convenience of readers, we restate these conditions below, adapting
the notation slightly to our setting.  Let $\tint=[a,b]$ and
$\sint=\otimes_{k=1}^{d}[a_{k},b_{k}]$ be one and $d$-dimensional
closed intervals in $\Rp$ and $\Rd$, respectively\footnote{In this paper,
unless otherwise stated we take $0\le a,a_{k}<1$ and $b=b_{k}=1$, for all $k$.}.
Let $\gamma=(\gamma_{1},\cdots,\gamma_{d+1}) \in(0,1]^{d+1}$ be a fixed vector,
and denote by $\rho$ the metric on $\Rp\times\Rd$ given by
 \beq\lbl{Xiaorho}
 \rho(s,t;x,y)=|t-s|^{\gamma_{1}}+\sum_{j=2}^{d}\lab x_{j}-y_{j}\rab^{\gamma_{j}};
 \ \ s,t\in\Rp,\ x,y\in\Rd.
 \eeq
\ben
\item[(C1)]  (Spatio-temporal bounds) There exist positive and finite constants
$c_{2,1}$ and $c_{2,2}$ such that
\beqs \label{Eq:C1b}
c_{2,1}\rho^{2}(s,t;x,y)\le \lbk\jsig(s,t;x,y)\rbk^{2}\le c_{2,2}\rho^{2}(s,t;x,y)
\eeqs
for all $s,t\in\tint$ and $x,y\in\sint$.
\item[(C2)]  (SLND) There exists a constant $c_{2,3}>0$ such that for all integers
$n\ge1$ and all $p,p^{(1)},\ldots,p^{(n)}\in\tint\times\sint$
\beqs \label{Eq:C2b}
\var \lpa U(p) | U(p^{(1)}),\ldots,U(p^{(n)})\rpa \ge c_{2,3} \sum_{j=1}^{d+1}
\min_{0\le k\le n}\lab p_{j}-p^{(k)}_{j}\rab^{2\gamma_{j}}.
\eeqs
\een
\brm\lbl{gammaexp}
In this article, $\gamma_{1}$ is the least upper bound for the temporal H\"older
exponents for our SPDEs/SPIDEs, and $\gamma_{2}=\cdots=\gamma_{d+1}$ are the least
upper bound for the spatial H\"older exponents for our SPDEs/SPIDEs. By Theorem 1.1
in Allouba \cite{Alksspde,Abtbmsie} $\gamma_{1}=(4-d)/8$ for the L-KS SPDE and for
the SPIDE \eqref{spide} when $\beta=1/2$; and, by Theorem 1.2 in Allouba \cite{Atfhosie}, $\gamma_{1}=(2\beta^{-1}-d)/{4{\beta^{-1}}}$ for the SPIDE \eqref{spide} for $\beta
=1/2^{k}, k\in\N$. Also, $\gamma_{j}=[({4-d})/{2}]\wedge 1$, $j=2,\ldots,d+1$, and
$d=1,2,3$ for all SPDEs/SPIDEs in this article by Theorem 1.1 and Theorem 1.2
\cite{Abtbmsie,Atfhosie} and by Theorem 1.1 in \cite{Alksspde}.
\erm
\brm\lbl{gammaexp2}
The main results of this article\footnote{The last result of this article is \thmref{SHcom}, which uses change of measure to give equivalence in law and to transfer regularity between linear L-KS SPDEs and time-fractional SPIDEs and their nonlinear versions, with Allen-Cahn type and polynomial nonlinearities, on compact time-space sets.}, Theorems \ref{tempmodlks}--\ref{thm:gradtempspide}, establish exact uniform and local moduli of continuity for the solutions of the L-KS SPDE and the SPIDE in the time variable $t$ and space variable $x$, separately.  Also, \thmref{tempmodlks} gives a Chung's law of iterated logarithm for simple L-KS SPDEs.  For proving these theorems, we will only use Conditions (C1) and (C2) for two special cases: either $x=y$ or $s=t$, respectively. Hence, the spectral conditions in Xiao
\cite{X07} can be applied to verify these conditions. Moreover, Theorems 2.1 and 2.5
in \cite{X07} allow us to prove more general properties by replacing the power functions
$|t-s|^{\gamma_{1}}$ and $|x_j-y_j|^{\gamma_{j}}$ in \eqref{Xiaorho} 
by regularly varying functions of $|s-t|$ or $|x-y|$ with regularity exponents smaller
than 1. Such an extension does not affect the proofs in Meerschaert,
Wang, and Xiao \cite{MWX}, hence the theorems in Sections 4 and 5 of
\cite{MWX} are still applicable.
\erm

\brm\lbl{gammaexp3}
It would be interesting to study analytic and geometric properties of the
solutions of L-KS SPDE and the SPIDE in both time and space variables $t$ and
$x$ simultaneously. For this purpose, the full strength of Conditions (C1)
and (C2) will be needed. The problems are more complicated and some new techniques
will be required. We will pursue this line of research in a separate article.
\erm

\subsection{Rigorous kernel stochastic integral equations formulations}\lbl{rigdisc}  
For the L-KS SPDE \eqref{nlksa}, as done in \cite{Alksspde,Abtpspde}, we use the linearized Kuramoto-Sivashinsky kernel introduced in \cite{Aks,Abtpspde,Alksspde} to define their rigorous mild SIE formulation.  This L-KS kernel is the fundamental solution to the deterministic version of \eqref{nlksa} ($a\equiv0$), as shown in
\cite{Aks,Abtpspde,Alksspde}, and is given by:
\beq\lbl{vepvthLKS}
\bsp
\KKSepthtxy&=\int_{-\infty}^{0}\df{\e^{\i\vth s} \e^{-|x-y|^2/2\i s}}{{\lpa2\pi \i s \rpa}^{d/2}}\peptsz ds+\int_{0}^{\infty}\df{\e^{\i\vth s} \e^{-|x-y|^2/2\i s}}{{\lpa2\pi \i s \rpa}^{d/2}}\peptsz ds,
\\&=(2\pi)^{-d}\int_{\Rd}\e^{-\frac{\vep t}{8}\left( -2\vth+\lab\xi\rab^{2} \right) ^{2}}\e^{\i \langle\xi, x-y\rangle}d\xi;
\\&= (2\pi)^{-d}\int_{\Rd}\e^{-\frac{\vep t}{8}\left( -2\vth+\lab\xi\rab^{2} \right) ^{2}}\cos\lpa{\langle\xi, x-y\rangle}\rpa d\xi;\ \ \vep>0,\ \vth\in\R.
\end{split}
\eeq
Let $b:\R\to\R$ be Borel measurable.  The rigorous L-KS kernel SIE (mild) formulation of the nonlinear drift-diffusion L-KS SPDE
\begin{equation} \label{nlksab}
 \begin{cases} \displaystyle\frac{\partial U}{\partial t}=
-\tfrac\vep8\lpa\lap+2\vth\rpa^{2}U+b(U)+a(U)\frac{\partial^{d+1} W}{\partial t\partial x}, \quad & (t,x)\in\Rpop\times\Rd;
\cr U(0,x)=\unx, & x\in\Rd,
\end{cases}
\end{equation}
is the stochastic integral equation 
\beq\lbl{ibtbapsol}
\bsp
U(t,x)=&\intrd\KKSepthtxy \uny dy
\\&+ \intrdzt \KKSepthtsxy \lbk b(U(s,y))dsdy+ a(U(s,y))\sW(ds\times dy)\rbk
\end{split}
\eeq
(see \cite[p.~530]{Abtpspde} and \cite[Definition 1.1, Eq.~ (1.11)]{Alksspde}).  Of course, the mild formulation of \eqref{nlks} is then obtained by setting $a\equiv1$ and $b\equiv0$ in \eqref{ibtbapsol}.

For the time fractional SPIDE \eqref{spidea}, we first explain heuristically the role of the extra time fractional integral $I^{1-\beta}$ in the formal formulation.  Succinctly, It compensates for the $\beta$ time fractional derivative ${}^{\mbox{\tiny C}}\pa_{t}^{\beta}$ so as to end up with a standard stochastic integral in time term with respect to space-time white noise, and wind up with a simpler (and smoother) SIE formulation.  To see this quickly before we give the formal computation, we first observe heuristically that to get a formulation with $U_{\beta}$ on the left hand side, we only need to get rid of the $\beta$ fractional time derivative ${}^{\mbox{\tiny C}}\pa_{t}^{\beta}$ by applying  a $\beta$ fractional integral $I_{t}^{\beta}$ to it.  This means we have to apply $I_{t}^{\beta}$ to the right side of the SPIDE \eqref{spidea} too.  So, if we want the time integral of the noise term to be of order $1$ (nonfractional), we need to have started already with a fractional integral $I_{t}^{1-\beta}$ of the noise so that $I_{t}^{\beta}\circ I_{t}^{1-\beta}=I_{t}^{1}$.  To put this heuristic on a firm ground and to get the SIE formulation of time fractional SPIDEs, we use Umarov's fractional Duhamel principle (see Theorem 3.6 in \cite{UmSa}), which we now proceed to describe.

If we replace the bracketed terms in the nonlinear drift-diffusion time-fractional  SPIDE
\beq\lbl{spideab}
\begin{cases} \displaystyle
{}^{\mbox{\tiny C}}\pa_{t}^{\beta}U_{\beta}=\tfrac12\D U_{\beta}+I_{t}^{1-\beta}\lbk b(U_{\beta})+a(U_{\beta})\frac{\partial^{d+1} W}{\partial t\partial x}\rbk, \quad & (t,x)\in\Rpop\times\Rd;
\cr U_{\beta}(0,x)=\unx, & x\in\Rd,
\end{cases}
\eeq	 
by a nice forcing term $f(t,x)$; then, using Theorem 3.6 in \cite{UmSa}, we obtain
\beq\lbl{umarovcomp}
\bsp
U_{\beta}(t,x)&=\intrd\Kbetatxy \uny dy+\int_{\Rd}\int_{0}^{t}\Kbetatsxy \lpa{}^{\mbox{\tiny R}}\pa_{t}^{1-\beta} I_{t}^{1-\beta} f(s,y)\rpa ds dy
\\&=\intrd\Kbetatxy \uny dy+\int_{\Rd}\int_{0}^{t}\Kbetatsxy  f(s,y) ds dy,
\end{split}
\eeq
where ${}^{\mbox{\tiny R}}\pa_{t}^{\alpha}$ is the Riemann-Liouville fractional derivatives of order $\alpha$:
\beq\lbl{RLder}
{}^{\mbox{\tiny R}}\pa_{t}^{\alpha}:=\bc\df{1}{\Gamma(1-\alpha)}\df{d}{dt}\int_{0}^{t}\df{f(\tau)}{(t-\tau)^{\alpha}}d\tau& \mbox{ if } 0<\alpha<1,\\
\df{d}{dt}f(t);& \mbox{ if } \alpha=1,
\ec
\eeq
where we used the fact that ${}^{\mbox{\tiny R}}\pa_{t}^{\alpha} I_{t}^{\alpha}=\mbox{Id}$, and where $\Kbetatx$ is the solution to the time-fractional PDE:
\beq\lbl{tfpde}
\begin{cases} \displaystyle
{}^{\mbox{\tiny C}}\pa_{t}^{\beta}U_{\beta}=\tfrac12\D U_{\beta}, \quad & (t,x)\in\Rpop\times\Rd;
\cr U(0,x)=\delta(x), & x\in\Rd,
\end{cases}
\eeq	
where $\delta(x)$ is the usual Dirac delta function.  These fundamental solutions $\Kbetatx$ are the densities of an inverse stable L\'evy time Brownian motion $B^{x}\lpa\Lambda_{\beta}(t)\rpa$, at time $t$, in which the inverse stable L\'evy motion $\Lambda_{\beta}$ of index $\beta$ acts as the time clock for an independent $d$-dimensional Brownian motion $B^{x}$ (see \cite{Atfhosie,Bert,MainLuch,MBS,MS}).  Thus,
\beq\lbl{eq:fundsolbeta}
\Kbetatx=\int_{0}^{\infty}\psx\KItzs ds;\ 0<\beta<1.
\eeq
In the case $\beta=1/2$, the kernel $\Kbetatx$ is the density of the Brownian-time Brownian motion as in \cite{Abtp1,Abtp2,Abtpspde}; and when $\beta\in\lbr1/2^{k};k\in\N\rbr$, the kernel $\Kbetatx$ is the density of $k$-iterated BTBM  as detailed in \cite{Atfhosie}.  Namely, denote by $$\mathbb{B}^{x}_{\underset{i=1}{\overset{k}{\bigcirc}} B_{i}}(t):=B^{x}\lpa\lab B_{k}\lpa\cdots B_{2}\lpa\lab B_{1}(t)\rab\rpa\cdots\rpa\rab\rpa$$ a $k$-iterated Brownian-time Brownian motion at time $t$; where $\lbr B_{i}\rbr_{i=1}^{k}$ are independent copies of a one dimensional scaled Brownian motion starting at zero, with density $\frac{1}{\sqrt{4\pi t}}\exp\left(-\frac{z^2}{4t}\right)=\lpa1/\sqrt2\rpa K^{\text{\tiny{\sc{BM}}}}_{t;0,z/\sqrt2}$, and independent from the standard $d$-dimensional Brownian motion $B^{x}$, which starts at $x\in\Rd$.  When
  $\beta^{-1}\in\lbr2^k;k\in\N\rbr$,  the density $\Kbetatx$ of $\mathbb{B}^{x}_{\underset{i=1}{\overset{k}{\bigcirc}} B_{i}}(t)$ is given by\footnote{We are using the convention $\prod_{i=0}^{-1}c_{i}=1$ for any $c_{i}$ and the convention $\int_{\Rp^{0}}f(s)ds=f(s)$, for every $f$.  Also, we use the convention that the case $k=0$ ($\beta=1$) is also the standard $d$-dimensional Brownian motion case.}.
\beq\lbl{betaker}
\Kbetatx=2^{\frac{k}{2}}\int_{\Rpop^{k}} \psoxz \ptzsk\prod_{i=0}^{k-2}\ptzskiimo ds_{1}\cdots ds_{k}.
\eeq
Now, denoting the white noise formally by $\dot{W}$ and replacing the nice forcing term $f$ by the bracketed terms in \eqref{spideab}, we see that \eqref{umarovcomp} becomes
\beq\lbl{umarovcompspide}
\bsp
U_{\beta}(t,x)&=\intrd\Kbetatxy \uny dy
\\&+\int_{\Rd}\int_{0}^{t}\Kbetatsxy \lbk b(U_{\beta}(s,y))dsdy+ a(U_{\beta}(s,y))\dot{W}(s,y) ds dy\rbk ,
\end{split}
\eeq
which is rigorously written as
\beq\lbl{isltbmsie}
\bsp
U_{\beta}(t,x)&=\intrd\Kbetatxy \uny dy
\\&+ \int_{\Rd}\int_{0}^{t}\Kbetatsxy \lbk b(U_{\beta}(s,y))dsdy+ a(U_{\beta}(s,y))\sW(ds\times dy)\rbk.
\end{split}
\eeq
Equation \eqref{isltbmsie}---with $a\equiv1$ and $b\equiv0$---is what we rigorously mean by the SPIDE \eqref{spide}, and it is the equation we work with.  Here, we call the stochastic integral equation in \eqref{isltbmsie} $\beta$-time-fractional SIE.  We stress here that, in the case $\beta=1/2$, equation \eqref{isltbmsie} is exactly equation (3.3) in \cite{Abtpspde} (when $a\equiv1$ and $b\equiv0$), the equation in \cite[top of p.~530]{Abtpspde}, and equation  (1.3) in \cite{Abtbmsie}; and, for $0<\beta<1$, equation \eqref{isltbmsie} is exactly equation (1.14) in \cite{Atfhosie}.  In the important case of $\beta=1/2$, we call \eqref{isltbmsie} the BTBM SIE since in this case $\Kbetatx$ is the density of a Brownian-time Brownian motion.
\bnt\lbl{notations}
Unless explicitly otherwise stated, $c$ and $C$ will denote constants whose value may change from a statement to another.   We refer the reader to the convenient end-of-paper list of notations.
\ent

\section{Kernels Fourier transforms}

The following lemma gives the spatial Fourier transform\footnote{In space, we are using the symmetric form of the Fourier transform $\hat{f}(\xi)=\lpa2\pi\rpa^{-\frac d2}\int_{\Rd}f(x)\e^{-\i\xi\cdot x}dx$.} of the $\beta$-time-fractional (including the $\beta=1/2$ BTBM case), and the $(\vep,\vth)$ L-KS kernels.
\blm[Spatial Fourier transforms]\lbl{FT}
Let $\KKSepthtx$ and $\Kbetatx$ be the $(\vep,\vth)$ LKS kernel and the $\beta$-time-fractional kernel,  respectively.
\ben\rencomrom
\item The spatial Fourier transform of the $(\vep,\vth)$ LKS kernel in \eqref{vepvthLKS} is given by
\beq\lbl{LKSF}
\FKKSepthtxi=\lpa2\pi\rpa^{-\frac d2}\e^{-\frac{\vep t}{8}\left( -2\vth+\lab\xi\rab^{2} \right) ^{2}};\ \ \ \vep>0,\ \vth\in\R.
\eeq

\item
Let $0<\beta<1$. The spatial Fourier transform of the $\beta$-time-fractional kernel is given by
\beq\lbl{betaf}
\FKbetatxi=\lpa2\pi\rpa^{-\frac d2}E_{\beta}\Big(-\frac{\lab\xi\rab^{2}}{2}t^{\beta}\Big),
\eeq
where
\beq\lbl{ML}
E_{\beta}(x)= \sum_{k=0}^{\infty}\frac{x^{k}}{\Gamma(1 + \beta k)},
\eeq
is the well known Mittag-Leffler function\footnote{See Haubold, Mathai, and Saxena \cite{HauMathSax11} and \cite{MathHau07} for the necessary background.}.  In particular, the Fourier transform of the BTBM density (the case $\beta=1/2$) is given by\footnote{\lbl{footFT}Strictly speaking, the $\beta=1/2$ BTBM Fourier transform
in \lemref{FT} is that of a BTBM in which the inner BM is time scaled.  The  Fourier transform of a standard BTBM is  $$\FKBtxi=\lpa2\pi\rpa^{-\frac d2}\e^{\frac t8\lab\xi\rab^{4}}\lbk\frac{2}{\sqrt{\pi}}\int_{\frac{\sqrt{2t}\lab\xi\rab^{2}}{4}}^{\infty}
\e^{-\tau^{2}}d\tau\rbk$$} 
\beq\lbl{btbmf}
\hat{\K}^{(1/2,d)}_{t,\xi}
=\lpa2\pi\rpa^{-\frac d2}\e^{\frac {t} 4\lab\xi\rab^{4}} \lbk\frac{2}{\sqrt{\pi}}\int_{\frac{\sqrt{t}\lab\xi\rab^{2}}{2}}^{\infty}\e^{-\tau^{2}}d\tau\rbk.
\eeq
\een
\elm
\bpf
The proof in the BTBM (the case $\beta=1/2$ or $k=1$) and the $(\vep,\vth)$
LKS kernels cases is given in \cite[Lemma 2.1]{Alksspde}.
We now prove the general $\beta$ case.  The kernel $\Kbetatx$ is 
given by \eqref{eq:fundsolbeta}.  Since the Laplace transform of $\KItzs$ in time is particularly simple and is given by $\LKIlzs=\theta^{\beta-1}\e^{-s\theta^{\beta}}$,
we easily get the Fourier transform of $\Kbetatx$ by first applying the Laplace transform in time and the Fourier transform in space to get
the following Laplace-Fourier transforrm for $\Kbetatx$
\beq\lbl{betaLFcomp}
\bsp
&\LFKbetathxi=\lpa2\pi\rpa^{-\frac d2}\int_{0}^{\infty}\lbr\int_{\Rd}\lbk\int_{0}^{\infty}\psx\KItzs ds\rbk \e^{-\i\xi\cdot x}dx\rbr \e^{-t\theta}dt
\\& =\lpa2\pi\rpa^{-\frac d2}\int_{0}^{\infty}\lbr\int_{\Rd}\psx\e^{-\i\xi\cdot x}dx\int_{0}^{\infty}\KItzs \e^{-t\theta}dt \rbr ds
\\&=\lpa2\pi\rpa^{-\frac d2}\int_{0}^{\infty} \theta^{\beta-1} \e^{-s\big(\theta^{\beta}+\frac{\lab\xi\rab^{2}}{2}\big)}ds\\&=\lpa2\pi\rpa^{-\frac d2}
\frac{\theta^{\beta-1}}{\theta^{\beta}+\frac{\lab\xi\rab^{2}}{2}},
\end{split}
\eeq
where $\tilde{}$ and $\hat{}$ denote the Laplace transform in time and the Fourier transform in space, respectively.
Taking the inverse Laplace transform we get that the Fourier transform of $\Kbetatx$ is given by
\beq\lbl{Fbeta}
\FKbetatxi =\lpa2\pi\rpa^{-\frac d2}E_{\beta}\Big(-\frac{\lab\xi \rab^{2}}{2}t^{\beta}\Big).
\eeq
The proof is complete.
\epf

We end this section with three facts about the Mittag-Leffler function $E_{\beta}(x)$ that will be applied in Section \ref{sec:tfspidepfs}.  The first two give upper and lower bounds as well as  the asymptotic behavior for these important functions (see  \cite[Theorem 4 and equation (6.6)]{TS}).  For $\beta\in(0,1)$ and all $x>0$, we have
\beq\lbl{elemineqML}
\frac1{1+\Gamma(1-\beta)x}\le E_{\beta}(-x)\le \frac1{1+\lbk\Gamma(1+\beta)\rbk^{-1}x};
\eeq
and
\beq\lbl{mlas}
E_{\beta}(-x)\sim x^{-1}\Gamma(1-\beta);\ \mbox{as }x\to\infty.
\eeq
The third Mittag-Leffler property we need relates to its Fourier transform, which we give next.   For clarity and convenience, we will, in the next lemma, use the notations $\La[f(t)](\theta)$ and $\Fr[f(t)](\tau)$ for the $t\mapsto\theta$ Laplace transform and the $t\mapsto\tau$ Fourier transform, respectively.   
\blm[Mittag-Leffler Fourier transform]\lbl{lm:MLFT} 
Assume that $x>0$. 
\ben\rencomrom
\item If $\beta>0$, then for any $\sigma>0$
\beq\lbl{eq:FrMLg}
\Fr\lbk\ind_{ \{ t > 0 \}}E_{\beta}\lpa-x t^{\beta}\rpa\e^{-t\sigma}\rbk(\tau)=\frac{(\sigma+\i\tau)^{\beta-1}}{(\sigma+\i\tau)^{\beta}+x}.
\eeq
\item If $\beta\in\lbr1/2^{k};k\in\N\rbr$, then \eqref{eq:FrMLg} hold for any $\sigma\ge0$.  In particular,
\beq\lbl{eq:FrML2k} 
\Fr\lbk\ind_{ \{ t > 0 \}}E_{\beta}\lpa-x t^{\beta}\rpa\rbk(\tau)=\frac{(\i\tau)^{\beta-1}}{(\i\tau)^{\beta}+x}.
\eeq
\een
\elm
\brm\lbl{rm:2kcap012}
The case $\beta\in\lbr1/2^{k};k\in\N\rbr$ is important and useful since it captures the behavior of our SPIDEs for all $0<\beta\le1/2$ while also representing the case where the kernels $\Kbetatx$ is the fundamental solution to higher order PDEs with memory (see e.g., \cite{Atfhosie,Abtp1,Abtp2} and the references therein for details).
\erm
\bpf
Let $\theta=\sigma+\i\tau$, and suppose $x>0$.  Then, 
\beq\lbl{eq:LapFrML}
\bsp
\frac{(\sigma+\i\tau)^{\beta-1}}{(\sigma+\i\tau)^{\beta}+x}&=\La\lbk E_{\beta}\lpa-x t^{\beta}\rpa\rbk(\theta)
\\&=\int_{0}^{\infty}E_{\beta}\lpa-x t^{\beta}\rpa\e^{-t\theta}dt
\\&=\int_{-\infty}^{\infty}\lbk\ind_{ \{ t > 0 \}}E_{\beta}\lpa-x t^{\beta}\rpa\e^{-t\sigma}\rbk\e^{-\i t\tau}dt
\\&=\Fr\lbk\ind_{ \{ t > 0 \}}E_{\beta}\lpa-x t^{\beta}\rpa\e^{-t\sigma}\rbk(\tau)
\end{split}
\eeq
If $\beta\in\lbr1/2^{k};k\in\N\rbr$, there are no poles if we set $\sigma=0$ in the ratio 
\beqs\lbl{eq:LaML}
\frac{(\sigma+\i\tau)^{\beta-1}}{(\sigma+\i\tau)^{\beta}+x}.
\eeqs
In this case, the radius of convergence of the Laplace transform in      
\eqref{eq:LapFrML} is $\Re(\theta)\ge0$.  Setting $\sigma=0$ in \eqref{eq:LapFrML}, we thus obtain the Mittag-Leffler Fourier transform in \eqref{eq:FrML2k}.  
The proof is complete.  
\epf

\section{The L-KS SPDEs: Proofs of Theorems 1.1--1.3}
Let $U$ be the solution to the L-KS SPDE \eqref{nlks}.  In \cite{Alksspde}, Allouba
obtained the temporal and spatial H\"older exponent $\gamma_{t}\in(0,(4-d)/8)$ and $\gamma_{s}\in(0,((4-d)/2)\wedge1)$, respectively, by establishing---in
\cite[Lemma 3.4]{Alksspde}---the following sharp dimension-dependent upper bounds
\begin{equation}\lbl{ublkssol}
\bc
\E \lbk U (t,x) - U(s,x)\rbk^{2q}&\le {C}_{d}\lab t-s\rab^{\tf{(4-d)q}{4}},\\
\E \lbk U (t,x) - U(t,y)\rbk^{2q}&\le {C}_{d}|x-y|^{2q\alpha_{d}};\quad  \alpha_{d}\in J_{d},\end{cases}
\end{equation}
for the more general nonlinear L-KS SPDE \eqref{nlksa}, with Lipschitz condition on $a$, for all $x,y \in \Rd$,
for all $t,s \in[0,T]$, for $q\ge1$, for $1\le d\le 3$, and for
\beq\lbl{intervals}
J_{d}=\bc (0,1];&d=1,\\
(0,1);&d=2,\\
(0,\frac12);&d=3.
\ec
\eeq

These H\"older exponents determine the temporal and spatial differences exponents
in the temporal and spatial moduli expressions in \thmref{tempmodlks} and \thmref{spatmodlks},
respectively.  They are also useful for getting a sharp upper bound for the uniform
spatio-temporal moduli of continuity for our L-KS SPDE \eqref{nlks}.  Rather than
complementing the upper bounds in \eqref{ublkssol} with corresponding lower bounds,
we take a harmonic/spectral analytic route combined with a useful decomposition of
our solution $U$ to get the exact uniform and local moduli of continuity in \thmref{tempmodlks}.
This approach, which we also use for our time-fractional SPIDEs, builds on
the results of Xiao in \cite{X07,X09} and  Meerschaert, Wang, and Xiao in \cite{MWX}.

Assume without loss of generality that $u_{0} = 0$, then the L-KS SPDE solution is given by
\begin{equation}\lbl{zeroinit}
U(t,x) = \int_{\mathbb R^{d}} \int_{0}^{t} \KKSepthtrxy\sW(dr\times dy), \quad t \ge 0, x \in \Rd.
\end{equation}

\subsection{Temporal modulus}\lbl{tmpmodlkssec}
Throughout this subsection, let $x\in\mathbb R^{d}$ be fixed but arbitrary.  We first introduce
the following auxiliary Gaussian process $\{ X(t,x), t \in \mathbb{R}_{+}\}$:
\begin{equation}\lbl{auxgauss}
X(t,x) = \int_{\mathbb{R}^{d}} \int_{\mathbb{R}} \lpa\KKSepthtrposxy-\KKSepthmrposxy \rpa \sW(dr \times dy),
\end{equation}
where $a_{+} = \max \{ a,0 \}$ for all $a \in \mathbb{R}$. Then the L-KS SPDE solution $U$ may be decomposed
as  $U(t,x) = X(t,x) - V(t,x) $, where
\begin{equation}
V(t,x) = \int_{\mathbb{R}^{d}} \int_{\mathbb{R}_{-}} \lpa\KKSepthtrposxy-\KKSepthmrposxy \rpa \sW(dr \times dy).
\end{equation}
This idea of decomposition originated in Mueller and Tribe \cite{MuTr02} in the
second order SPDEs setting; and it has been applied in Wu and Xiao \cite{WX06} and in Tudor and Xiao \cite{TuX15},
also in the second order heat SPDE setting. See also Mueller and Wu \cite{MW12} for related results on stochastic
heat equation.

We first prove our results on the moduli of continuity for the auxiliary process $X$, then using the aforementioned
decomposition of $U$, in terms of $X$ and a smooth process $V$, we transfer them to our L-KS SPDE solution $U$.
The following result is pivotal.
\bthm\lbl{auxdecompspden}  Assume the spatial dimension $d\in\{1,2,3\}$. Let $X$ be as
defined in \eqref{auxgauss} and $x \in \R$ be fixed.
\begin{enumerate}
\renewcommand{\labelenumi}{(\roman{enumi})}
\item The Gaussian process $\lbr X(t,x); t\ge 0\rbr$ has stationary temporal increments.  Moreover, we have
\begin{equation*}
\mathbb{E} \lbk X(t,x) - X(s,x)\rbk^{2} = 2 \int_{\mathbb{R}} [1-\cos((t-s) \tau)] \Delta(\tau) d\tau,
\end{equation*}
where the spectral density $\Delta$ is given by
\begin{equation*}
\Delta(\tau) = (2\pi)^{-d} \int_{\mathbb{R}^{d}} \frac{d\xi}{\tau^{2} + \frac{\varepsilon^{2}}{64}(-2\theta +|\xi|^{2})^{4}}.
\end{equation*}
\item 
%
For each $k \ge 1$, there exists a modification of $\{V(t,x), t\in \mathbb{R}_{+} \}$ such that its (temporal)
sample function is almost surely continuously $k$-times differentiable on $(0,\infty)$.
\item Let $\gamma_1 = \frac{4-d}{8}$. There is a finite constant $C$ such that
\begin{equation}\label{Eq:modY}
\lim_{\varepsilon \to 0} \sup _{s, t\in [0,\, \varepsilon]} \frac{\left|V(t,x) - V(s,x) \right|}
{|t-s|^{\gamma_1} \sqrt{\log \log (1/|t-s|)} } \leq C  \ \  \mbox{ a.s. }
\end{equation}
\end{enumerate}
\ethm
\bpf
To verify (i), we apply Parseval's identity to the integral in $y$ to get that for any $0 < s < t$:
\begin{equation}\lbl{tempdiffL2}
\begin{split}
\mathbb{E} & \lbk X(t,x) - X(s,x)\rbk^{2}= \int_{\mathbb{R}^{d}} \int_{\mathbb{R}} \lab\KKSepthtrposxy-\KKSepthsrposxy \rab^{2} dr dy \\
& = \int_{\R}  \int_{\Rd}  \lab\FKKSepthtrposxxi-\FKKSepthsrposxxi \rab^{2} d\xi dr.
\end{split}
\end{equation}
Since
\begin{equation}\lbl{tempom}
\FKKSepthtrposxxi = (2\pi)^{-d/2} \cdot \e^{-\i \langle x,\xi\rangle - \frac{\varepsilon (t-r)}{8} (-2\theta + |\xi|^{2})^{2}} \ind_{ \{t > r\} },
\end{equation}
equation \eqref{tempdiffL2} becomes
\begin{equation}
\begin{split}
&\E \lbk X(t,x) - X(s,x)\rbk^{2}\\ & =  \int_{\mathbb{R}^{d}}  \int_{\mathbb{R}}  \df{\lab \e^{- \frac{\varepsilon (t -r)}
{8} (-2 \theta + | \xi |^{2} )^{2}} \ind_{ \{ t > r \}}
- \e^{- \frac{\vep (s -r)}{8} ( -2 \theta +| \xi |^{2} )^{2}} \ind_{ \{ s > r \}} \rab^{2}}{(2\pi)^{d}} dr d\xi.
\end{split}
\end{equation}
Now, we apply Parseval's identity to the inner integral in $r$. To this end, let
\begin{equation*}
\phi (r,\xi) = \e^{-\frac{\varepsilon (t-r)}{8} (-2\theta + | \xi |^{2})^{2}} \ind_{\{ t > r \} }
 - \e^{-\frac{\varepsilon(s - r)}{8} (-2\theta + |\xi |^{2} )^{2}} \ind_{ \{ s > r \}}.
 \end{equation*}
 Its Fourier transform in $r$ is
 \begin{equation*}
 \widehat{\phi} (\tau,\xi) = (\e^{\i\tau t} - \e^{\i\tau s} ) \frac1{\i\tau + \frac{\varepsilon}{8} (-2\theta + | \xi |^{2})^{2}}.
 \end{equation*}
Hence, by Parseval's identity, we get
\begin{equation}
\begin{split}
\E & \lbk X(t,x) - X(s,x)\rbk^{2} = (2\pi)^{-d} \int_{\mathbb{R}^{d}} \int_{\mathbb{R}} \big| \widehat{\phi} (\tau,\xi) \big|^{2} d\tau d\xi\\
& = 2(2\pi)^{-d}  \int_{\mathbb{R}} (1 - \cos ((t-s)\tau)) \int_{\mathbb{R}^{d}} \frac{d\xi}{\tau^{2} + \frac{\varepsilon^{2}}{64} (-2\theta
+ | \xi|^{2} )^{4}} d\tau.
\end{split}
\end{equation}
The proof of (i) is complete.

The proof of part (ii) is similar to \cite[Theorem 4.8]{XX11}, but is more complicated in our higher order case and its corresponding kernel.
For completeness, we give the main steps of the proof.

We start with the case $k=1$.  The mean square derivative of $V$ at $t \in (0, \infty)$ is given by
\[
\pat V(t,x) = \int_{\mathbb{R}^{d}} \int_{\mathbb{R}_{-}} \pat\KKSepthtrxy \sW(dr\times dy).
\]
This can be verified by checking the covariance function.  For every $s, t \in (0, \infty)$ with $s\leq t$ we have
\beq\label{Eq:deriva}
\bsp
&\mathbb{E}\lab\pat V(t,x)-\pa_{s} V(s,x) \rab ^{2}
\\& =\mathbb{E} \lab\int_{\mathbb{R} ^{d}} \int_{-\infty} ^{0}  \left( \pa_{t}\KKSepthtrxy -\pa_{s}\KKSepthsrxy \right) \sW(dr\times dy)\rab ^{2}
\\& = \int_{-\infty} ^{0}  \int_{\Rd}  \lab\pat\FKKSepthtrxxi-\pa_{s}\FKKSepthsrxxi \rab^{2} d\xi dr
\\& =C \int_{\mathbb{R}^d} (-2\theta + |\xi|^{2})^{4}  \int_{0} ^{\infty}
\lab \e^{- \frac{\varepsilon (t+r)}{8} (-2\theta + |\xi|^{2})^{2}}- \e^{- \frac{\varepsilon (s+r)}{8} (-2\theta + |\xi|^{2})^{2}}\rab^{2}dr d\xi,
\end{split}
\eeq
where we have used Parseval's identity to the integral in $y$ and the fact that the Fourier transform of the function
$y \mapsto \pa_{t}
\KKSepthtpry$ is
\[
\pat
\FKKSepthtprxi= -\frac{\varepsilon(-2\theta + |\xi|^{2})^{2}}{8(2\pi)^{d/2}}
\e^{- \frac{\varepsilon (t+r)}{8} (-2\theta + |\xi|^{2})^{2}}.
\]
Let $\psi(r,\xi) = \big( \e^{- \frac{\varepsilon (t+r)}{8} (-2\theta + |\xi|^{2})^{2}}- \e^{- \frac{\varepsilon (s+r)}{8} (-2\theta + |\xi|^{2})^{2}}\big)
\ind_{ \{r > 0\} }$.
Then, its Fourier transform in $r$ is given by
\beqs
\bsp
\widehat{\psi} (\tau,\xi)&=\lpa\e^{- \frac{\varepsilon }{8} t(-2\theta + |\xi|^{2})^{2}}- \e^{- \frac{\varepsilon }{8} s(-2\theta + |\xi|^{2})^{2}}\rpa
\int _0 ^\infty  \e^{-r\lbk\frac{\varepsilon}{8}(-2\theta + |\xi|^{2})^{2} +i\tau\rbk} dr
\\&= \lpa\e^{- \frac{\varepsilon }{8} t(-2\theta + |\xi|^{2})^{2}}- \e^{- \frac{\varepsilon }{8} s(-2\theta + |\xi|^{2})^{2}}\rpa \frac1{\i\tau +
\frac{\varepsilon}{8} (-2\theta + | \xi |^{2})^{2}}.
\end{split}
\eeqs
Thus, for any $0 < a < b < \infty$, we see that for each $s, t \in [a,b]$ with $s < t$ equation \eqref{Eq:deriva} becomes
\beq\label{Eq:deriva2}
\bsp
&\mathbb{E}\lab\pat V(t,x)-\pa_{s} V(s,x) \rab ^{2}
\\ &=  C \int_{\mathbb{R}^d} (-2\theta + |\xi|^{2})^{4}
\lab \e^{- \frac{\varepsilon t}{8} (-2\theta + |\xi|^{2})^{2}}- \e^{- \frac{\varepsilon s}{8} (-2\theta + |\xi|^{2})^{2}}\rab^{2}
\\&\qquad\times\int_{\R} \frac{1}{\tau^{2} + \frac{\varepsilon^{2}}{64} (-2\theta + | \xi|^{2} )^{4}} d\tau d\xi
\\&= C \int_{\mathbb{R}^d} (-2\theta + |\xi|^{2})^{2} \e^{- \frac{\varepsilon s}{4} (-2\theta + |\xi|^{2})^{2} } \lab1- \e^{- \frac{\varepsilon (t-s)}{8}
(-2\theta + |\xi|^{2})^{2}} \rab^{2} d\xi
\\&\le C\lab t-s\rab^{2}\int_{\mathbb{R}^d} (-2\theta + |\xi|^{2})^{6} \e^{- \frac{\varepsilon a}{4} (-2\theta + |\xi|^{2})^{2} } d\xi
\le C\lab t-s\rab^{2}.
\end{split}
\eeq
Now, by using Kolmogorov's continuity theorem, we can find a modification of $V$ such that $V(t,x)$ is continuously differentiable on
in $t$ on $[a, b]$ (see e.g. \cite{XX11}).
This proves (ii) for $k =1$.  For $k =2$, we apply the above argument to the Gaussian
process $\{\pat^{2}V(t,x), t \ge 0\}$, where, for each $t> 0$, $\pat^{2}V$ is the
second order mean-square derivative,  and we
find a modification of $V$ whose temporal sample paths are twice continuously
differentiable on $[a, b]$. Iterating this procedure finishes the proof of (ii).

To prove (iii), we will apply the metric-entropy method (cf. e.g., \cite{MarRos}). It can be verified that
$\E \big[V^2(t,x)\big] \asymp t^{2 \gamma_1}$  for  $ t \in [0, 1]$ and that $\E \big[V^2(t,x)\big] \sim C t^{2 \gamma_1}$
as $|t|\searrow0$.  Recall that $\gamma_1 = (4-d)/8$. For any $0 < s < t$, we proceed
similarly to part (ii) above to get
\beq\label{Eq:regatzero}
\bsp
&\mathbb{E}\lab V(t,x)- V(s,x) \rab ^{2}
\\ &=  C \int_{\mathbb{R}^d}
\lab \e^{- \frac{\varepsilon t}{8} (-2\theta + |\xi|^{2})^{2}}- \e^{- \frac{\varepsilon s}{8} (-2\theta + |\xi|^{2})^{2}}\rab^{2}
\\&\qquad\times\int_{\R} \frac{1}{\tau^{2} + \frac{\varepsilon^{2}}{64} (-2\theta + | \xi|^{2} )^{4}} d\tau d\xi
\\&= C \int_{\mathbb{R}^d} \e^{- \frac{\varepsilon s}{4} (-2\theta + |\xi|^{2})^{2} } \frac{\lab1- \e^{- \frac{\varepsilon (t-s)}{8}
(-2\theta + |\xi|^{2})^{2}} \rab^{2}}{(-2\theta + |\xi|^{2})^{2}} d\xi
\\&\le C\lab t-s\rab^{2}\int_{\mathbb{R}^d} (-2\theta + |\xi|^{2})^{2} \e^{- \frac{\varepsilon s}{4} (-2\theta + |\xi|^{2})^{2} } d\xi
\\&\le C s^{2\gamma_1-2}\lab t-s\rab^{2}.
\end{split}
\eeq
Thus, the canonical metric of $V$ is given by
\[
d_V(s, t) = \sqrt{\mathbb E |V(t,x)-V(s,x) |^{2}} \le C\left\{\begin{array}{cc}
t^{\gamma_1}\qquad \ \ &\hbox{ if } 0=s < t,\\
s^{ \gamma_1-1} |t-s| &\hbox{ if } 0 < s < t.
\end{array}
\right.
\]
It follows from the Gaussian isoperimetric inequality (cf. Lemma 2.1 in \cite{Tal95})
that there is a constant $C \ge 1$ such that for any constant $\varepsilon  > 0$, and $u > 0$,
\[
\mathbb {P}\Big(\max_{0 \le s, t \le \varepsilon}|V(t,x)- V(s,x)| \ge u\Big)
\le C \exp \Big(- \frac{u^2}{C \varepsilon^{2 \gamma_1}}\Big).
\]
A standard Borel-Cantelli argument yields that for some positive and finite constant $C$,
\[
\lim_{\varepsilon \to 0} \sup_{0 \le s, t \le \varepsilon} \frac{|V(t,x) - V(s,x)|}
{|t-s|^{\gamma_1} \sqrt{\log \log (1/|t-s|)}} \le C,
\quad \hbox{ a.s.}
\]
This proves (iii).
\epf

\bcr\lbl{tmplksspdeas}
Assume the spatial dimension $d\in\{1,2,3\}$.  The spectral density $\Delta$ is
asymptotically given by
\begin{equation}
\lbl{sdastimelks}
\Delta(\tau) \sim \frac{(2\pi)^{-d}}{\tau^{2-d/4}} \int_{\mathbb{R}^{d}}
\frac{d\xi}{1+\frac{\varepsilon^{2}}{64} |\xi|^{8}},
\mbox{ as } \tau \to \infty.
\end{equation}
\ecr

By combining the asymptotic behavior of the spectral density $\Delta$ in
\coref{tmplksspdeas} and Theorem 2.1 in Xiao \cite{X07},  we obtain the following
strong local nondeterminism and double-sided bounds for $\E [X(t,x)-X(s,x)]^{2}$.

\bcr[Temporal SLND and double-sided bounds for X]\lbl{SLNDtlks}
Assume the spatial dimension $d\in\{1,2,3\}$.  For any $T>0$, there is a
positive constant $c$ such that for all $t\in(0,T]$
and all $0 < r \le 1\wedge|t|$
\beq\lbl{SLNDtlkseq}
\var \lpa X(t,x) | X(s,x);s\in[0,T],|t-s|\ge r\rpa \ge c\, r^{\frac{4-d}{4}}
\eeq
Also,
\beq\lbl{dstemplks}
\E [X(t,x)-X(s,x)]^{2} \asymp|t-s|^{\frac{4-d}{4}};\ \ \ \forall s,t\in[0,T].
\eeq
Here and on the sequel, the notation $f\asymp g$ on $S$ means
$c_{l}g(x)\le f(x)\le c_{u}g(x)$ for all $x\in S$ for some constants $c_{l},c_{u}$.
\ecr
\bpf
\coref{tmplksspdeas} and Theorem 2.1 in \cite{X07} imply \eqref{SLNDtlkseq}.
\coref{tmplksspdeas} and Theorem 2.5 in \cite{X07} imply \eqref{dstemplks}
\epf
From \coref{SLNDtlks}, the Gaussian process $\{X(t, x), t \ge 0\}$ satisfies
conditions (C1) and (C2) above ((A1) and (A2) in \cite{MWX}). Hence we can apply
the results in \cite{MWX} on the uniform and local moduli of continuity to get
the following theorem on the time regularity of $X$.
\begin{thm}\lbl{auxmod}
Let $x\in \mathbb{R}^{d}$, $d=1,2,3$, be fixed. Let $\{ X(t,x), t \ge 0 \}$ be
defined as above and let $\gamma_1 = \frac{4-d}8$.  Then,
\ben\rencomrom
\item (Uniform Modulus of Continuity)  for every compact interval $\tint\subseteq \mathbb{R}_{+}$
\begin{equation}\lbl{umctaux}
\lim_{\delta\searrow 0} \sup_{\substack{|s-t| < \delta \\ s, t \in I} } \frac{|X(t,x) - X(s,x)|}{|s-t|^{\gamma_1}\sqrt{\log \frac1{|t-s|}}} = k_{1}^{(d)}; \mbox{ a.s.}
\end{equation}
\item (Local Modulus of Continuity) and for every fixed $t \ge 0$
\begin{equation}\lbl{lmctaux}
\lim_{\delta\searrow 0 } \frac{\sup_{|s-t| <\delta} | X(s,x) - X(t,x)|}{\delta^{\gamma_1} \sqrt{\log\log \frac1{\delta}}} = k_{2}^{(d)}; \mbox{ a.s.}
\end{equation}
\een
In the above, $0<k_{i}^{(d)}<\infty$ ($i = 1, 2$) are $d$-dependent constants, independent of $x\in\Rd$.
\end{thm}
\bpf The uniform modulus of continuity of $X$ in \eqref{umctaux} follows from Theorem 4.1 in \cite{MWX}.  The local modulus
of continuity of $X$ in \eqref{lmctaux} follows upon applying \cite[Theorem 5.1]{MWX}.
The constants $k_{i}^{(d)}$ ($i = 1, 2$) in \thmref{auxmod} do not depend on $x\in\Rd$ since the distribution of the process
$\{X(t,x), t \ge 0\}$ does not depend on $x$, see \eqref{sdastimelks}.
\epf

We believe that $k_{1}^{(d)} = k_{2}^{(d)}$ because the large deviation behavior of
the tail probabilities of the maxima
$\sup_{s, t \in [0, b], \vert s-t\vert \leq \varepsilon }  \left| X(t,x)- X(s, x) \right| $
and, for fixed $t$,
$\sup _{\vert s-t\vert \leq \varepsilon}  \left| X(s,x)- X(t, x) \right| $ are the same.
However, the method in \cite{MWX} is not
enough for proving $k_{1}^{(d)} = k_{2}^{(d)}$, a different argument may be needed.

We are now ready to use the the decomposition $U(t, x) = X(t, x) - V(t, x)$ ($t\ge 0$),
 Theorems 3.1 and 3.2  to prove part (i) of  \thmref{tempmodlks}.

\bpfs{Proof of \thmref{tempmodlks} (i)} In order to derive the temporal uniform modulus
of continuity for our L-KS SPDE solution process $U$, we use  part (iii) of
\thmref{auxdecompspden} to see that, almost surely, there exists $\varepsilon_0 > 0$
such that for all $0 < \varepsilon \le \varepsilon_0$,
\begin{equation}\label{Eq:o}
\sup _{s, t\in [0,\, \varepsilon]} \frac{\left|V(t,x) - V(s,x) \right|}
{|t-s|^{\gamma_1} \sqrt{\log \log (1/|t-s|)} } \leq C.
\end{equation}
By splitting the interval $I_{time} = [0, b] \subset \R_+$ into $[0, \varepsilon_0]\cup 
[\varepsilon_0, b]$ and applying (\ref{Eq:o}) and part (ii) of \thmref{auxdecompspden}
respectively, one can see that the solution process $U$ and $\{X(t,x), t \ge 0\}$ have 
the same exact uniform modulus of continuity on  $I_{time}$. Hence, it follows from 
Theorem 3.2 (i) that \thmref{tempmodlks} (i) (a) holds almost surely.

To prove \thmref{tempmodlks}  (i) (b),  we see that, for any $t > 0$, (\ref{Eq:time-LIL}) 
follows from part (ii) of \thmref{auxdecompspden} and part (ii) of \thmref{auxmod}.   
When $t = 0$,  \thmref{auxmod} does not imply   (\ref{Eq:time-LIL}) because the local 
oscillation $V(t, x)$ at the origin may be of the same order. We can prove (\ref{Eq:time-LIL}) 
for $t=0$ by using the comparison result in Lemma 7.1.10 and Remark 7.1.11 in \cite{MarRos}. 
Since this is very similar to the proof of Proposition 2 in \cite{TuX15}, we omit
the details. This finishes the proof of \thmref{tempmodlks} part (i).
\epfs
\subsection{The bifractional Brownian motion link: the case $\vth=0$}  We now turn 
to proof of the L-KS SPDE
bifractional Brownian morion link.

\bpfs{Proof of part (ii) of \thmref{tempmodlks}}
Using Parseval's identity to compute
the covariance function of $U$, we get
\begin{equation}
\begin{split}
\mathbb{E} & [U(t,x) U(s, x)] = \int_{\mathbb{R}^{d}} \int_{0}^{s} 
\KKSepthtrxy  \KKSepthsrxy dr dy \\
& =\int_{0}^{s}  \int_{\mathbb{R}^{d}} \FKKSepthtrxxi  \overline{\FKKSepthsrxxi}d\xi dr\\
& = (2\pi)^{-d} \int_{0}^{s}  \int_{\mathbb{R}^{d}} \e^{-\frac{\varepsilon(t-r)}{8} 
(-2\theta + |\xi|^{2} )^{2} - \frac{\varepsilon(s-r)}{8}
(-2\theta + |\xi|^{2} )^{2}} d\xi dr\\
& = (2\pi)^{-d} \int_{0}^{s} \int_{\mathbb{R}^{d}} \e^{-\frac{\varepsilon (t+s-2r)}{8} 
(-2\theta + |\xi|^{2} )^{2}} d\xi dr.
\end{split}
\end{equation}
When $\theta = 0$,  the above becomes:
\begin{equation}
\begin{split}
\mathbb{E} & \lbk U(t,x)U(s,x)\rbk = (2\pi)^{-d} \int_{0}^{s}  \int_{\mathbb{R}^{d}} \e^{-\frac{\varepsilon(t+s-2r)}{8} \cdot |\xi|^{4}} d\xi dr\\
& = \lbk(2\pi)^{-d} \Big(\frac{8}{\varepsilon}\Big)^{d/4} \frac1{2-d/2} 
\int_{\mathbb{R}^{d}} \e^{-|\xi|^{4}} d\xi\rbk \lbk(t+s)^{1-\frac{d}{4}} - (t-s)^{1-\frac{d}{4}}\rbk.
\end{split}
\end{equation}
Hence, up to a constant, the mean zero Gaussian process $\lbr U(t,x), t\ge 0 
\rbr \ (x \in \mathbb{R}^{d} \mbox{ fixed})$ is a bifractional Brownian motion 
with indices $H=\frac12$ and $K= 1 - \frac{d}4$.  More precisely, $U(\cdot,x)\eqL c_{d}B^{\lpa\frac12,\frac{4-d}{4}\rpa}$, where
\beq
c_{d}=(2\pi)^{-d/2} \Big(\frac{8}{\varepsilon}\Big)^{d/8} \frac{2^{(d-4)/8}}{\sqrt{2-d/2}}\sqrt{\int_{\mathbb{R}^{d}} \e^{-|\xi|^{4}} d\xi}.
\eeq
Hence many sample path properties of  $\lbr U(t,x), t\ge 0 \rbr$, including Chung's 
law of the iterated logarithm in \eqref{chunglil},
can be derived from Tudor and Xiao \cite{TuX07} directly\footnote{When $\theta \in 
\mathbb{R} \backslash \{0\}$, it is not as
simple to obtain an explicit expression in terms of $s$ and $t$. }.\epfs
\subsection{Spatial modulus}
Recall our standing assumption that $\un=0$, and the solution is given by 
\eqref{zeroinit}. Let $t > 0$ be fixed, we consider the L-KS Gaussian random 
field $\{U(t,x), x \in \mathbb{R}^{d} \}$.
Our results are based on the following lemma.
\begin{lem}[L-KS SPDE spatial spectral density]\lbl{spdenlksspatlem}
Assume the spatial dimension $d\in\{1,2,3\}$.  The centered Gaussian random field $\{ U(t,x), x \in \mathbb{R}^{d}\}$
is stationary with spectral density
\begin{equation*}
S(\xi) = \frac{4}{\varepsilon(2\pi)^{d}} \cdot \frac{1-\e^{-\frac{\varepsilon t}{4}
(-2\theta + |\xi|^{2})^{2}}}
{(-2\theta + |\xi|^{2})^{2}}, \quad \forall \xi \in \mathbb{R}^{d}.
\end{equation*}
\end{lem}
\bpf Using Parseval's identity, we compute the covariance
\begin{equation}
\begin{split}
\mathbb{E} & \lbk U(t,x) U(t,y)\rbk = \int_{0}^{t} dr \int_{\mathbb{R}^{d}} 
\KKSepthtrxz \KKSepthtryz dz \\
& = \int_{0}^{t} dr \int_{\mathbb{R}^{d}}  \FKKSepthtrxxi\overline{\FKKSepthtryxi}d\xi \\
& = (2\pi)^{-d} \int_{0}^{t} dr \int_{\mathbb{R}^{d}} \e^{\i \langle\xi, x-y\rangle} \cdot \e^{-\frac{\varepsilon(t-r)}{4} (-2\theta + |\xi|^{2})^{2}} d\xi\\
& = (2\pi)^{-d} \int_{\mathbb{R}^{d}} \e^{\i \lang \xi, x-y\rang} d\xi
\cdot \int_{0}^{t} \e^{- \frac{\varepsilon (t-r)}{4} (-2\theta + |\xi|^{2} )^{2}} dr \\
& = \frac4{\varepsilon}(2\pi)^{-d}  \int_{\mathbb{R}^{d}} \e^{\i \lang \xi, x-y\rang} \frac{1-\e^{-\frac{\varepsilon t}{4} (-2\theta + |\xi|^{2} )^{2}}}
{(-2\theta + |\xi|^{2} )^{2}} d\xi.
\end{split}
\end{equation}
Thus, the conclusions of the Lemma follows.
\epf
As an immediate consequence, we get
\bcr\lbl{lkssdas}  Assume the spatial dimension $d\in\{1,2,3\}$.
The spectral density $S$ satisfies $0 < S(0) < \infty$ and has
the asymptotic behavior
\begin{equation*}
S(\xi) \sim \frac{C_{\varepsilon, \theta}}{|\xi|^{d+2\gamma_2}},
\end{equation*}
as $|\xi| \to \infty$, where $\gamma_2 = 2-\frac{d}{2}$.
\ecr

\subsubsection{The case $d=1$: gradient spatial H\"older and modulus of continuity}
We now complete the proof of the one dimensional case in \thmref{spatmodlks} (i).
\bpfs{Proof of \thmref{spatmodlks} (i)}
Fix $t>0$.  We start with the H\"older assertion for the gradient.  By applying
\lemref{spdenlksspatlem}, we can show that the mean square gradient $\pa_{x} U(t,x)$ exists and
\beq\lbl{spderlkseq}
\bsp
&\E\Big|\pa_{x} U(t,x) -\pa_{y}U(t,y) \Big|^2
=\int_{\R} \xi^2 \big|\e^{\i x\xi} - \e^{\i y\xi}\big|^2 S(\xi)d\xi\\
&=\frac{4}{\varepsilon \pi }\int_{\R} \xi^2 \lbk1- \cos(\xi(x-y))\rbk \frac{ 1-\e^{-\frac{\varepsilon t}{4} (-2\theta + |\xi|^{2} )^{2}}}
{(-2\theta + |\xi|^{2} )^{2}} d\xi. 
\end{split}
\eeq
Assume, without loss of generality, that $\vep=\vth=1$.  We also assume 
that $|x-y| \le 1/2$. Proceeding as in the proof of
\cite[Lemma 3.3]{Alksspde}, we split the last integral over three sets 
$\B_{1}:=\lbr\xi\in\R: \lab\xi\rab< 2\rbr$,  $\B_{2}:=
\big\{\xi\in\R: 2\le \lab\xi\rab< \frac 1 {|x-y|}\big\}$ and $\B_{3}
:=\big\{\xi\in\R: \lab\xi\rab \ge \frac 1 {|x-y|}\big\}$. We will
make use of the following elementary inequalities:
\beq\lbl{elemineq}
\bsp
&\mbox{(a)}\quad1-\cos z \le 2 \wedge  z^2 ,
\\&\mbox{(b)}\quad \frac{1-\e^{-\frac{t}{4} (-2+\lab\xi\rab^{2})^{2}}}
{(-2+\lab\xi\rab^{2})^{2}}
\le \left\{\begin{array} {ll}
\tfrac t4, \ \ \ \ &\hbox{ on } \B_1,\\
\frac{c}{|\xi|^4},\ \  &\hbox{ on } \B_2 \cup \B_3.
\end{array}
\right.
\end{split}
\eeq
It follows from \eqref{spderlkseq} and \eqref{elemineq} that
\beq\lbl{spderlkseq2}
\bsp
&\E\Big|\pa_{x} U(t,x) -\pa_{y}U(t,y) \Big|^2 \\
&\le C \lbk\int_{\B_1} \xi^4 |x-y|^2  d\xi +  \int_{\B_2}  |x-y|^2  \frac{1}{|\xi|^2} d\xi
+ \int_{\B_{3}} \frac 1 {\xi^2}  d\xi\rbk
\\&\le C |x-y|.
\end{split}
\eeq
Thus, Kolmogorov's continuity theorem gives us the spatial local $\gamma$-H\"older
continuity for the L-KS gradient, $\pa_{x}U$, for $\gamma\in(0,1/2)$.

Turning now to the exact uniform and local spatial continuity moduli for the L-KS gradient,
$\pa_{x}U$, in \thmref{spatmodlks} (i)
(a) and (b). We first compute the gradient covariance as follows
\beq\lbl{lksgradcov}
\bsp
\mathbb{E} & \lbk \pa_{x}U(t,x)\pa_{y} U(t,y)\rbk = \int_{0}^{t}  \int_{\mathbb{R}}  \xi^{2}\FKKSepthtrxxi\overline{\FKKSepthtryxi}d\xi dr\\
& = (2\pi)^{-d} \int_{0}^{t}  \int_{\mathbb{R}}\xi^{2} \e^{\i \langle\xi, x-y\rangle} \cdot \e^{-\frac{\varepsilon(t-r)}{4} (-2\theta + |\xi|^{2})^{2}} d\xi dr\\
& = \frac4{\varepsilon}(2\pi)^{-d}  \int_{\mathbb{R}} \e^{\i \lang \xi, x-y\rang} \frac{\xi^{2}[1-\e^{-\frac{\varepsilon t}{4} (-2\theta + |\xi|^{2} )^{2}}]}
{(-2\theta + |\xi|^{2} )^{2}} d\xi.
\end{split}
\eeq
This means that the spatial spectral density of $\pa_{x}U$, denoted by $f$, and
its asymptotic behavior are given by
\beq\lbl{lksgradspdensp}
f(\xi)=\frac4{\varepsilon}(2\pi)^{-d}\frac{\xi^{2}[1-\e^{-\frac{\varepsilon t}{4}
(-2\theta + |\xi|^{2} )^{2}}]}{(-2\theta + |\xi|^{2} )^{2}} d\xi\sim\frac{C}{|\xi|^{2}},
\quad \hbox{ as } |\xi| \to \infty.
\eeq
Equation \eqref{lksgradspdensp} and Theorem 2.1 in \cite{X07} imply that, for every fixed
$t>0$, the gradient of the L-KS SPDE solution $\{\pa_{x}U(t,x), x \in \R\}$ is spatially
strongly locally nondeterministic. More precisely, for every $M>0$, there exists a finite
constant $c>0$ (depending on $t$ and $M$) such that for every $n\geq 1$ and
for every $ x,y_{1},..., y_{n} \in [-M,M]$,
\begin{equation}\lbl{SLNDlksgradsp}
{\rm Var}  \left[\pa_{x}U(t,x) |\pa_{x}U(t, y_{1}),\ldots,\pa_{x} U(t, y_{n} ) \right]
\geq c\min _{ 0\leq j\leq n} \{ \vert x -y_{j} \vert \},
\end{equation}
where $y_{0}=0$.  Also,  \eqref{lksgradspdensp} and Theorem 2.5 in \cite{X07} imply
the double sided second moment bounds
\begin{equation}\lbl{dsblksgradsp}
\E [\pa_{x}U(t,x)-\pa_{y}U(t,y)]^{2} \asymp|x-y|;\quad  \forall x, y\in[-M,M].
\end{equation}
Thus, the uniform modulus of continuity of $\pa_{x}U$ in \thmref{spatmodlks} (i) (a)
follows from Theorem 4.1 in \cite{MWX}.
The local modulus of continuity of $\pa_{x}U$ in \thmref{spatmodlks} (i) (b)
follows upon applying \cite[Theorem 5.1]{MWX}.
\epfs

\subsubsection{The fractal cases $d=2,3$}  We now turn to the rougher two and
three dimensional cases.  Starting with the $d=3$ case, we first obtain the
strong local nondeterminism property and double-sided second moment bounds in
space for the L-KS SPDE solution $\lbr U(t,x); x \in \mathbb{R}^{3}\rbr$.

\blm[Spatial SLND and double-sided bounds for L-KS SPDEs] \label{SLNDslks}
For every fixed $t>0$, the L-KS SPDE solution $\{U(t,x), x \in \R^{3} \}$ is
spatially strongly locally nondeterministic.
Namely, for every $M>0$, there exists a finite constant $c>0$ (depending on
$t$ and $M$) such that for every $n\geq 1$ and
for every $ x, y_{1},..., y_{n} \in [-M,M] ^{3}$,
\begin{equation}\lbl{SLNDslkseq}
{\rm Var}  \left[ U(t,x) | U(t, y_{1}),\ldots, U(t, y_{n} ) \right] \geq
c\min _{ 0\leq j\leq n} \{ \vert x -y_{j} \vert \},
\end{equation}
where $y_{0}=0$.  Also,
\begin{equation}\lbl{dsbslks}
\E [U(t,x)-U(t,y)]^{2} \asymp|x-y|;\quad \forall x,y\in[-M,M]^{3}.
\end{equation}
\elm
\bpf
When $d = 3$, \coref{lkssdas} implies that the condition (2.17) in \cite{X07}
is satisfied with $\underline{\alpha} = \overline{\alpha} = \gamma_2
= \frac 1 2$. Hence, the conclusions in \eqref{SLNDslkseq} and \eqref{dsbslks}
follow from Theorem 2.5 in \cite{X07} with $\phi(r) = r$.
\epf

Now we can obtain the exact spatial uniform and local continuity moduli in
\thmref{spatmodlks} (ii) for the three dimensional case.
\bpfs{Proof of \thmref{spatmodlks} (ii)}  With \lemref{SLNDslks} in hand, 
the uniform modulus of continuity of $U$ in \thmref{spatmodlks}
(ii) (a) follows from Theorem 4.1 in \cite{MWX}; and the local modulus of 
continuity of $U$ in \thmref{spatmodlks} (ii) (b) follows upon
applying \cite[Theorem 5.1]{MWX}.
\epfs

Finally we turn to the proof of the upper bounds on the uniform and local 
continuity moduli in the critical two dimensional case in
\thmref{spatmodlks} (iii).

\bpfs{Proof of \thmref{spatmodlks} (iii)}
Similarly to (\ref{spderlkseq}), we apply Lemma 3.1 to derive that for $d=2$
\beq\lbl{d2a}
\bsp
\E\big| U(t,x) - U(t,y) \big|^2
&=\int_{\R^2} \big| \e^{\i \langle x, \xi \rangle } - 
\e^{\i \langle y, \xi \rangle}\big|^2 S(\xi)d\xi\\
&=\frac{2}{\varepsilon \pi^{2}}\int_{\R^2} \lbk1- \cos \langle \xi, x-y \rangle \rbk 
\frac{ 1-\e^{-\frac{\varepsilon t}{4} (-2\theta + |\xi|^{2} )^{2}}}
{(-2\theta + |\xi|^{2} )^{2}} d\xi. 
\end{split}
\eeq
As in the proof of Theorem 1.2 (i), we assume $\vep=\vth=1$ and $|x-y| \le 1/2$. 
Let $\B_{1}:=\lbr\xi\in\R^2: \lab\xi\rab< 2\rbr$,  $\B_{2}:=
\big\{\xi\in\R^2: 2\le \lab\xi\rab< \frac 1 {|x-y|}\big\}$ and $\B_{3}
:=\big\{\xi\in\R^2: \lab\xi\rab \ge \frac 1 {|x-y|}\big\}$. By splitting 
the last integral in (\ref{d2a}) over three sets $\B_1$, $\B_2$, $\B_3$ and 
by using the inequalities in (\ref{elemineq}), one can  derive \beq\lbl{d2b}
\bsp
\E\big| U(t,x) - U(t,y) \big|^2 \le C |x-y|^2\, \int_{\B_2} \frac{d\xi} 
{|\xi|^2}  \le C |x-y|^2 \log \Big(\frac1 {|x -y|}\Big). 
\end{split}
\eeq
The desired upper bounds for the uniform and local continuity moduli for the
sample function $x \mapsto U(t,x) $ in $d=2$ follow from the Gaussian 
isoperimetric inequality and a Borel-Cantelli argument. Since this is the 
same as that in the proof of part (iii) of Theorem 3.1, we omit the details.
\epfs

It is natural to expect that (\ref{Eq:d2univ}) and (\ref{lsm2}) hold with 
``$\le$'' replaced by ``='', which would give the exact uniform and local 
continuity moduli for $x \mapsto U(t,x) $ in $d=2$. However, substantial 
extra work is needed for proving these statements. In particular, in order 
to apply the method in \cite{MWX}, one will have to establish the property 
of strong nondeterminism for $U(t, \cdot)$. Unfortunately the method in 
\cite{X07} does not seem  useful anymore and some new ideas may be needed. 

\subsection{The L-KS gradient temporal H\"older and modulus of continuity}
\lbl{sec:lksgradtemp}
We prove the temporal regularity of the spatial gradient $\pax U$ in 
\thmref{thm:gradtemplks}.
\bpfs{Proof of \thmref{thm:gradtemplks}}
Let $d=1$.  We start with the H\"older assertion for the gradient.  
Recall that $U(t,x)=X(t,x)-V(t,x)$ and that the temporal regularity 
of $U$ is totally determined by the rougher process $X$.  Similarly, 
the temporal regularity of the gradient $\pax U$ is entirely determined 
by the gradient of the rougher auxiliary process $X$ ($\pax X$)\footnote{\lbl{fn:gradV}
It can be shown that the smoothness assertions in \thmref{auxdecompspden} (ii) 
and (iii) (with $\gamma_{1}=1/8$) hold for $\pax V$.  Since the proof 
follows the same steps as the one for \thmref{auxdecompspden} (ii) and 
(iii) with straightforward modifications, we leave it to the interested 
reader.}.  Here, $\pax X$ plays the role of the auxiliary process for $\pax U$.  
We start with Parseval's identity to the integral in $y$ to get:
\begin{equation}\lbl{gradtempdiffL2}
\begin{split}
\mathbb{E} & \lbk \pax X(t,x) - \pax X(s,x)\rbk^{2}= \int_{\mathbb{R}} \int_{\mathbb{R}} \lab\pax\KKSepthtrposxy-\pax\KKSepthsrposxy \rab^{2} dr dy \\
& = \int_{\R}  \int_{\R} \xi^{2} \lab\FKKSepthtrposxxi-\FKKSepthsrposxxi \rab^{2} d\xi dr
\end{split}
\end{equation}
Since
\begin{equation}\lbl{gradtempom}
\FKKSepthtrposxxi = (2\pi)^{-d/2} \cdot e^{-\i \langle x,\xi\rangle - \frac{\varepsilon (t-r)}{8} 
(-2\theta + |\xi|^{2})^{2}} \ind_{ \{t > r\} },
\end{equation}
equation \eqref{gradtempdiffL2} becomes
\begin{equation}
\begin{split}
\E & \lbk \pax X(t,x) - \pax X(s,x)\rbk^{2}\\& =  \int_{\mathbb{R}}  
\int_{\mathbb{R}}  \df{\xi^{2}\lab e^{- \frac{\varepsilon (t -r)}{8} (-2 \theta 
+ | \xi |^{2} )^{2}} \ind_{ \{ t > r \}} - e^{- \frac{\vep (s -r)}{8} ( -2 \theta 
+| \xi |^{2} )^{2}} \ind_{ \{ s > r \}} \rab^{2}}{(2\pi)^{d}} dr d\xi.
\end{split}
\end{equation}
Now, we apply Parseval's identity to the inner integral in $r$. To this end, let
\begin{equation*}
\phi (r,\xi) = e^{-\frac{\varepsilon (t-r)}{8} (-2\theta + | \xi |^{2})^{2}} \ind_{\{ t > r \} }
 - e^{-\frac{\varepsilon(s - r)}{8} (-2\theta + |\xi |^{2} )^{2}} \ind_{ \{ s > r \}}.
 \end{equation*}
 Its Fourier transform in $r$ is
 \begin{equation*}
 \widehat{\phi} (\tau,\xi) = ( e^{\i\tau t} - e^{\i\tau s} ) \frac1{\i\tau +
  \frac{\varepsilon}{8} (-2\theta + | \xi |^{2})^{2}}.
 \end{equation*}
Hence, by Parseval's identity, inequality \eqref{elemineq} (a), its related inequality
\beq\lbl{eq:shcosineq}
1-\cos\lpa z\cdot\tau\rpa\le2\lpa1\wedge\lab z\rab^{2\alpha}\rpa\lbk1-\cos(\tau)\rbk;\ 0<\alpha\le1,
\eeq
and the asymptotic
\beq\lbl{eq:lksgradspdenas}
\bsp
\Delta(\tau):=(2\pi)^{-1}\int_{\R}\frac{\xi^{2}}{\tau^{2}+\frac{\vep^{2}}{64}(-2\theta + |\xi|^{2} )^{4}}d\xi\sim\frac{C}{\lab\tau\rab^{5/4}}, \mbox{ as }\lab\tau\rab\nearrow\infty,
\end{split}
\eeq
we get, for a large enough $N$, that
\begin{equation}\lbl{eq:tempspdengradlks}
\begin{split}
\E & \lbk\pax X(t,x) -\pax X(s,x)\rbk^{2} = (2\pi)^{-1} \int_{\mathbb{R}} \int_{\mathbb{R}} \xi^{2}\lab \widehat{\phi} (\tau,\xi) \rab^{2} d\tau d\xi\\
& = \frac1\pi \int_{\mathbb{R}} (1 - \cos ((t-s)\tau)) \int_{\mathbb{R}} \frac{\xi^{2} d\xi}{\tau^{2} + \frac{\varepsilon^{2}}{64} (-2\theta + | \xi|^{2} )^{4}} d\tau
\\&\le C|t-s|^{2\alpha}\lbk\int_{0}^{N}{(1 - \cos (\tau))}\Delta(\tau)d\tau+\int_{N}^{\infty}\lab\tau\rab^{2\alpha-\tf54}d\tau\rbk
\\&\le C|t-s|^{2\alpha},\ 0<\alpha<1/8.
\end{split}
\end{equation}
It follows that $\pax X(\cdot,x)$ is $\gamma$-H\"older continuous in time, with 
$\gamma\in(0,1/8)$.  This, together with the gradient decomposition
\beq\lbl{eq:graddec}
\pax U(t,x)=\pax X(t,x)-\pax V(t,x),
\eeq
and the fact that $\pax V$ is temporally smooth (see footnote \ref{fn:gradV}) establish 
the H\"older regularity assertion for $\pax U$ in \thmref{thm:gradtemplks}.

Turning now to the uniform and local spatial continuity moduli results for the L-KS 
gradient, $\pa_{x}U$, in \thmref{thm:gradtemplks}.  Equation \eqref{eq:tempspdengradlks} 
means that $\pax X$ has stationary increments and the spatial spectral density 
of $\pa_{x}X$ and its asymptotic behavior are given by \eqref{eq:lksgradspdenas}.

Equation \eqref{eq:lksgradspdenas} and Theorem 2.1 in \cite{X07} imply that, for 
every fixed $x\in\R$, the gradient $\{\pa_{x}X(t,x), t\ge0\}$ is temporally strongly 
locally nondeterministic. Namely, for any $T>0$, there is a positive constant $c$ 
such that for all $t\in(0,T]$ and all $0 < r \le 1\wedge|t|$
\beq\lbl{SLNDtlkseq}
\var \lpa\pax X(t,x) | \pax X(s,x);s\in[0,T],|t-s|\ge r\rpa \ge c r^{\frac{1}{4}}
\eeq
Also, \eqref{lksgradspdensp} and Theorem 2.5 in \cite{X07} imply the double sided 
second moment bounds
\beq\lbl{dstemplks}
\E [\pax X(t,x)-\pax X(s,x)]^{2} \asymp|t-s|^{\frac{1}{4}};\ \forall s,t\in[0,T].
\eeq
 Thus, the uniform modulus of continuity of $\pax X$
 \beq\lbl{eq:tempumlksauxgrad}
\P\lbk\lim_{\delta\searrow0} \sup_{\substack{|t-s|<\delta\\t,s\in\tint} }\frac{\lab\pa_{x} X(t,x)-\pa_{x}X(s,x)\rab}{|t-s|^{{1}/{8}}\sqrt{\log\lbk1/|t-s|\rbk}}=k\rbk=1,
\eeq
for every compact interval $\tint\subset\Rp$ and for some constant $k>0$, follows 
from Theorem 4.1 in \cite{MWX}.  The local modulus of continuity of $\pax X$
\beq\lbl{eq:templmlksauxgrad}
\P\lbk\lim_{\delta\searrow0} \sup_{\substack{|t-s|<\delta} }\frac{\lab\pa_{x} X(t,x)-\pa_{x}X(s,x)\rab}{\delta^{1/{8}}\sqrt{\log\log\lbk1/\delta\rbk}}=k\rbk=1,
\eeq
follows upon applying \cite[Theorem 5.1]{MWX}.  The corresponding continuity moduli 
assertions for the gradient $\pax U$ in \thmref{thm:gradtemplks} follow from those 
of the auxiliary process $\pax X$ (\eqref{eq:tempumlksauxgrad} and 
\eqref{eq:templmlksauxgrad}), the decomposition \eqref{eq:graddec}, and 
the smoothness of $\pax V$ (see footnote \ref{fn:gradV}).
\epfs

\section{The time-fractional SPIDEs: Proofs of Theorems 1.4--1.7}\lbl{sec:tfspidepfs}
As with the L-KS SPDE case, Allouba obtained in \cite{Abtbmsie,Atfhosie},  the time and space H\"older exponents
$\gamma_{t}\in(0,(2\beta^{-1}-d)/4\beta^{-1})$ and $\gamma_{s}\in(0,((4-d)/2)\wedge1)$, respectively, after establishing
the sharp dimension-and-$\beta$-dependent upper bounds
\begin{equation}\lbl{ubtfsol}
\bc
\E \lbk U_{\beta}(t,x) - U_{\beta}(s,x)\rbk^{2q}&\le C_{d,\beta}\lab t-s\rab^{\tf{(2\beta^{-1}-d)q}{2{\beta^{-1}}}},\\
\E \lbk U_{\beta} (t,x) - U_{\beta}(t,y)\rbk^{2q}&\le {C}_{d}\, |x-y|^{2q\alpha_{d}};\quad \alpha_{d}\in J_{d},
\ec
\end{equation}
for the more general nonlinear time-fractional SPIDE \eqref{spidea}, with Lipschitz condition on $a$, for all $x \in \Rd$, $t,s \in[0,T]$,
$q\ge1$, $1\le d\le 3$, $\beta\in\{1/2^{k};k\in\N\}$, and for the intervals $J_{d}$ as in \eqref{intervals}.  Now, we take the
spectral/harmonic analysis and solution decomposition route we took in \secref{tmpmodlkssec}---with the time-fractional kernel
$\Kbetatx$ replacing the L-KS one---to get the exact dimension-dependent temporal and spatial uniform and local moduli of
continuity in \thmref{tempmodtf} and in \thmref{spatmodtf}.

Assume without loss of generality that $u_{0} = 0$, then the $\beta$ time-fractional SPIDE solution is given by
\begin{equation}\lbl{zin}
U_{\beta}(t,x) = \int_{\mathbb R^{d}} \int_{0}^{t} \Kbetatsxy\sW(ds\times dy), \quad t \ge 0, \, x \in \Rd.
\end{equation}

\subsection{Temporal modulus}
Throughout this subsection, let $x\in\mathbb R^{d}$ be fixed but arbitrary.
Let $U_{\beta}$ be the solution to the time-fractional SPIDE \eqref{spide}, given by (\ref{zin}).
Following the template used in the L-KS proofs, we first introduce the following auxiliary Gaussian process
$\{ X_{\beta}(t,x), t \in \mathbb{R}_{+}\}$:
\begin{equation}\lbl{auxgauss2}
X_{\beta}(t,x) = \int_{\mathbb{R}^{d}} \int_{\mathbb{R}} \lpa\Kbetatrposxy-\Kbetamrposxy \rpa \sW(dr \times dy),
\end{equation}
where $x \in \mathbb R^{d}$ is arbitrary but fixed.  Then the solution $U_{\beta}$ may be decomposed as
$U_{\beta}(t,x) = X_{\beta}(t,x) - V_{\beta}(t,x) $, where
\begin{equation}\label{Vbeta}
V_{\beta}(t,x) = \int_{\mathbb{R}^{d}} \int_{\mathbb{R}_{-}} \lpa\Kbetatrposxy-\Kbetamrposxy \rpa \sW(dr \times dy).
\end{equation}

We start by proving the following crucial result for the auxiliary process $X_{\beta}$ and the
smoothness of $V_{\beta}$. 
\bthm\lbl{auxdecompspiden}
Assume $d\in\{1,2,3\}$ and $0<\beta\le1/2$.  Let $X_{\beta}$ be as defined in \eqref{auxgauss2}.
\begin{enumerate}
\renewcommand{\labelenumi}{(\roman{enumi})}
\item The Gaussian process $\lbr X_{\beta}(t,x); t\ge 0\rbr$ has stationary temporal increments.
Moreover we have
\begin{equation*}
\mathbb{E} \lbk X_{\beta}(t,x) - X_{\beta}(s,x)\rbk^{2} = 2 \int_{\mathbb{R}} [1-\cos((t-s) \tau)] \Delta_{\beta}(\tau) d\tau,
\end{equation*}
where the spectral density $\Delta_{\beta}$ is given by
\begin{equation}\lbl{tspdentf}
\Delta_{\beta}(\tau) = (2\pi)^{-d} \frac{1}{|\tau|^{2-(\beta d)/2}}\int_{\mathbb{R}^{d}} \frac{d\xi}{1+\lab\xi\rab^{2}\cos\big(\frac{\pi\beta}
{2}\big)+\tf14\lab\xi\rab^{4}}.
\end{equation}
\item For each $k \ge 1$, there exists a modification of $\{V_{\beta}(t,x), t\in \mathbb{R}_{+}\}$
such that its (temporal) sample function
is almost surely continuously $k$-times differentiable on $(0,\infty)$.
\item Let $ H =\frac{2 - \beta d}{4 }$. There is a finite constant $C$ such that
\begin{equation}\label{Eq:modVb}
\lim_{\varepsilon \to 0} \sup _{s, t\in [0,\, \varepsilon]} \frac{\left| V_\beta(t,x) - V_\beta(s,x) \right|}
{|t-s|^H \sqrt{\log \log (1/|t-s|)} } \leq C  \ \  \mbox{ a.s. }
\end{equation}
\end{enumerate}
\ethm
\bpf
To verify (i), we apply Parseval's identity to the integral in $y$ to get:
\begin{equation}\lbl{tempduffauxL2}
\begin{split}
\mathbb{E} & \lbk X_{\beta}(t,x) - X_{\beta}(s,x)\rbk^{2}= \int_{\mathbb{R}^{d}} \int_{\mathbb{R}} \lab\Kbetatrposxy-\Kbetasrposxy \rab^{2} dr dy \\
& = \int_{\R} \int_{\Rd}  \lab\FKbetatrposxxi-\FKbetasrposxxi \rab^{2} d\xi dr.
\end{split}
\end{equation}
Since
\begin{equation}\lbl{tempom2}
\FKbetatrposxxi = (2\pi)^{-d/2} \cdot \e^{-\i \langle x,\xi\rangle} E_{\beta}\bigg(-\frac{\lab\xi\rab^{2}}{2}(t-r)^{\beta}\bigg) \ind_{ \{t > r\} },
\end{equation}
equation \eqref{tempduffauxL2} becomes
\begin{equation}
\begin{split}
\E \lbk X_{\beta}(t,x) - X_{\beta}(s,x)\rbk^{2}=(2\pi)^{-d}  \int_{\mathbb{R}^{d}}
\int_{\mathbb{R}}  \lab\phi (r,\xi)\rab^{2} dr d\xi,
\end{split}
\end{equation}
where \begin{equation*}
\phi (r,\xi) = E_{\beta}\bigg(-\frac{\lab\xi\rab^{2}}{2}(t-r)^{\beta}\bigg)
\ind_{ \{ t > r \}} - E_{\beta}\bigg(-\frac{\lab\xi\rab^{2}}{2}(s-r)^{\beta}\bigg)
\ind_{ \{ s > r \}}.
 \end{equation*}
Now, we apply Parseval's identity to the inner integral in $r$. To this end,
assume for simplicity and without loss of generality
that $\beta\in\{1/2^{k};k\in\N\}$.  In this case, using \lemref{lm:MLFT} above, the Fourier transform of $\phi$ in $r$
is
 \begin{equation}\lbl{FTrtf}
 \widehat{\phi} (\tau,\xi) = (\e^{\i\tau t} - \e^{\i\tau s})\frac{\i^{\beta-1}\tau^{\beta-1}}{\i^{\beta}\tau^{\beta}+\tf12\lab\xi\rab^{2}}.
  \end{equation}
Hence, by Parseval's identity, we get
\begin{equation}
\begin{split}
\E& \lbk X_{\beta}(t,x) - X_{\beta}(s,x)\rbk^{2}=(2\pi)^{-d} \int_{\mathbb{R}^{d}} \int_{\mathbb{R}} \big| \widehat{\phi} (\tau,\xi) \big|^{2} d\tau d\xi\\
& =  2(2\pi)^{-d}  \int_{\mathbb{R}} (1 - \cos ((t-s)\tau)) \int_{\mathbb{R}^{d}} \frac{\tau^{2(\beta-1)}}
{\tau^{2\beta}+\lab\xi\rab^{2}\tau^{\beta}\Re\lpa\i^{\beta}\rpa+\tf14\lab\xi\rab^{4}} d\xi d\tau
\\&=2(2\pi)^{-d}  \int_{\mathbb{R}} (1 - \cos ((t-s)\tau))\frac{d\tau}{|\tau|^{2-(\beta d)/2}} \int_{\mathbb{R}^{d}} \frac{d\xi}{1+\lab\xi\rab^{2}
\cos\big(\frac{\pi\beta}{2}\big)+\tf14\lab\xi\rab^{4}}.
\end{split}
\end{equation}
The proof of (i) is complete.  The proof of parts (ii) and (iii) is very similar to the proof of \thmref{auxdecompspden} (ii) and (iii), with
now obvious modifications.  We leave the details to the interested reader.
\epf

Using the asymptotic behavior of the spectral density $\Delta_{\beta}$ in \eqref{tspdentf} \thmref{auxdecompspiden} (i), 
we proceed as in \secref{tmpmodlkssec}
to obtain the following SLND and two-sided bounds for $X_{\beta}$.
\bcr[Temporal SLND and double-sided bounds for $X_{\beta}$]\lbl{SLNDtspides}
Let $d\in\{1,2,3\}$ and let $0<\beta\le1/2$.  For any $T>0$, there is a positive constant $c$ such that for all $t\in(0,T]$ and all $0 < r \le
1\wedge|t|$ such that
\beq\lbl{SLNDttfeq}
\var \lpa X_{\beta}(t,x) | X_{\beta}(s,x);s\in[0,T],|t-s|\ge r\rpa \ge c\, r^{\frac{2- \beta d}{2 }},
\eeq
and
\beq\lbl{dstemptf}
\E [X_{\beta}(t,x)-X_{\beta}(s,x)]^{2} \asymp|t-s|^{\frac{2 - \beta d}{2}};\quad \forall s,t\in[0,T].
\eeq
Moreover, the function $\sigma_\beta^2(h) = \E [X_{\beta}(t+h,x)-X_{\beta}(t,x)]^{2} $ is regularly varying at $h=0$ of order
$ (2- \beta d)/{2}$.
\ecr
\bpf
The property of the spectral density $\Delta_{\beta}$ in \thmref{auxdecompspiden} (i)
and Theorem 2.1 in \cite{X07} imply \eqref{SLNDttfeq}. Similarly,
\thmref{auxdecompspiden} (i) and Theorem 2.5 in \cite{X07} imply \eqref{dstemptf}. Finally, since the spectral density
$\Delta_\beta$ is regularly varying of order $-(2 - \frac{\beta d} 2)$ at $\infty$,  the last conclusion of the corollary follows
from Theorem 1 in \cite{Pitman68}.
\epf

From \coref{SLNDtspides}, conditions (C1) and (C2) above ((A1) and (A2) in \cite{MWX}) hold.
Now, applying the results in \cite{MWX} on the uniform and local continuity moduli for the
auxiliary Gaussian processes to $\{X_{\beta}(t,x),t \ge 0\}$, we get the following theorem
on the time regularity of $X_{\beta}$. Recall that $H = \frac{2- \beta d}{4 }$.
 \begin{thm}\lbl{auxmod2}
Let $x\in \mathbb{R}^{d}$, $d=1,2,3$, be fixed; let $0<\beta\le1/2$; and let
$\{ X_{\beta}(t,x), t \ge 0 \}$ be defined as in \eqref{auxgauss2} above.  Then,
 \ben\rencomrom
\item (Uniform Modulus of Continuity)  for every compact interval $\tint \subseteq
\mathbb{R}_{+}$
\begin{equation}\lbl{umctauxtf}
\lim_{\delta\searrow 0} \sup_{\substack{|s-t|< \delta \\ s, t \in I} } \frac{|X_{\beta}(t,x) - X_{\beta}(s,x)|}{|s-t|^{H}\sqrt{\log \frac1{|t-s|}}}
= k_{6}^{(\beta,d)}; \mbox{ a.s.}\end{equation}
\item (Local Modulus of Continuity) and for every fixed $t \ge 0$
\begin{equation}\lbl{lmctauxtf}
\lim_{\delta\searrow 0 } \frac{\sup_{|s-t| <\delta} | X_{\beta}(s,x) - X_{\beta}(t,x)|}{\delta^{H} \sqrt{\log\log \frac1{\delta}}}
= k_{7}^{(\beta,d)}; \mbox{ a.s.},
\end{equation}
\een
where $k_{i}^{(\beta,d)}$ ($i = 6,7$) are positive and finite constant that depend on $d$ and $\beta$, but are  independent of $x$.
 \end{thm}
\bpf The uniform modulus of continuity of $X_{\beta}$ in \eqref{umctauxtf} follows from Theorem 4.1 in \cite{MWX}.  The local modulus
of continuity of $X_{\beta}$ in \eqref{lmctauxtf} follows upon applying \cite[Theorem 5.1]{MWX}. The constant $k_{i}^{(\beta,d)}$ ($i = 6, 7$)
do not depend on $x\in\Rd$ since the distribution of the process $\{X_{\beta}(t,x), t \ge 0\}$ does not depend on $x$, see \eqref{sdastimelks}.
\epf

Now we use the decomposition $U_{\beta}(t,x) = X_{\beta}(t,x) - V_{\beta}(t,x)$ to prove  \thmref{tempmodtf} (i).
\bpfs{Proof of \thmref{tempmodtf} (i)} As in the proof of of Theorem 1.1 (i), we see that (\ref{Eq:SPI_umod}) and (\ref{Eq:SPI_LIL})
follow from the aforementioned decomposition and Theorems  4.1 and 4.2.
\epfs

\subsection{Time-fractional SPIDEs are \emph{not} bifractional Brownian motions}
%
Let $U_{\beta}$ be the solution to the time-fractional SPIDE \eqref{spide}, given in (\ref{zin}).
We now characterize the law of $\lbr U_{\beta}(t,x);t\ge0\rbr$---which we call the $\beta$
time-fractional SPIDE law---and we show that, unlike the L-KS SPDE, it's fundamentally different
from the  bifractional Brownian motion law.

\bpfs{Proof of \thmref{tempmodtf} (ii)} For any $0<s<t$, we use Parseval's identity to get
\begin{equation*}\lbl{tfcovpars}
\bsp
&\E  [U_{\beta}(t,x) U_{\beta}(s, x)] = \int_{\mathbb{R}^{d}} \int_{0}^{s} \Kbetatrxy\Kbetasrxy dr dy
\\& =\int_{0}^{s}  \int_{\mathbb{R}^{d}} \FKbetatrxxi  \overline{\FKbetasrxxi}d\xi dr\\
& = (2\pi)^{-d} \int_{0}^{s}  \int_{\mathbb{R}^{d}} E_{\beta}\Big(-\frac{\lab\xi\rab^{2}}{2}(t-r)^{\beta}\Big) 
E_{\beta}\Big(-\frac{\lab\xi\rab^{2}} {2}(s-r)^{\beta}\Big) d\xi dr\\
& = (2\pi)^{-d} \int_{\mathbb{R}^{d}}\int_{0}^{s}  \sum_{k=0}^{\infty}\Bigg[\sum_{j=0}^{k}\frac{(t-r)^{\beta j}(s-r)^{\beta(k-j)}}
{\Gamma(1+\beta j)\Gamma(1+\beta(k-j))}\Bigg]\frac{(-1)^{k}\lab\xi\rab^{2k}}{2^{k}}drd\xi
\\&=(2\pi)^{-d} \int_{\mathbb{R}^{d}}  \sum_{k=0}^{\infty}\Bigg[\sum_{j=0}^{k}\frac{t^{\beta j}s^{\beta(k-j)+1}\;   _{2}
F_{1}\lpa1,-\beta j;2+\beta(k-j);\frac{s}t\rpa}{\lbk\beta(k-j)+1\rbk\Gamma(1+\beta j)\Gamma(1+\beta(k-j))}
\Bigg]\frac{(-1)^{k}\lab\xi\rab^{2k}}{2^{k}}d\xi,
\end{split}
\end{equation*}
which proves the covariance assertion of \thmref{tempmodtf} (ii). Moreover, we see from the above that 
for any constant $c > 0$, 
$$\E  [U_{\beta}(ct,x) U_{\beta}(cs, x)] = c^{\frac{2 - \beta d} 2}\E  [U_{\beta}(t,x) U_{\beta}(s, x)].$$
Hence the Gaussian process $U_\beta = \{U_{\beta}(t,x), t \ge 0\}$ is self-similar with index $(2 - \beta d)/4$.

To show that $U_\beta$  does not have the same law as any bifractional Brownian motion, and to give 
an alternative form of the covariance  function $\E  [U_{\beta}(t,x) U_{\beta}(s, x)]$, we exploit the 
form of the kernels $\Kbetatx$ directly rather than using their Fourier transforms.  Computing
the covariance of $U_{\beta}$, directly we obtain
\beq\lbl{tfcov}
\begin{split}
&\E  [U_{\beta}(t,x) U_{\beta}(s, x)] = \int_{\mathbb{R}^{d}} \int_{0}^{s} \Kbetatrxy\Kbetasrxy dr dy
\\&=2^{k}\int_{0}^{s}\lbr\int_0^\infty\int_0^\infty\lbk\int_{\Rd} \psoxy \puoxy dy\rbk\right.
\\&\times \lpa\int_{\Rpop^{k-1}} \ptrzsk\prod_{i=0}^{k-2}\ptzskiimo ds_{2}\cdots ds_{k}\rpa\\&\left.\times
\lpa\int_{\Rpop^{k-1}} \psrzuk\prod_{i=0}^{k-2}\ptzukiimo du_{2}\cdots du_{k}\rpa ds_{1}du_{1}\rbr dr
\\&= \int_{0}^{s}\lbr\int_0^\infty\int_0^\infty\lbk\df{2^{k}}{\lbk2\pi (s_{1}+u_{1})\rbk^{d/2}}\rbk\right.
\\&\times \lpa\int_{\Rpop^{k-1}} \ptrzsk\prod_{i=0}^{k-2}\ptzskiimo ds_{2}\cdots ds_{k}\rpa\\&\left.\times
\lpa\int_{\Rpop^{k-1}} \psrzuk\prod_{i=0}^{k-2}\ptzukiimo du_{2}\cdots du_{k}\rpa ds_{1}du_{1}\rbr dr.
\end{split}
\eeq
Gathering the two inside integrals and transforming to polar coordinates
$(s_{i},u_{i})\mapsto(\rho_{i},\theta_{i})$, $i=1,\ldots,k$,
letting  $\underline{\rho}=(\rho_{1},\ldots,\rho_{k})$ and
$\underline{\theta} =(\theta_{1},\ldots,\theta_{k})$, letting $\I_{\pi}=(0.\pi/2)$,
and noticing that all $\rho_{i}$ for $i=2,3,\ldots,k$ cancel when $k\ge2$; equation
\eqref{tfcov} gives us the covariance $\E  [U_{\beta}(t,x) U_{\beta}(s, x)]$ as
\beq\lbl{sumofsqofisltbmden}
\bsp
C_{\beta,d}\int_{0}^{s}\int\limits_{\I_{\pi}^{k}}\int\limits_{\Rpop^{k}}\df{\df{\e^{-\rho_{k}^2
\big[\frac{\cos^{2}(\theta_{k})}{4(t-r)}+\frac{\sin^{2}(\theta_{k})}{4(s-r)}\big]}}
{\sqrt{(t-r)(s-r)}}\ds\prod_{i=0}^{k-2} \e^{-\frac{\rho_{k-i-1}^2}{4\rho_{k-i}}
\big[\frac{\cos^{2}(\theta_{k-i-1})}{\cos(\theta_{k-i})}+\frac{\sin^{2}(\theta_{k-i-1})}
{\sin(\theta_{k-i})}\big]}}
{\rho_{1}^{\frac{d}{2}-1}\big[sin(\theta_{1})+\cos(\theta_{1})\big]^{\frac{d}{2}}\ds
\prod_{i=0}^{k-2}\sqrt{\sin(\theta_{k-i})\cos(\theta_{k-i})}}d\underline{\rho} d\underline{\theta}dr.
\end{split}
\eeq
To simplify our computations, it is enough for our purposes to assume that $k=1$ or $\beta=1/2$
(the Brownian-time Brownian motion case) and take $d=2$. The integrals with respect to $\rho$
and then $r$ in equation \eqref{sumofsqofisltbmden} then give
\beq\lbl{covbtbmsierhor}
\bsp
&\int_{0}^{s}\int\limits_{\Rpop}\df{\e^{-\rho^2\big[\frac{(t-r)-(t-s)\cos^{2}(\theta)}{4(t-r)(s-r)}\big]}}
{\sqrt{(t-r)(s-r)}}d\rho dr=\int_{0}^{s}\frac{\sqrt\pi}{\sqrt{t-r-(t-s)\cos^{2}(\theta)}}dr
\\&=2\sqrt\pi\lbk\sqrt{t\sin^{2}(\theta)+s\cos^{2}(\theta)}-\sqrt{t-s}\lab\sin(\theta)\rab\rbk.
\end{split}
\eeq
Finally, in the BTBM $\beta=1/2$ and $d=2$ case, the covariance \eqref{sumofsqofisltbmden}  becomes
\beq\lbl{btbmsiecovd2}
\bsp
&2\sqrt\pi C_{\beta,d}\int_{0}^{\pi/2}\frac{\sqrt{t\sin^{2}(\theta)+s\cos^{2}(\theta)}
-\sqrt{t-s}\sin(\theta)} {\sin(\theta)+\cos(\theta)}d\theta
\\&=2\sqrt\pi C_{\beta,d}\lbk\int_{0}^{\pi/2}\frac{\sqrt{(t-s)\sin^{2}(\theta)+s}}{\sin(\theta)
+\cos(\theta)}d\theta-\frac\pi4\sqrt{t-s}\rbk
\\&=2\sqrt\pi C_{\beta,d}\lbk\frac {1}{8} \left(-2\tanh^{-1}\lpa{\frac {  \sqrt {2}s-\sqrt {
2}t-2s  }{\sqrt {2t(s+t)}}}\rpa\sqrt {t+s}\right.\right.
\\&\quad\left.\left.-2\sqrt {t+s}\lbr\tanh^{-1}\lpa{\frac {  \sqrt {2}s-\sqrt {2}t+2s }{\sqrt {2t(s+t)}}}\rpa-\Re\lpa\tanh^{-1}\lpa\frac {2s+t}
{2\sqrt {s(s+t)}}\rpa\rpa\rbr\right.\right.
\\&\quad\left.\left.+\sqrt{t-s}\lbr-2\sin^{-1}\left({\frac {2s-t}{t}}\right)-4\ln
\left( \sqrt {t}+\sqrt {t-s} \right)+2\ln\left(s \right)+\pi\rbr  \right.\right.\Bigg)
\\&\quad-\frac{\pi}{4}\sqrt{t-s}\Bigg].
\end{split}
\eeq
It can now be easily verified that the bracketed term is not equal to
\beq\lbl{corbifbm}
C\lbk\sqrt{t+s}-\sqrt{t-s}\,\rbk
\eeq
for any constant $C$. Thus the law of the BTBM SPIDE is not a bifractional Brownian motion in $d=2$.
The cases $d=1,3$ and $\beta<1/2$ are similar and we omit them.
\epfs

\subsection{Spatial modulus}
Without loss of generality, we again assume that $\un=0$, and the random field
solution $U_{\beta}$ is given by \eqref{zin}. Fix an arbitrary $t>0$ throughout this subsection.
Our spatial results for this case crucially depend on the following Lemma.
\begin{lem}[Time-fractional SPIDEs spatial spectral density]\lbl{sspdenspide}
Let $d=1,2,3$ and $0<\beta\le1/2$. The centered Gaussian random field
$\{ U_{\beta}(t,x), x \in \mathbb{R}^{d}\}$ is stationary with spectral density
\begin{equation*}
\bsp
S_{\beta}(\xi)&=(2\pi)^{-d}\int_{0}^{t} E^{2}_{\beta}\bigg(-\frac{\lab\xi\rab^{2}}{2}(t-r)^{\beta}\bigg) dr
\\&= (2\pi)^{-d}\sum_{k=0}^{\infty}\df{(-1)^{k}a_{k}\lab\xi\rab^{2k}t^{\beta k+1}}{2^{k}(\beta k+1)},
\quad \forall \xi \in \mathbb{R}^{d},
\end{split}
\end{equation*}
where
$$a_{k}=\sum_{j=0}^{k}\df{1}{\Gamma(1+\beta j)\Gamma(1+\beta(k- j))}.$$
\end{lem}
\bpf
Computing the covariance of $U_{\beta}$, we use (4.2) and Parseval's identity to get
\begin{equation}
\begin{split}
\mathbb{E} & \lbk U_{\beta}(t,x) U_{\beta}(t,y)\rbk = \int_{0}^{t}  \int_{\mathbb{R}^{d}}
\Kbetatrxz\Kbetatryz dz dr\\
& = \int_{0}^{t}  \int_{\mathbb{R}^{d}}  \FKbetatrxxi\overline{\FKbetatryxi}d\xi dr\\
& = (2\pi)^{-d} \int_{0}^{t}  \int_{\mathbb{R}^{d}} \e^{\i \langle\xi, x-y\rangle} E^{2}_{\beta}\bigg(-\frac{\lab\xi\rab^{2}}{2}(t-r)^{\beta}\bigg) d\xi dr\\
& = (2\pi)^{-d} \int_{\mathbb{R}^{d}} \e^{\i \lang \xi, x-y\rang} \int_{0}^{t} E^{2}_{\beta}\bigg(-\frac{\lab\xi\rab^{2}}{2}(t-r)^{\beta}\bigg) dr d\xi\\
& = (2\pi)^{-d} \int_{\mathbb{R}^{d}} \e^{\i \lang \xi, x-y\rang} \sum_{k=0}^{\infty}\df{(-1)^{k}a_{k}\lab\xi\rab^{2k}t^{\beta k+1}}{2^{k}(\beta k+1)} d\xi.\\
\end{split}
\end{equation}
Thus, the conclusions of the Lemma follows.
\epf
\subsubsection{Spectral asymptotics for $0<\beta<1/2$}
We need the asymptotic behavior of $S_{\beta}$ at $\infty$, which is captured in
the next lemma for the case $0<\beta<1/2$.\footnote{Another approach is used---and
a different result is obtained---for the case $\beta=1/2$, which we provide next.}
\blm\lbl{tfsdasnothalf}  Fix an arbitrary $t>0$ and $d=1,2,3$. If and $0<\beta<1/2$, then
the spectral density $S_\beta$ satisfies $0 < S_\beta(0) < \infty$ and has the asymptotic behavior
\begin{equation}\lbl{spidespcas}
S_{\beta}(\xi) \sim \frac{C_{t,\beta,d}}{|\xi|^{d+2\gamma}};\ \ \mbox{ as }|\xi|\to\infty,
\end{equation}
for some finite constant $C_{t, \beta,d}$ where $\gamma = 2-\frac{d}{2}$. Moreover,
$S_{\beta}(\xi) \le C_{t,\beta,d}|\xi|^{-(d+2\gamma)}$ for all $\xi \in \Rd\backslash\{0\}$.
\elm
\bpf  Let $d\in\{1,2,3\}$, $0<\beta<1/2$, and $\gamma = 2-\frac{d}{2}$. Clearly,  $0 < S_\beta(0) < \infty$
follows from Lemma 4.1.  Moreover,
by using the asymptotic property of the Mittag-Leffler function in \eqref{mlas}, we get that
as $|\xi| \to \infty$,
\beq\lbl{astftime1}
\bsp
S_{\beta}(\xi)&=\df{\ds\int_{0}^{t}E^{2}_{\beta}\Big(-\frac{\lab\xi\rab^{2}}{2}(t-r)^{\beta}\Big) dr}{(2\pi)^{d}}
\sim\df{\ds\int_{0}^{t}\Big(\frac{\lab\xi\rab^{2}}{2}(t-r)^{\beta}\Big)^{-2}
\Gamma^{2}(1-\beta)dr}{(2\pi)^{d}}\\
&=\frac{4\Gamma^{2}(1-\beta)}{(2\pi)^{d}}\,{\frac {{t}^{1-2\beta}}{\left( 1-2\beta \right) {\lab\xi\rab}^{d+2\gamma}}},
\end{split}
\eeq
and \eqref{spidespcas} follows with
$$C_{\beta,d, t}=\,{\frac {4\Gamma^{2}(1-\beta){t}^{1-2\,\beta }}{ (2\pi)^{d}\left( 1-2\,\beta \right)}}.$$
Finally, the upper bound for $S_{\beta}(\xi)$  follows from the first equation in \eqref{astftime1}
and the upper bound for $M_\beta(-x)$ in \eqref{elemineqML}. The proof is complete.
\epf
\subsubsection{Spectral asymptotics for the critical fraction $\beta=1/2$}
Since the second integral in \eqref{astftime1} diverges at $\beta=1/2$, the proof of
the case $0<\beta<1/2$ above does not work for the case $\beta=1/2$.  The reason,
as is clear from \thmref{spatmodtf} and \thmref{spatmodtfhalf}, is that
the case $\beta=1/2$ has a rougher modulus than that of  $\beta<1/2$.  This is captured
in the following lemma.
\blm[Spectral asymptotic behavior at $\beta=1/2$]\label{Lem:sd2}
\lbl{tfsdashalf} Fix an arbitrary $t>0$,
and let $\beta=1/2$.  As $|\xi| \to \infty$, the spectral density has the asymptotic behavior
\begin{equation}\lbl{spidespcashalf}
S_{1/2}(\xi) \sim \frac{C_{t,d}}{|\xi|^{d+2\gamma}}\log\lab\xi\rab,
\end{equation}
for some finite constant $C_{t,d}$ where $\gamma = 2-\frac{d}{2}$, $d=1,2,3$. Moreover,
$S_{1/2}(\xi) \le C_{t,d}|\xi|^{-(d+2\gamma)}\log\lab\xi\rab$ for all $\xi \in \Rd\backslash\{0\}$.
\elm
\bpf By \lemref{sspdenspide}, \lemref{FT}, and footnote \ref{footFT}, we have
\beq\lbl{astftimehalf}
\bsp
S_{1/2}(\xi)&=(2\pi)^{-d}\int_{0}^{t}\lpa\e^{\frac {r}4\lab\xi\rab^{4}}\lbk\frac{2}{\sqrt{\pi}}\int_{\frac{\sqrt{r}\lab\xi\rab^{2}}{2}}^{\infty}
\e^{-\tau^{2}}d\tau\rbk\rpa^{2}dr
\\&=\frac{(2\pi)^{-d}}{\lab\xi\rab^{4}}\int_{0}^{t\lab\xi\rab^{4}}\e^{\frac {\rho}2}\lpa\lbk\frac{2}{\sqrt{\pi}}\int_{\frac{\sqrt{\rho}}{2}}^{\infty}
\e^{-\tau^{2}}d\tau\rbk\rpa^{2}d\rho
\\&\sim\frac{(2\pi)^{-d}}{\lab\xi\rab^{4}}\int_{1}^{t\lab\xi\rab^{4}}\e^{\frac {\rho}2}\lpa\lbk\frac{2}{\sqrt{\pi}}\int_{\frac{\sqrt{\rho}}{2}}^{\infty}
\e^{-\tau^{2}}d\tau\rbk\rpa^{2}d\rho, \ \ \mbox{ as }|\xi|\to\infty,
\end{split}
\eeq
where we have used the change of variable $\rho=r|\xi|^{4}$. Now, using the standard
asymptotic for Mills' ratio for the standard normal random variable,
$m(x)=\frac{\int_{x}^{\infty}\e^{-u^{2}/2 }du}{\e^{-x^{2}/2}}\sim1/x$, we get
\beq\lbl{eq:astftimehalf}
\bsp
S_{1/2}(\xi)&\sim\frac{C_{d}}{\lab\xi\rab^{4}}\int_{1}^{t\lab\xi\rab^{4}}\e^{\frac {\rho}
2}\bigg(\frac{2\e^{\frac{-\rho}{4}}}{\sqrt{\rho}}\bigg)^{2} d\rho
\\&=\frac{C_{d}}{\lab\xi\rab^{4}}\int_{1}^{t\lab\xi\rab^{4}}\frac{d\rho}{\rho}\sim\frac{C_{t, d}
\log|\xi|}{\lab\xi\rab^{4}}, \ \ \mbox{ as }|\xi|\to\infty.
\end{split}
\eeq
This proves (\ref{spidespcashalf}). The last conclusion follows from the above proof
by using the upper bound in Mills' ratio to the inner integral $d\tau$ in
\eqref{astftimehalf}. 
The lemma is now proved.
\epf
\subsubsection{Finishing the proofs of Theorems \ref{spatmodtf} and \ref{spatmodtfhalf}}
We are now ready to finish the proof of \thmref{spatmodtf} and \thmref{spatmodtfhalf}.
We start with the case $d=1$.
\bpfs{Proof of  part (i) of \thmref{spatmodtf} and \thmref{spatmodtfhalf}}
Fix $t>0$ and assume $0<\beta \le 1/2$.  We first find the spectral density of the
gradient as follows:  we use (4.2) and Parseval's identity to get
\beq\lbl{tfgradcov}
\bsp
\mathbb{E} & \lbk \pa_{x}U_{\beta}(t,x)\pa_{y} U_{\beta}(t,y)\rbk =
\int_{0}^{t}  \int_{\mathbb{R}}  \xi^{2}\FKbetatrxxi\overline{\FKbetatryxi}\,
d\xi dr\\
& = C\int_{\mathbb{R}} \e^{\i \langle\xi, x-y\rangle}  \xi^{2}\int_{0}^{t} E^{2}_{\beta}
\bigg(-\frac{\lab\xi\rab^{2}}{2}(t-r)^{\beta}\bigg) drd\xi.
\end{split}
\eeq
This means that the spatial spectral density of $\pa_{x}U_{\beta}$ is
$\widetilde{S}_\beta (\xi) = \xi^2 S_\beta (\xi)$, where $S_\beta (\xi)$
is given in Lemma 4.1. As $|\xi| \to \infty$, the asymptotic behavior of $\widetilde{S}_\beta (\xi)$
is---upon using \lemref{tfsdasnothalf}, \lemref{tfsdashalf}, and \eqref{tfgradcov}---given by
\beq\lbl{tfgradspden}
\widetilde{S}_\beta(\xi)
\sim \bc
\df{C}{|\xi|^{2}},\ \ &\hbox{ if } 0<\beta<\frac12;\\
\df{C\log|\xi|}{|\xi|^{2}},&\hbox{ if } \beta=\frac1 2.
\ec
\eeq

We start with the H\"older assertion for the gradient in \thmref{spatmodtf} (i).
When $0<\beta<1/2$, we apply \lemref{sspdenspide}, inequality \eqref{elemineq} (a),
and the Mittag-Leffler upper bound in \eqref{elemineqML} to obtain
\beq\lbl{spderlkseq7}
\bsp
&\E\Big|\pa_{x} U_{\beta}(t,x) -\pa_{y}U_{\beta}(t,y) \Big|^2
=\int_{\R} \xi^2 \big|\e^{i x\xi} - \e^{i y\xi}\big|^2 S_{\beta}(\xi)d\xi\\
&=C \int_{\R} \xi^2 \lbk1- \cos(\xi(x-y))\rbk \int_{0}^{t} E^{2}_{\beta}
\bigg(-\frac{\lab\xi\rab^{2}}{2}(t-r)^{\beta}\bigg) dr d\xi
\\
\\&\le C \int_{\R} \lbk1- \cos(\xi(x-y))\rbk \frac{d\xi} {\xi^2} \int_{0}^{t} (t-r)^{-2\beta} dr
\\& =C  t^{1-2\beta}|x-y|,
\end{split}
\eeq
where the last equality follows from a change of variable in the integral $d\xi$ (or
the well-known formula for the variance of fractional Brownian motion). Kolmogorov's
continuity theorem gives us the spatial local $\gamma$-H\"older continuity for the
$\beta$-time-fractional SPIDEs gradient, $\pa_{x}U_{\beta}$, for $\gamma\in(0,1/2)$
and $0<\beta<1/2$. For the critical $\beta=1/2$ case in \thmref{spatmodtfhalf} (i),
we use the last statement in Lemma \ref{Lem:sd2} together with the second equality in
\eqref{spderlkseq7} and inequality \eqref{elemineq} (a) to obtain
\beq\lbl{spderlkseq77}
\bsp
\E\Big|\pa_{x} U_{1/2}(t,x) -\pa_{y}U_{1/2}(t,y) \Big|^2
&\le C \int_{\R} \lbk 1- \cos(\xi(x-y))\rbk  \frac{\log|\xi|} {\xi^2} d\xi
\\&\le C|x-y|\,\log \frac 1 {|x-y|}
\end{split}
\eeq
for all $x, y \in \R$ with $|x-y|\le 1/2$, where the last inequality follows 
from a change of variable. Hence the same H\"older assertion holds for the case 
of $\beta = 1/2$.

Turning now to the exact uniform and local spatial continuity moduli of the
$\beta$-time-fractional SPIDEs $\pa_{x}U_{\beta}$,
in \thmref{spatmodtf} and \thmref{spatmodtfhalf} (i) (a) and (b).

Combining the property of the spectral density $\widetilde{S}_\beta$ in 
\eqref{tfgradspden} and Theorems 2.1 and 2.5 in \cite{X07}, we can verify 
that the following hold: Given any constant $M>0$, there exists a finite 
constant $c>0$ (depending on $t$ and $M$)
such that for every $n\geq 1$ and
for every $ x, y_{1},..., y_{n} \in [-M,M]$,
\begin{equation}\lbl{SLNDpa1}
{\rm Var}  \left[\pa_{x}U_{\beta}(t,x) |\pa_{x} U_{\beta}(t, y_{1}),\ldots, 
\pa_{x}U_{\beta}(t, y_{n} )
\right] \geq c\min _{ 0\leq j\leq n} \varphi_\beta(\vert x -y_{j} \vert),
\end{equation}
where $y_{0}=0$, and $\varphi_\beta$ is defined on $(0, \infty)$ by
\begin{equation}\label{def:phibeta}
\varphi_\beta(r) = \left\{\begin{array}{ll}
r; \quad &\hbox{ if } 0 < \beta < 1/2,\\
r |\log r|, &\hbox{ if } \beta = 1/2.
\end{array}\right.
\end{equation}
 Also,
\begin{equation}\lbl{variogrampa1}
\E [\pa_{x} U_\beta(t,x)- \pa_{x} U_\beta(t,y)]^{2} \asymp \varphi_\beta(|x-y|);
\quad  \forall x, y\in[-M,M].
\end{equation}
Hence, $\{\pa_{x} U_\beta(t,x), x \in \R\}$ satisfies Condition (C1) and (C2) 
(or slight variants when $\beta = 1/2$). Consequently, the desired uniform 
continuity in (i) (a) of Theorems \ref{spatmodtf} and  \ref{spatmodtfhalf}
follow from Theorem 4.1 in \cite{MWX}.  The local modulus of continuity of 
$\pa_{x}U_{\beta}$ Theorems \ref{spatmodtf} and \ref{spatmodtfhalf} (i) (b) 
follow upon applying \cite[Theorem 5.1]{MWX}, completing the proof of 
Theorems \ref{spatmodtf} and \ref{spatmodtfhalf} part (i).
\epfs

We now turn to the rougher spatial regularity in two and three dimensional 
fractal cases for the SPIDEs \eqref{spide}. First, we start with the $d=3$ 
case in Theorems \ref{spatmodtf} and \ref{spatmodtfhalf} (ii),
for the cases $0<\beta<1/2$ and $\beta=1/2$, respectively.

The following lemma provides the strong local nondeterminism property and 
double-sided second moment bounds in space for the $\beta$-time-fractional 
SIPDE solution $\lbr U_{\beta}(t,x); x \in \mathbb{R}^{3}\rbr$.

\blm[Spatial SLND and double-sided bounds for time-fractional SPIDEs] 
\label{SLNDstf}
Suppose $0<\beta\le1/2$ and $d = 3$. For every fixed $t>0$, the time-fractional 
SIPDE solution $\lbr U_{\beta}(t,x); x \in \mathbb{R}^{3}\rbr$ is spatially 
strongly locally nondeterministic. Namely, for every $M>0$, there exists a 
finite constant $c>0$ (depending on $t$ and $M$) such that for every
$n\geq 1$ and for every $ x, y_{1},..., y_{n} \in [-M,M]^{3} $,
\begin{equation}\lbl{SLNDstfeq}
{\rm Var}\left[ U_{\beta}(t,x) | U_{\beta}(t, y_{1}),\ldots, U_{\beta}(t, y_{n} ) 
\right] \geq c\,\min _{ 0\leq j\leq n} \varphi_\beta( \vert x -y_{j} \vert ),
\end{equation}
where $y_{0}=0$ and the function $\varphi_\beta$ is defined in \eqref{def:phibeta}.  
Also,
\begin{equation}\lbl{dsbstf}
\E [U_\beta(t,x)-U_\beta(t,y)]^{2} \asymp \varphi_\beta(|x-y|);\quad  
\forall x, y\in[-M,M]^{3}.
\end{equation}
Moreover, as $|x-y| \to 0$, ``$\asymp$'' in \eqref{dsbstf} can be replaced 
by $\sim$ [up to a constant factor].
\elm

\bpf
The conclusions in \eqref{SLNDstfeq} and \eqref{dsbstf} follow from Lemmas 
\ref{tfsdasnothalf} and \ref{tfsdashalf} together with Theorems 2.1 and 2.5 
in \cite{X07}. Finally, the last statement of the lemma follows from 
\eqref{spidespcas}, (\ref{spidespcashalf}) and Theorem 1 of Pitman \cite{Pitman68}.
\epf

Next, we  prove the results on spatial uniform and local continuity moduli 
in \thmref{spatmodtf} and \thmref{spatmodtfhalf} (ii)-(iii).

\bpfs{Proof of \thmref{spatmodtf} and \thmref{spatmodtfhalf} part (ii)-(iii)}  
With \lemref{SLNDstf} in hand, the uniform modulus of continuity of $U_{\beta}$ 
in \thmref{spatmodtf} and \thmref{spatmodtfhalf} (ii) (a) follow from Theorem 4.1 
in \cite{MWX}; while the local modulus of continuity of $U_{\beta}$ in both 
\thmref{spatmodtf} and \thmref{spatmodtfhalf} (ii) (b)
follow upon applying \cite[Theorem 5.1]{MWX}.

To prove part (iii) of \thmref{spatmodtf} and \thmref{spatmodtfhalf}, we start by 
deriving sharp upper bounds for $\E [U_\beta(t,x)- U_\beta(t,y)]^{2} $. This is 
similar to the proof of Theorem 1.2 (iii), which can be obtained by using
the upper bounds for the spectral density function $S_\beta$ in Lemma 4.2 and 
4.3, respectively. More precisely, we can verify that 
\begin{equation}\label{Th14-3}
\E [U_\beta(t,x)- U_\beta(t,y)]^{2}  \le c\, \left\{ \begin{array}{ll}
|x- y|^{2}; \quad &\hbox{ if } 0 < \beta < 1/2,\\
|x- y|^{2} \big|\log |x-y| \big|^2, &\hbox{ if } \beta = 1/2,
\end{array}\right.
\end{equation}
for all $x, y \in \R$ with $|x-y| \le 1/2$. In the above, $c \in (0, \infty)$ 
is a constant. The rest of the proof is similar to that of part (iii) of 
Theorem 1.2 and is omitted.
 \epfs


\subsection{The time-fractional SPIDE gradient temporal H\"older and modulus of continuity}\lbl{sec:spidegradtemp}
We prove the temporal regularity of the spatial gradient $\pax U_{\beta}$ in \thmref{thm:gradtempspide}.
\bpfs{Proof of \thmref{thm:gradtempspide}}
Let $d=1$.  We start with the H\"older assertion for the gradient.  Recall $U_{\beta}(t,x)=X_{\beta}(t,x)-V_{\beta}(t,x)$, where $X_{\beta}$ is the rougher process.  Proceeding as in the proof of \thmref{thm:gradtemplks} for L-KS SPDEs,
Parseval's identity applied to the integral in $y$ gives
\begin{equation}\lbl{gradtempdiffL2spide}
\begin{split}
\mathbb{E}  \lbk \pax X_{\beta}(t,x) - \pax X_{\beta}(s,x)\rbk^{2}&= \int_{\mathbb{R}} \int_{\mathbb{R}} \lab\pax\Kbetatrposxy-\pax\Kbetasrposxy \rab^{2} dr dy \\
& = \int_{\R}  \int_{\R} \xi^{2} \lab\FKbetatrposxxi-\FKbetasrposxxi  \rab^{2} d\xi dr
\end{split}
\end{equation}
By \eqref{tempom2}, equation \eqref{gradtempdiffL2spide} becomes
\begin{equation}
\begin{split}
&\quad\E  \lbk \pax X_{\beta}(t,x) - \pax X_{\beta}(s,x)\rbk^{2}\\& =  \int_{\R}\int_{\R} \df{\xi^{2}  \lab E_{\beta}\lpa-\frac{\lab\xi\rab^{2}}{2}(t-r)^{\beta}\rpa \ind_{ \{ t > r \}} - E_{\beta}\lpa-\frac{\lab\xi\rab^{2}}{2}(s-r)^{\beta}\rpa \ind_{ \{ s > r \}}\rab^{2}}{{2\pi}}dr d\xi.
\end{split}
\end{equation}
Taking Fourier transform in $r$, assuming without loss of generality $\beta\in\{1/2^{k};k\in\N\}$, using \lemref{lm:MLFT} above, proceeding as in \eqref{FTrtf} and immediately after, and using the inequalities \eqref{elemineq} (a) and \eqref{eq:shcosineq}, and the asymptotic
\beq\lbl{eq:spidegradspdenas}
\bsp
\Delta_{\beta}(\tau)&:=(2\pi)^{-1}\int_{\mathbb{R}} \frac{\lab\tau\rab^{2(\beta-1)}\xi^{2}}{\lab\tau\rab^{2\beta}+\lab\xi\rab^{2}\lab\tau\rab^{\beta}\cos\lpa\frac{\pi\beta}{2}\rpa+\tf14\lab\xi\rab^{4}} d\xi
\\&\sim\frac{C}{\lab\tau\rab^{({4\beta^{-1}-3})/{2\beta^{-1}}}}, \mbox{ as }\lab\tau\rab\nearrow\infty,
\end{split}
\eeq
we get, for a large enough $N$, that
\begin{equation}\lbl{eq:tempspdengrad}
\begin{split}
\E & \lbk\pax X_{\beta}(t,x) -\pax X_{\beta}(s,x)\rbk^{2} = (2\pi)^{-1} \int_{\mathbb{R}} \int_{\mathbb{R}} \xi^{2}\lab \widehat{\phi} (\tau,\xi) \rab^{2} d\tau d\xi\\
& = 2(2\pi)^{-1}  \int_{\mathbb{R}} (1 - \cos ((t-s)\tau)) \int_{\mathbb{R}} \frac{\lab\tau\rab^{2(\beta-1)}\xi^{2} d\xi}{\lab\tau\rab^{2\beta}+\lab\xi\rab^{2}\lab\tau\rab^{\beta}\cos\lpa\frac{\pi\beta}{2}\rpa+\tf14\lab\xi\rab^{4}} d\tau
\\&\le C|t-s|^{2\alpha}\lbk\int_{0}^{N}{(1 - \cos (\tau))}\Delta_{\beta}(\tau)d\tau+\int_{N}^{\infty}\lab\tau\rab^{2\alpha-({4\beta^{-1}-3})/{2\beta^{-1}}}d\tau\rbk
\\&\le C|t-s|^{2\alpha},\ 0<\alpha<(2\beta^{-1}-3)/4\beta^{-1}.
\end{split}
\end{equation}
It follows that $\pax X_{\beta}(\cdot,x)$ is $\gamma$-H\"older continuous in time, with $\gamma\in\lpa0,(2\beta^{-1}-3)/4\beta^{-1}\rpa$.  This, together with the gradient decomposition
\beq\lbl{eq:graddecspide}
\pax U_{\beta}(t,x)=\pax X_{\beta}(t,x)-\pax V_{\beta}(t,x),
\eeq
and the fact that $\pax V_{\beta}$ is temporally smooth\footnote{\lbl{fn:gradVbeta}As in the L-KS SPDE case, it can be shown that the smoothness assertions in \thmref{auxdecompspiden} (ii) and (iii) (with $(H=2\beta^{-1}-3)/4\beta^{-1}$) hold for $\pax V$.  Since the proof follows the same steps as the one for \thmref{auxdecompspden} (ii) and (iii) with straightforward modifications, we leave it to the interested reader.} establish the H\"older regularity assertion for $\pax U_{\beta}$ in \thmref{thm:gradtemplks}.

Turning now to the uniform and local spatial continuity moduli results for the time-fractional SPIDE 
gradient, $\pa_{x}U_{\beta}$, in \thmref{thm:gradtempspide}.   Equation \eqref{eq:tempspdengrad} 
means that $\pax X_{\beta}$ has stationary increments and the spatial spectral density of $\pa_{x}
X_{\beta}$ and its asymptotic behavior are given by \eqref{eq:spidegradspdenas}.

Equation \eqref{eq:spidegradspdenas} and Theorem 2.1 in \cite{X07} imply that, for every fixed $x\in\R$, 
the gradient $\{\pa_{x}X_{\beta}(t,x), t\ge0\}$ is temporally strongly locally nondeterministic. Namely, 
for any $T>0$, there is a positive constant $c$ such that for all $t\in(0,T]$ and all $0 < r \le 1\wedge|t|$
\beq\lbl{SLNDtlkseq}
\var \lpa X_{\beta}(t,x) | X_{\beta}(s,x);s\in[0,T],|t-s|\ge r\rpa \ge c r^{\frac{2\beta^{-1}-3}{2\beta^{-1}}}
\eeq
Also, \eqref{lksgradspdensp} and Theorem 2.5 in \cite{X07} imply the double sided second moment bounds
\beq\lbl{dstemplks}
\E [X_{\beta}(t,x)-X_{\beta}(s,x)]^{2} \asymp|t-s|^{\frac{2\beta^{-1}-3}{2\beta^{-1}}};\ \forall s,t\in[0,T].
\eeq
Thus, the uniform modulus of continuity of $\pax X_{\beta}$
 \beq\lbl{eq:tempumspideauxgrad}
\P\lbk\lim_{\delta\searrow0} \sup_{\substack{|t-s|<\delta\\t,s\in\tint} }\frac{\lab\pa_{x} X_{\beta}(t,x)-\pa_{x}X_{\beta}(s,x)\rab}{|t-s|^{(2\beta^{-1}-3)/4\beta^{-1}}\sqrt{\log\lbk1/|t-s|\rbk}}=k\rbk=1,
\eeq
for every compact interval $\tint\subset\Rp$ and for some constant $k>0$, follows from Theorem 4.1 in \cite{MWX}.  The local modulus of continuity of $\pax X$
\beq\lbl{eq:templmspideauxgrad}
\P\lbk\lim_{\delta\searrow0} \sup_{\substack{|t-s|<\delta} }\frac{\lab\pa_{x} X_{\beta}(t,x)-\pa_{x}X_{\beta}(s,x)\rab}{\delta^{(2\beta^{-1}-3)/4\beta^{-1}}\sqrt{\log\log\lbk1/\delta\rbk}}=k\rbk=1,
\eeq
follows upon applying \cite[Theorem 5.1]{MWX}.  The corresponding continuity moduli assertions for the gradient $\pax U_{\beta}$ in \thmref{thm:gradtempspide} follow from those of the auxiliary process $\pax X_{\beta}$ (\eqref{eq:tempumspideauxgrad} and \eqref{eq:templmspideauxgrad}), the decomposition \eqref{eq:graddecspide}, and the smoothness of $\pax V_{\beta}$ (see footnote \ref{fn:gradVbeta}).
\epfs

\section{From linear to time-fractional Allen-Cahn and Swift-Hohenberg equations via measure change}
We quickly remark in this section that, at their core, the space-time change of measure theorems in \cite{Acom,Acom1,Acom2} are ``noise'' results that are independent of both the type and order of the SPDE under consideration.  This makes them conveniently adaptable to
different SPDEs settings.  As was done in \cite{Alksspde} for L-KS SPDEs, we can extend the  results in \cite{Acom1,Acom2} to our
$\beta$-time-fractional SPIDEs \eqref{spide}.  The almost sure $L^{2}$ condition in \cite{Acom1,Acom2,Alksspde}, which is much
weaker than the usual Novikov condition typically found in change-of-measure results, allows us to state an equivalence in law---and thus in all almost sure regularity results---between both \eqref{nlks} and \eqref{spide} and their nonlinear versions the Swift-Hohenberg SPDEs
\begin{equation} \label{SH}
 \begin{cases} \displaystyle\frac{\partial U}{\partial t}=
-\tfrac\vep8\lpa\lap+2\vth\rpa^{2}U+b(U)+\frac{\partial^{d+1} W}{\partial t\partial x}, & (t,x)\in(0,T]\times\S;
\cr U(0,x)=\unx, & x\in\S,
\end{cases}
\end{equation}
and the $\beta$-time-fractional Allen-Cahn SPIDE
\beq\lbl{ACspide}
\begin{cases} \displaystyle
{}^{\mbox{\tiny C}}\pa_{t}^{\beta}U_{\beta}=\tfrac12\D U_{\beta}+I_{t}^{1-\beta}\lbk b(U)+ \frac{\partial^{d+1} W}
{\partial t\partial x}\rbk,& (t,x)\in(0,T]\times\S;
\cr U_{\beta}(0,x)=\unx, & x\in\S,
\end{cases}
\eeq	
respectively, where $T>0$ is fixed and arbitrary, and where
\beq\lbl{gSHcom}
b(u)=\sum_{k=0}^{2p-1}c_{k}u^{k},\ \S=\prod_{i=1}^{d}[0,L_{i}], \mbox{ and with $p\in\N$,   $c_{2p-1}<0$, and $d=1,2,3$.}
\eeq
Let $\T:=[0,T]$.  We supplement \eqref{SH} and \eqref{gSHcom} with suitable boundary conditions\footnote{E.g., boundary conditions of Neumann type $\pa U/\pa n=\pa\D U/\pa n=0$ or Dirichlet type conditions  $U=\D U=0$ on $\pa\S$ and $d=1,2,3$.}, the nature of which is irrelevant to our next change of measure result.  For concreteness, we assume Dirichlet boundary conditions throughout this section.  We also modify the kernels $\KKSepthtx$ and $\Kbetatx$ in the mild kernel formulations \eqref{ibtbapsol} and \eqref{isltbmsie} (with $a\equiv1$) to account for the boundary conditions\footnote{E.g., in the Neumann (Dirichlet) case, the propagator ${ \e^{-|x-y|^2/2\i s}}/{{\lpa2\pi \i s \rpa}^{d/2}}$ in the definition of the $(\vep,\vth)$ L-KS kernel $\KKSepthtx$ \eqref{vepvthLKS} is replaced with the propagator with reflection (absorption) at $\pa\S$, respectively. Similar comments apply to the outside $d$-dimensional BM density $\psx$ in the definition of $\Kbetatx$ \eqref{eq:fundsolbeta}.}, and we replace $\Rd$ with $\S$.  The linear-nonlinear equivalence result is now stated.  For completeness, we restate the Swift-Hohenberg conclusions from \cite{Alksspde}.

\bfr
\bthm[Swift-Hohenberg and time-fractional Allen-Cahn law equivalence to their linear counterparts]\lbl{SHcom}
Fix $T>0$.  Let $\S=\prod_{i=1}^{d}[0,L_{i}]$, $d=1,2,3,$ and assume that $u_0$ satisfies \eqref{init} with $\Rp\times\Rd$ replaced by $\T\times\S$. The generalized Swift-Hohenberg SPDE, \eqref{SH} and \eqref{gSHcom}, admits uniqueness in law and is law equivalent to the $b\equiv0$ version of \eqref{SH} on $\sB\lpa\C(\T\times\S;\R)\rpa;$ consequently,  it has the same H\"older continuity and modulus of continuity regularity as the linear L-KS SPDE on $\T\times\S$.  The same uniqueness assertion holds for the $\beta$-time-fractional Allen-Cahn SPIDE \eqref{ACspide} and \eqref{gSHcom}.  Also, the law---and hence the H\"older continuity and the continuity modulus regularity---equivalence hold between the $\beta$-time-fractional Allen-Cahn SPIDE, \eqref{ACspide} and \eqref{gSHcom}, and its zero-drift ($b\equiv0$) version.  The uniqueness in law and the law equivalence (hence regularity equivalence) conclusions above all hold if $b(u)=\sum_{k=0}^{l}c_{k}u^{k}$, for $l\in\N\cup\{0\}$ and $c_{k}\in\R$. 
\ethm
\efr
The proof of the uniqueness and law equivalence assertions for L-KS SPDEs (in both the Allen-Cahn nonlinearity $b$ \eqref{gSHcom} and the general polynomial $b$ cases) was given in \cite[Theorem 1.3 and Corollary 1.1]{Alksspde}.  The proof is exactly the same  for the linear-to-nonlinear time-fractional SPIDEs case, and we omit it\footnote{The pathwise uniqueness, and hence uniqueness in law, of solutions trivially follows in the linear $b\equiv0$ case from the mild formulation.}.

\section{Concluding remarks}
\subsection{Time-fractional SPIDEs are \emph{different} from time fractional SPDEs and their equivalent high-order memoryful SPDEs}  Here, we make a brief but important distinction that was emphasized in \cite{Atfhosie,Abtbmsie}, and that the astute reader will note.  For completeness, we incorporate the SPIDE coinage of our present paper here in making our point, which we now state and discuss.  The time-fractional SPIDEs \eqref{spidea} are \emph{not} equivalent to (their rigorous SIEs formulation \eqref{isltbmsie} are \emph{not} the mild form of) their \emph{rougher} and \emph{fundamentally different} relatives: (1) the time-fractional SPDEs
\beq\lbl{fracbetaspde}\bc
\ds\p^{\beta}_{t}{U_{\beta}}=\tfrac12\Delta U_{\beta}+a(U_{\beta})\frac{\partial^{d+1} W}{\partial t\partial x}& (t,x)\in\Rpop\times\Rd\\
U_{\beta}(0,x)=u_0(x),&x\in\Rd,
\ec
\eeq
and (2) the $2\beta^{-1}=2\nu$ order, $\nu\in\lbr2^k;k\in\N\rbr$, memoryful SPDEs
\beq\lbl{2nuspde}\bc
\ds\df{\p U_{\beta}}{\p t}=\sum_{\kappa=1}^{{\nu}-1}\frac{C_{\beta,\kappa}\Delta^{\kappa}\unx}{2^{\kappa}t^{1-\kappa/\nu}}{}+\frac{\Delta^{{\nu}}U_{\beta}}{2\nu}+a(U_{\beta})\frac{\partial^{d+1} W}{\partial t\partial x}& (t,x)\in\Rpop\times\Rd\\
U_{\beta}(0,x)=u_0(x),&x\in\Rd
\ec
\eeq
where $C_{\beta,\kappa}=\frac{\E\lpa\Lambda_{\beta}(1)\rpa^{\kappa}}{\kappa!}$;  the process $\ILb$ is the $\beta$-inverse-stable-L\'evy motion, as in \secref{rigdisc} above, which arises in the work of Meerschaert et al.~\cite{MBS,MS} as scaling limits of continuous time random walks and which is reviewed, along with its link to $k$-iterated Brownian-time Brownian motion, in \cite{Atfhosie}.

 In \cite{Abtbmsie}, Allouba showed that in the case $\beta=1/2$ (the Brownian-time Brownian motion case), the $\beta$-time-fractional SIE in \eqref{isltbmsie} is \emph{not} the mild formulation of the $\beta=1/2$ of \eqref{2nuspde}; but rather it is an integral formulation of what he called parametrized BTBM SPDE, evaluated at the diagonals (see \cite{Abtbmsie} pp.~428--431, Lemma 1.2, and footnote 3 p.~416 for the details).  Similarly, as stressed in \cite{Atfhosie}, for general $\beta\in\{1/2^{k};k\in\N\}$, the $\beta$-time-fractional SIE in \eqref{isltbmsie} is \emph{not} the mild formulation of \eqref{2nuspde}.  This is contrary to what was erroneously stated in \cite{MN}.   Thus, unfortunately, Theorem 4 of \cite{MN} claiming the equivalence between the time-fractional SPIDE\footnote{In \cite{MN} time-fractional SPIDEs are called time-fractional SPDEs, which is less precise since the name ignores the crucial smoothing effect of the time-fractional integral $I^{1-\beta}$ in the formal formulation \eqref{spidea}.  We reserve the name time-fractional SPDEs for \eqref{fracbetaspde}.} \eqref{spidea}---corresponding to the rigorous SIE form \eqref{isltbmsie}---and the memoryful high order SPDEs \eqref{2nuspde} is incorrect.  In  fact, in \cite{Atfhosie} it was repeatedly emphasized that the $\beta$-time-fractional SIE in \eqref{isltbmsie} (and hence its corresponding SPIDE \eqref{spidea}) is a \emph{different and smoother} stochastic version of the deterministic PDEs---obtained by setting $a\equiv0$ in either \eqref{fracbetaspde} or \eqref{2nuspde}---than the \emph{rougher} equivalent SPDEs \eqref{fracbetaspde} and \eqref{2nuspde} (see the discussions in \cite{Atfhosie} right before equations (1.4) and (1.8) and the discussion following equation (1.15), including footnote 15).   

The formal SPDEs in \eqref{fracbetaspde} and \eqref{2nuspde} require rigorous formulations quite different from \eqref{isltbmsie}.  This is handled, and the equivalence between \eqref{fracbetaspde} and \eqref{2nuspde} for suitably regular initial data $\un$, is shown in an upcoming separate article.  Of course, even formally, it is obvious that the time-fractional SPIDEs \eqref{spidea} are \emph{not} equivalent to the time-fractional SPDEs \eqref{fracbetaspde}, which lack the fractional integral and its smoothing effect.

\subsection{Other remarks}\lbl{fnlrem}\hspace{1mm}
Further properties on the local times and fractal behavior of the solution process for both L-KS SPDEs and time-fractional SPIDEs $\{U(t, x), x \in \R^d\}$, when $t >0$ is fixed,
can now be derived from \cite{Xiao96,X07,X09}. It is also possible to investigate sample path properties of the Gaussian random field  $\{U(t, x), t \ge 0, x \in \R^d\}$ in both time and space variables. We will carry out this in subsequent work.
\appendix

\section{Glossary of frequently used acronyms and notations}\lbl{glossary}
\begin{enumerate}\renewcommand{\labelenumi}{\Roman{enumi}.}
\item {\textbf{Acronyms}}\vspace{2mm}
\bit
\item BM: Brownian motion.
\item bifBM: bifractional BM.
\item BTBM: Brownian-time Brownian motion.
\item SIE: Stochastic integral equation.
\item SLND: Strong local nondeterminism.
 \item SPIDE: Stochastic partial integro-differential equation.
\item KS: Kuramoto-Sivashinsky.
\eit
\vspace{2.5mm}
\item {\textbf{Notations}}\vspace{2mm}
\bit
\item $B^{(H,K)}$: bifractional BM with indices $H$ and $K$.
\item $\N$: The usual set of natural numbers $\lbr1,2,3,\ldots\rbr$.
\item $\T$: The time interval $[0,T]$ for some arbitrary ficed $T>0$.
\item $\ptsz$: The  density of a $1$-dimensional BM, starting at $0$.
\item $\KKSepthtx$:  The generalized $(\vep,\vth)$ L-KS kernel.
\item $\KBtx$: The kernel or density of a $d$-dimensional Brownian-time Brownian motion.
\item $\E$: The expectation operator.
\item $E_{\beta}$: The Mittag-Leffler function. 
\item $\C^{k,\gamma}(\R,\R)$: The set of $k$-continuously differentiable functions on $\R$ whose $k$-th derivative is locally H\"older continuous, with H\"older exponent $\gamma$.
\item $\H^{\gamma_{*}^{-}}(\Rp;\R)$:  The space of locally H\"older continuous functions $f:\Rp\to\R$ whose H\"older exponent $\gamma\in(0,\gamma_{*})$.
\item $\pa^{n}_{x_{i}}f(x_{1},\ldots,x_{N})=\pa^{n} f/\pa x^{n}_{i}$, $i=1,\ldots,N$ and $n\in\N$.
\item $f(x)\asymp g(x)$ on $\S$ means $c_{l}g(x)\le f(x)\le c_{u}g(x)$ for some constants $c_{l},c_{u}$ for every $x\in\S$.
\vspace{0.5 mm}
\eit
\end{enumerate}


\begin{thebibliography}{100}\lbl{ref}
\bibitem{Alksspde} H.~Allouba,  \emph{L-Kuramoto-Sivashinsky SPDEs in one-to-three dimensions: L-KS kernel, 
sharp Hšlder regularity, and Swift-Hohenberg law equivalence.}  J. Differential Equations {\bf 259} (2015), no. 11, 6851--6884.
\bibitem{Atfhosie} H.~Allouba, \emph{Time-fractional and memoryful $\Delta^{2^{k}}$ SIEs on $\Rp\times\Rd$: 
how far can we push white noise?}  Illinois J. Math. {\bf 57} (2013), no. 4, 919--963.
\bibitem {Abtbmsie}  H.~Allouba,  \emph{Brownian-time Brownian motion SIEs on $\Rp\times\Rd$: ultra regular direct and lattice-limits solutions and fourth order SPDEs links.}  Discrete Contin. Dyn. Syst. {\bf 33} (2013), no. 2, 413--463.  
\href{http://www.ams.org/mathscinet-getitem?mr=MR2975119}{MR2975119}
\bibitem{Abtbs} H.~Allouba,
\emph{From Brownian-time Brownian sheet to a fourth order and a Kuramoto-Sivashinsky-variant interacting PDEs systems}.
\newblock Stoch. Anal. Appl. \textbf{29} (2011), 933--950. \href{http://www.ams.org/mathscinet-getitem?mr=MR2847330}{MR2847330}
\bibitem{Abtpspde}  H.~Allouba,   \textit{A Brownian-time excursion into fourth-order PDEs,
linearized Kuramoto-Sivashinsky, and BTP-SPDEs on ${\mathbb R}\sb +\times{\mathbb R}\sp d$}.
 {Stoch. Dyn.} {\bf 6} (2006), no. 4, 521--534.
\href{http://www.ams.org/mathscinet-getitem?mr=MR2285514}{MR2285514}
\bibitem{Aks}  H.~Allouba,  \textit{A linearized Kuramoto-Sivashinsky PDE via an imaginary-Brownian-time-Brownian-angle process.} 
{C. R. Math. Acad. Sci. Paris} {\bf 336} (2003), no. 4, 309--314. \href{http://www.ams.org/mathscinet-getitem?mr=MR1976309}{MR1976309}
\bibitem{Abtp2} H.~Allouba,  \textit{Brownian-time processes: the PDE connection II and the corresponding Feynman-Kac formula.}
 {Trans. Amer. Math. Soc.}  {\bf 354} (2002), no. 11, 4627--4637 (electronic).
\href{http://www.ams.org/mathscinet-getitem?mr=MR1926892}{MR1926892}
\bibitem{Acom2} H.~Allouba, \textit{SPDEs law equivalence and the compact support property: applications to the Allen-Cahn SPDE.} {C. R. Acad. Sci. Paris Sér. I Math.} {\bf 331} (2000), no. 3, 245--250.
\href{http://www.ams.org/mathscinet-getitem?mr=MR1781835}{MR1781835}
\bibitem{Acom1}  H.~Allouba, \textit{Uniqueness in law for the Allen-Cahn SPDE via change of measure.} {C. R. Acad. Sci. Paris Sér. I Math.} {\bf 330} (2000), no. 5, 371--376.
\href{http://www.ams.org/mathscinet-getitem?mr=MR1751673}{MR1751673}
\bibitem{Acom}   H.~Allouba, \textit{Different types of SPDEs in the eyes of Girsanov's theorem.}  {Stochastic Anal. Appl.} {\bf 16} (1998), no. 5, 787--810.  \href{http://www.ams.org/mathscinet-getitem?mr=MR1643116}{MR1643116}
\bibitem{AN} H.~Allouba and E.~Nane,
\emph{Interacting time-fractional and $\Delta^\nu$ PDEs systems via Brownian-time and Inverse-stable-L\'evy-time Brownian sheets}. 
{Stoch. Dyn.~\textbf{13} (2013),  no.1  (31 pages).} \href{http://www.ams.org/mathscinet-getitem?mr=MR3007250}{MR3007250}
\bibitem{Abtp1}  H.~Allouba and W.~Zheng, \textit{Brownian-time processes: the PDE connection and the half-derivative generator.}  {Ann. Probab.}  {\bf 29} (2001), no. 4, 1780--1795. \href{http://www.ams.org/mathscinet-getitem?mr=MR1880242}{MR1880242}

\bibitem{Bert}  J. Bertoin, L\'evy Processes. Cambridge University Press, Cambridge, 1996.
\bibitem{LeC} L.~Chen, \emph{Nonlinear stochastic time-fractional diffusion equations on $\R$: moments, H\"older regularity and intermittency.}  \href{http://arxiv.org/abs/1410.1911}{arXiv1410.1911}.
\bibitem{LeCDal} L. Chen and R. C. Dalang, \emph{H\"older-continuity for the nonlinear stochastic heat equation with rough initial conditions.} Stoch. PDE: Anal. Comp. {\bf 2} (2014), 316--352.
\bibitem{ZCHKmPKm} Z.-Q. Chen, K.-H. Kim and P. Kim, \emph{Fractional time stochastic partial differential equations.} Stoch. Process. Appl. {\bf 125} (2015), no. 4, 1470--1499.
\bibitem{DOOToa14} M. D'Ovidio, E. Orsingher and B. Toaldo,  \emph{Time-changed processes governed by space-time fractional telegraph equations.} Stoch. Anal. Appl. 32 (2014), no. 6, 1009--1045.
\bibitem{DuanWei14}  J. Duan and W. Wei, Effective Dynamics of Stochastic Partial Differential Equations. 
Elsevier Insights, Elsevier, Amsterdam, ISBN 978-0-12-800882-9, 2014, xii+270 pp. MR3289240.
\bibitem{GOP15} R.~Garra, E.~Orsingher, and F.~Polito, \emph{Fractional diffusions with time-varying coefficients.} 
J. Math. Phys. 56 (2015), no. 9, 093301, 17 pp. 
\bibitem{HauMathSax11} H.J. Haubold, A.M. Mathai, R.K. Saxena, \textit{Review Article: Mittag-Leffler functions and their applications.} 
J. Appl. Math. 2011 (2011) 51.
\bibitem{HoudVilla03} C.~Houdr\'e and J.~Villa, \emph{An example of infinite dimensional quasi-helix.} 
Stochastic models (Mexico City, 2002), 195--201,
Contemp. Math., 336, Amer. Math. Soc., Providence, RI, 2003.
\bibitem{LN09} P. Lei and D. Nualart,
{\em A decomposition of the bifractional Brownian motion and some applications.}
        Statist. Probab. Lett. {\bf 79} (2009), 619--624.
\bibitem{}  V.~Keyantuo and C.~Lizama, \emph{On a connection between powers of operators 
and fractional Cauchy problems.} J. Evol. Equ. 12 (2012), no. 2, 245--265.
\bibitem{LX12}
N. Luan and Y. Xiao, \emph{Spectral conditions for strong local nondeterminism and exact Hausdorff measure 
of ranges of Gaussian random fields.} J. Fourier Anal. Appl. {\bf 18} (2012), no. 1, 118--145.
\bibitem{LX10} N. Luan and Y. Xiao, {\em Chung's law of the iterated logarithm
for anisotropic Gaussian random fields.}  Statist. Probab. Lett. {\bf 80} (2010), 1886--1895.

\bibitem{LunSin} A. Lunardi and E. Sinestrari, {\it An inverse problem in the theory of materials with memory.} 
Nonlin. Anal. Theory Meth. Appl. {\bf 12} (1988), 1317--1355.
\bibitem{Main1}  F. Mainardi, Fractional Calculus and Waves in Linear Viscoelasticity. 
Imperial College Press, London, 2010.
\bibitem{MainLuch} F. Mainardi, Y. Luchko and G. Pagnini,  {\it The fundamental solution of the space-time 
fractional diffusion equation.}  Fract. Calc. Appl. Anal., {\bf 4} (2001), (2): 153--192.
\bibitem{MarRos}  
M.B. Marcus and J. Rosen, Markov Processes, Gaussian Processes, and Local Times. 
Cambridge University Press, Cambridge, 2006.
\bibitem{MathHau07}   
A.M. Mathai, H.J. Haubold, Special Functions for Applied Scientists, Springer, 2007.
\bibitem{MBS} M.M. Meerschaert, D.A.  Benson, H.P. Scheffler and B. Baeumer, 
{\it Stochastic solution of space-time fractional diffusion equations.}  Phys. Rev. E 
{\bf 65} (2002), 1103--1106.\href{http://www.ams.org/mathscinet-getitem?mr=MR1917983}{MR1917983}
\bibitem{MeerSik}  M. M. Meerschaert and A. Sikorskii, Stochastic Models for Fractional Calculus. 
De Gruyter Studies in Mathematics 43, De Gruyter, Berlin, 2012, ISBN 978-3-11-025869-1.  
\href{http://www.ams.org/mathscinet-getitem?mr=MR2884383}{MR2884383}
\bibitem{MS}  M.M. Meerschaert and H.P. Scheffler, {\it Limit theorems for continuous time random walks 
with infinite mean waiting times.} J. Applied Probab. {\bf 41} (2004), 623--638.
\bibitem{MWX} M.M. Meerschaert, W. Wang and Y. Xiao, {\it Fernique-type inequalities 
and moduli of continuity for anisotropic Gaussian random fields.} Trans. Amer. Math. Soc. {\bf 365} (2013), 1081--1107.
\bibitem{MN}  J.B. Mijena and N. Erkan, \emph{Space-time fractional stochastic partial differential equations.} Stoch. Process. Appl. {\bf 125} (2015), 3301--3326.
\bibitem{MuTr02}  C. Mueller and R. Tribe, \emph{Hitting probabilities of a random string.} 
Electronic J. Probab. 7, (2002), Paper No. 10, 29 pp.
\bibitem{MW12} C. Mueller and Z. Wu, {\em Erratum: A connection between the stochastic
heat equation and fractional Brownian motion and a simple proof of a result of Talagrand.}
Electron. Commun. Probab. {\bf 17} (2012)(8), 10 pp.
\bibitem{Pitman68} E. J. G. Pitman, {\it On the behavior of the characteristic function of a probability 
distribution in the neighbourhood of the  origin.}  J. Australian Math. Soc. Series A {\bf 8} (1968), 422--443.
\bibitem{RT}F. Russo and C.A. Tudor,  {\em On bifractional Brownian motion.} Stoch. Proc. Appl.  {\bf 116} (2006), 830--856.
\bibitem{TS}  T.~Simon, \emph{Comparing Fr\'echet and positive stable laws}.  Electron. J. Probab. {\bf 19} (2014), no. 16, 1--25.
\bibitem{Tal95} M.~Talagrand,  \emph{Hausdorff measure of trajectories of multiparameter fractional Brownian motion.} Ann. Probab. {\bf 23} (1995), no. 2, 767--775.
\bibitem{T}  R.~Temam, Infinite-dimensional Dynamical Systems in Mechanics and Physics. Second edition. Applied Mathematical Sciences, 68, Springer-Verlag, New York, 1997. xxii+648 pp. ISBN: 0-387-94866-X.
\bibitem{TuX07}  
C.~Tudor and Y.~Xiao, \emph{Sample path properties of bifractional Brownian motion.}  Bernoulli
{\bf 13} (2007), 1023--1052.
\bibitem{TuX15}  C.~Tudor and Y.~Xiao, \emph{Sample paths of the solution to the fractional-colored stochastic heat equation}. 
Stoch. Dyn., to appear. \href{http://arxiv.org/abs/1501.06828}{arXiv1501.06828}.
\bibitem{WX06}  D. Wu and Y. Xiao, \textit{Fractal properties of random string processes. } IMS Lecture Notes-
Monograph Series High Dimensional Probability. 51, pp.128--147, Institute of
Mathematical Statistics, Beachwood, Ohio, U.S.A. 2006.
\bibitem{UmSa}  S. Umarov and E. Saydamatov, \emph{A fractional analog of the Duhamel principle.} 
Fract. Calc. Appl. Anal., {\bf 9} (2006) (1), 57--70.
\bibitem{W} J.B.~Walsh,  An Introduction to Stochastic Partial Differential Equations. \'Ecole d'\'et\'e de Probabilit\'es de Saint-Flour XIV. Lecture Notes in Math. 1180. Springer, New York.  1986.
\href{http://www.ams.org/mathscinet-getitem?mr=MR0876085}{MR0876085}
\bibitem{Xiao96}
Y. Xiao, {\em H\"older conditions for the local times and the
Hausdorff measure of the level sets of Gaussian random fields.}
 Probab. Th. Rel. Fields  {\bf 109} (1996), 129--157.
\bibitem{X07} Y. Xiao, \textit{Strong local nondeterminism and the sample path properties of Gaussian random fields. }
In: Asymptotic Theory in Probability and Statistics with Appli- cations (Tze Leung Lai, Qiman Shao, Lianfen Qian, editors), pp. 136--176, Higher Education Press, Beijing, 2007.
\bibitem{X09} Y. Xiao, \textit{Sample path properties of anisotropic Gaussian random fields.} In: A Mini- course on Stochastic Partial Differential Equations, D. Khoshnevisan, F. Rassoul-Agha, editors, Lecture Notes in Math. 1962, pp 145Ð212, Springer, New York, 2009.
\bibitem{XX11}
Y. Xue and Y. Xiao,  {\em Fractal and smoothness properties of space-time Gaussian models.}
Frontiers Math. China {\bf 6} (2011), 1217--1246.

\end{thebibliography}
\end{document}